\documentclass[UTF-8,reqno]{amsart}
\usepackage{enumerate, bbm, soul, ulem}
\usepackage{amssymb,url,color, booktabs}
\usepackage{mathrsfs}

\usepackage{color}
\usepackage[colorlinks=true]{hyperref}
\hypersetup{
    linkcolor=blue,          
    citecolor=red,        
    filecolor=blue,      
    urlcolor=cyan
}

\usepackage{color}
\definecolor{MyDarkBlue}{cmyk}{0.8,0.3,0.8,0.4}
\definecolor{yellow}{rgb}{0.99,0.99,0.70}
\definecolor{white}{rgb}{1.0,1.0,1.0}
\definecolor{black}{rgb}{0.00,0.00,0.00}

\newcommand{\red}{\color{red}}

\numberwithin{equation}{section}

\newcommand{\be}{\begin{eqnarray}}
\newcommand{\ee}{\end{eqnarray}}
\newcommand{\ce}{\begin{eqnarray*}}
\newcommand{\de}{\end{eqnarray*}}
\newtheorem{theorem}{Theorem}[section]
\newtheorem{lemma}[theorem]{Lemma}
\newtheorem{remark}[theorem]{Remark}
\newtheorem{definition}[theorem]{Definition}
\newtheorem{proposition}[theorem]{Proposition}
\newtheorem{Examples}[theorem]{Example}
\newtheorem{corollary}[theorem]{Corollary}

\def\nor{|\mspace{-3mu}|\mspace{-3mu}|}

\def\eps{\varepsilon}

\def\e{\mathrm{e}}

\def\p{\partial}

\def\[{{\Big[}}
\def\]{{\Big]}}
\def\<{{\langle}}
\def\>{{\rangle}}
\def\({{\Big(}}
\def\){{\Big)}}

\def\bx{{\mathbf{x}}}

\def\osc{{\rm osc}}

\def\dif{{\mathord{{\rm d}}}}

\def\min{{\mathord{{\rm min}}}}

\def\bbp{{\boldsymbol{p}}}
\def\bbr{{\boldsymbol{r}}}
\def\bbq{{\boldsymbol{q}}}

\def\bb2{{\boldsymbol{2}}}
\def\no{\nonumber}
\def\={&\!\!=\!\!&}

\def\bC{{\mathbf C}}

\def\cA{{\mathcal A}}
\def\cB{{\mathcal B}}
\def\cC{{\mathcal C}}
\def\cD{{\mathcal D}}
\def\cE{{\mathcal E}}

\def\cG{{\mathcal G}}

\def\cI{{\mathcal I}}
\def\cJ{{\mathcal J}}

\def\cM{{\mathcal M}}

\def\cP{{\mathcal P}}

\def\cR{{\mathcal R}}

\def\mB{{\mathbb B}}

\def\mD{{\mathbb D}}
\def\mE{{\mathbb E}}

\def\mI{{\mathbb I}}

\def\mL{{\mathbb L}}

\def\mN{{\mathbb N}}

\def\mP{{\mathbb P}}
\def\mQ{{\mathbb Q}}
\def\mR{{\mathbb R}}
\def\mS{{\mathbb S}}

\def\bP{{\mathbf P}}

\def\bE{{\mathbf E}}
\def\1{{\mathbf{1}}}

\def\b1{{\mathbbm 1}}

\def\sA{{\mathscr A}}

\def\sF{{\mathscr F}}

\def\sI{{\mathscr I}}

\def\sL{{\mathscr L}}
\def\sM{{\mathscr M}}

\def\sV{{\mathscr V}}

\def\geq{\geqslant}
\def\leq{\leqslant}

\def\div{\mathord{{\rm div}}}

\def\fint{{-\!\!\!\!\!\!\int\!\!}}

\def\eps{\varepsilon}

\def\e{\mathrm{e}}

\def\p{\partial}

\def\[{{\Big[}}
\def\]{{\Big]}}
\def\<{{\langle}}
\def\>{{\rangle}}
\def\({{\Big(}}
\def\){{\Big)}}

\def\bx{{\mathbf{x}}}

\def\osc{{\rm osc}}

\def\dif{{\mathord{{\rm d}}}}

\def\min{{\mathord{{\rm min}}}}

\def\no{\nonumber}
\def\={&\!\!=\!\!&}
\def\bt{\begin{theorem}}
\def\et{\end{theorem}}
\def\bl{\begin{lemma}}
\def\el{\end{lemma}}
\def\br{\begin{remark}}
\def\er{\end{remark}}
\def\bx{\begin{Examples}}
\def\ex{\end{Examples}}
\def\bd{\begin{definition}}
\def\ed{\end{definition}}
\def\bp{\begin{proposition}}
\def\ep{\end{proposition}}
\def\bc{\begin{corollary}}
\def\ec{\end{corollary}}

\def\geq{\geqslant}
\def\leq{\leqslant}

\def\div{\mathord{{\rm div}}}
\def\wt{\widetilde}

\def\bH{{\mathbf H}}

\def\bP{{\mathbf P}}

 \def\R{\mathbb R}
 \def\R{\mathbb R}

\def\<{\langle} \def\>{\rangle}

\allowdisplaybreaks

\def\wt{\widetilde}

\def\red{\color{red} }

\begin{document}


\title[Nonlocal nonlinear-diffusion processes with supercritical drifts]
{Harnack inequalities for nonlocal operators with supercritical drifts and their applications}

\author{Zhen-Qing Chen and Xicheng Zhang}

\address{Zhen-Qing Chen:
Department of Mathematics, University of Washington, Seattle, WA 98195, USA\\
Email: zqchen@uw.edu
 }

\address{Xicheng Zhang:
School of Mathematics and Statistics, Beijing Institute of Technology, Beijing 100081, China\\
Faculty of Computational Mathematics and Cybernetics, Shenzhen MSU-BIT University, 518172 Shenzhen, China\\
		Email: xczhang.math@bit.edu.cn
}

\thanks{
This work is supported by National Key R\&D program of China (No. 2023YFA1010103) and NNSFC grant of China (No. 12131019)  and the DFG through the CRC 1283 ``Taming uncertainty and profiting from randomness and low regularity in analysis, stochastics and their applications''. }

\begin{abstract}

In this paper, we investigate Harnack estimates for weak solutions to the following nonlocal equation:
$$
\partial_t u = \Delta^{\alpha/2} u + b \cdot \nabla u + f,
$$
where $\Delta^{\alpha/2}$ denotes the fractional Laplacian, $b$ is a divergence-free vector field in a critical or supercritical regularity regime, and $f$ is a distribution
 in a fractional Sobolev space with negative indices.
As applications of the analytical results obtained in this paper, we establish the well-posedness of critical stochastic quasi-geostrophic equations driven by additive Brownian noise, prove the existence of weak solutions to the two-dimensional fractional Navier--Stokes equations with measure-valued initial vorticity, and demonstrate the well-posedness of generalized martingale problems associated with critical stochastic differential equations.

\bigskip
\noindent 
\textbf{Keywords}: 
Harnack inequality, Nonlocal equations, Stochastic quasi-geostrophic equations,
Fractional 2D-Navier-Stokes equation, Critical SDEs\\

\noindent
  {\bf AMS 2020 Mathematics Subject Classification:}  Primary: 60H10, 35R09; Secondary: 60G51
\end{abstract}

\maketitle \rm

\tableofcontents

\section{Introduction}

In the celebrated work \cite{Di-Li}, DiPerna and Lions  established the existence, uniqueness, and stability of renormalized solutions to the transport equation on $\mathbb{R}^d$:
$$
\partial_t u = b \cdot \nabla u \quad \hbox{with } \  u(0) = \varphi,
$$
where   $b \in L^1_{loc}(\mR_+; W^{1,1}_{\text{loc}}(\mathbb{R}^d))$ 
has bounded divergence and $\varphi$ is a finite measurable function.
As a key application, they proved the well-posedness, in the Lebesgue sense,
of the ordinary differential equation
$$
\dot X_t(x) = b(t, X_t(x)) \quad \hbox{with } \   X_0(x) = x.
$$
These results have been subsequently extended in several directions.
Ambrosio \cite{Am1} generalized them to the case where the drift $b$ is of bounded variation, while the second author \cite{Zh1} treated stochastic differential equations.

It is well known that noise can exert a regularizing effect on irregular ODEs. For instance, consider the stochastic differential equation on $\R^d$
driven by a Brownian motion $B$:
$$
\mathrm{d} X_t = b(t, X_t) \mathrm{d}t + \mathrm{d} B_{2t}, \quad X_0 = x.
$$
Even when $b$ is merely bounded and measurable or possesses singularities, 
there exists a unique strong solution for every initial point $x$ (see \cite{Ve79, Kr-Ro} etc.).
This phenomenon is intimately connected to the existence and uniqueness of solutions to the second-order parabolic equation:
$$
\partial_t u = \Delta u + b \cdot \nabla u + f, \quad u(0) = \varphi,
$$
where $\Delta$ denotes the Laplacian, which is the infinitesimal generator of the Brownian motion $B_{2t}$.

In this paper, we are interested in the  
regularization phenomenon
under the influence of  $\alpha$-stable noises in the presence of a divergence-free drift $b$.  
Here, the parameter $\alpha$ ranges over $(0,2)$, where smaller values of $\alpha$ correspond to weaker noise.  
More precisely, we study the following fractional diffusion transport equation on $\mathbb{R}^d$:
\begin{align}\label{AA6} 
\partial_t u = \Delta^{\alpha/2} u + b \cdot \nabla u + f, \quad u(0) =\varphi,
\end{align}
as well as the following stochastic differential equation (SDE) driven by an isotropic $\alpha$-stable process:
\begin{align}\label{SDE1}
\mathrm{d} X_t = b(t, X_t) \, \mathrm{d}t + \mathrm{d} L^{(\alpha)}_t, \quad X_0 = x,
\end{align}
where $\alpha \in (0,2)$ and $\Delta^{\alpha/2} := -(-\Delta)^{\alpha/2}$ is the fractional Laplacian
on $\R^d$, which  is also
the infinitesimal generator of the 
$d$-dimensional
isotropic $\alpha$-stable process $L^{(\alpha)}$.

When $\alpha \in (1,2)$, the drift term can be controlled by the nonlocal diffusion term
$\Delta^{\alpha/2}$,
and the regularity of the PDE \eqref{AA6} has been extensively studied via Duhamel's representation (see, for example, \cite{Zh13a}). This regime is referred to as the subcritical case.  
When $\alpha = 1$,   
the drift term  $b \cdot \nabla $ and the fractional Laplacian $\Delta^{1/2}$  are of the same order,
leading to the critical case.  
For $\alpha \in (0,1)$, the drift term  is of higher order than the fractional Laplacian $\Delta^{\alpha/2}$, 
and this is termed the supercritical case (see \cite{CSZ18} and \cite{CZZ21} for studies related to the SDE \eqref{SDE1}).

In the subcritical case $\alpha \in (1,2)$ with $b$ being $\beta$-H\"older continuous where $\beta \in (1 - \frac{\alpha}{2}, 1)$, strong well-posedness of the SDE \eqref{SDE1} was established in \cite{P12} (see also \cite{Zh13a}).  
More recently, weak well-posedness for all $\alpha \in (0,2)$ and $\beta$-H\"older continuous $b$ with $\beta \in ((1 - \alpha)^+ \vee 0, 1]$ was shown in \cite{CSZ18} and \cite{CZZ21}.  
Two-sided heat kernel estimates for SDE \eqref{SDE1} were also obtained in \cite{MZ20}.  
We also mention that in \cite{ABM20}, the authors studied strong well-posedness for $\alpha \in (1,2)$ and $b \in \bC^\beta$ with $\beta > \frac{1 - \alpha}{2}$.
An important open question on the weak well-posedness in the critical case still remains
 where $\beta = 1 - \alpha$ for $\alpha \in (0,1)$.

We now perform a scaling analysis. Let $u = u(t, x)$ be a solution to \eqref{AA6} on $\mathbb{R}^{1+d}$. For $\lambda > 0$, define
$$
u_\lambda(t,x) := u(\lambda^\alpha t, \lambda x), \quad 
b_\lambda(t,x) := \lambda^{\alpha-1} b(\lambda^\alpha t, \lambda x), \quad 
f_\lambda(t,x) := \lambda^\alpha f(\lambda^\alpha t, \lambda x).
$$
Then
$$ 
\p_t u_\lambda=\Delta^{\alpha/2} u_\lambda+ b_\lambda \cdot\nabla u_\lambda+ f_\lambda.
$$
We introduce the space-time norm: for $p,q\in[1,\infty]$, 
$$
\|f\|_{\mL^q_t\mL^p_x}:=\left(\int_\mR\left(\int_{\mR^d}|f(t,x)|^p\dif x\right)^{q/p}\dif t\right)^{1/q}
$$
with the standard modification when $q=\infty$ or $p=\infty$.
Notice that by the change of variables,
$$
\|b_\lambda\|_{\mL^q_t\mL^p_x}=\lambda^{\alpha-1-\frac{d}{p}-\frac{\alpha}{q}}\|b\|_{\mL^q_t\mL^p_x}
$$
and
$$
\|f_\lambda\|_{\mL^q_t\mL^p_x}=\lambda^{\alpha-\frac{d}{p}-\frac{\alpha}{q}}\|f\|_{\mL^q_t\mL^p_x}.
$$
For the drift $b$, according to the sign of $\alpha - 1 - \frac{d}{p} - \frac{\alpha}{q}$, we 
make the following classification:
\begin{enumerate}[\rm (i)]
\item (Subcritical case): $\alpha\in(1,2)$ and $\frac{d}{p}+\frac{\alpha}{q}<\alpha-1$.

\item (Critical case): $\alpha\in[1,2)$ and $\frac{d}{p}+\frac{\alpha}{q}=\alpha-1$.

\item (Supercritical case): $\alpha\in(0,2)$ and $\alpha-1<\frac{d}{p}+\frac{\alpha}{q}<\alpha$.
\end{enumerate}

Since the critical and supercritical cases are of special interest in physical models, this paper focuses primarily on studying the well-posedness of the PDE \eqref{AA6}
and the associated SDEs \eqref{SDE1} in these regimes. We discuss them in the following subsections.

\subsection{Stochastic quasi-geostrophic equations}

The 2D quasi-geostrophic equation is an important model in fluid dynamics that describes the evolution of geophysical fluids. It takes the following form:
$$
\partial_t\theta = \mathcal{R}\theta \cdot \nabla\theta, \quad \theta(0) = \theta_0,
$$
where $\theta:\mathbb{R}^2\to\mathbb{R}$ represents the temperature and 
$$
\mathcal{R}\theta := (-\partial_2\Delta^{-1/2}\theta, \partial_1\Delta^{-1/2}\theta) =: (-\mathcal{R}_2\theta, \mathcal{R}_1\theta),
$$
with $\mathcal{R}_j$ denoting the $j$-th Riesz transform. The dissipative quasi-geostrophic equation is given by
$$
\partial_t\theta = \Delta^{\alpha/2}\theta + \mathcal{R}\theta \cdot \nabla\theta, \quad \theta(0) = \theta_0,
$$
where $\alpha\in(0,2)$. It is well known that for $\alpha>1$ (subcritical case), the above initial value problem has a global $C^\infty$-solution for any smooth initial value $\theta_0$ (see \cite{CW99, R95}).
The critical case $\alpha=1$, which arises in geophysical studies of strongly rotating fluid flows (cf. \cite{C02}), is solved by Kiselev-Nazarov-Volberg \cite{KNV07} and Caffarelli-Vasseur \cite{CV10} using completely different methods.

On the other hand, R\"ockner, Zhu and Zhu \cite{RZZ15} considered the following stochastic fractionally dissipative quasi-geostrophic equation on the torus $\mathbb{T}^2$:
\begin{align}\label{ACV7}
\mathrm{d}\theta = ( \Delta^{\alpha/2}\theta + \mathcal{R}\theta \cdot \nabla\theta)\, \mathrm{d}t + G(\theta(t))\, \mathrm{d}W_t, \quad \theta(0) = \theta_0.
\end{align}
Therein, the existence of martingale solutions for all $\alpha\in(0,2)$ and the existence and uniqueness of strong solutions and invariant measures in the subcritical case $\alpha\in(1,2)$ are established. 
Later, their results were extended to the whole space and more general cases by Brze\'zniak and Motyl \cite{BM19}.
However, the well-posedness of SPDE \eqref{ACV7} in the critical case $\alpha=1$ remains open.

We mention that Maekawa and Miura \cite{MM13, MM13b} studied the fundamental solution to the following nonlocal equation:
$$
\partial_t u = \Delta^{\alpha/2}u + b \cdot \nabla u,
$$
where $b$ is divergence-free and $\alpha \in (0,2)$. The existence and continuity of fundamental solutions, as well as pointwise upper bound estimates, are obtained when $b$ has critical regularity.
Moreover,  a priori H\"older regularity for the above advection fractional-diffusion equations was studied in 
\cite{CW08} for $\alpha \in (0,1)$ and in \cite{DS18} for $\alpha \in (1,2)$.
Their methods are based on De-Giorgi's iteration and Caffarelli-Silvestre's extension theorem \cite{CS07}
for the fractional Laplacian operator $\Delta^{\alpha/2}$.

\subsection{Fractional Navier-Stokes equation in $\mR^2$}

Consider the following two-dimensional fractional Navier-Stokes equation:
\begin{align}\label{2NS}
\partial_t \mathbf{u} = \Delta^{\alpha/2} \mathbf{u} + \mathbf{u} \cdot \nabla \mathbf{u} + \nabla p
 \quad \text{with} \quad \nabla \cdot \mathbf{u} = 0,
\end{align}
where $\mathbf{u} = (u_1, u_2)$ denotes the velocity field of the fluid, $p$ denotes the pressure, and $\alpha \in (0,2)$.
Here, $\Delta^{\alpha/2} := -(-\Delta)^{\alpha/2}$ is the fractional Laplacian of order $\alpha$.
It is known that if the initial velocity $\mathbf{u}_0$ is smooth, then there exists a unique smooth global solution to \eqref{2NS} (see \cite{Zh13b}).

Taking 
curl on
both sides of \eqref{2NS}, we find that the vorticity $\rho := \operatorname{curl} \mathbf{u} := \partial_1 u_2 - \partial_2 u_1$ satisfies
$$
\partial_t \rho = \Delta^{\alpha/2} \rho + \mathbf{u} \cdot \nabla \rho = \Delta^{\alpha/2} \rho + \nabla \cdot (\rho \mathbf{u}).
$$
Moreover, by the Biot--Savart law (cf.~\cite{Ma-Be}),
$$
\mathbf{u}(t, x) = \int_{\mathbb{R}^2} K_2(x - y) \rho(t, y) \, \mathrm{d} y =: (K_2 * \rho)(t, x),
$$
where
\begin{align}\label{K2}
K_2(x) := \frac{1}{2\pi} \left( \frac{-x_2}{|x|^2}, \frac{x_1}{|x|^2} \right).
\end{align}
Hence, $\rho$ satisfies the following nonlinear integro-differential equation:
\begin{align}\label{2VNS}
\partial_t \rho = \Delta^{\alpha/2} \rho + \nabla \cdot (\rho \cdot (K_2 * \rho)).
\end{align}

Suppose that $\rho(0, x) \geq 0$ and $\int_{\mathbb{R}^2} \rho(0, x) \, \mathrm{d} x = 1$.
By the maximum principle and integrating both sides of \eqref{2VNS} with respect to $x$, we find that for all $t > 0$,
$$
\rho(t, x) \geq 0 \quad \hbox{and} \quad 
\int_{\mathbb{R}^2} \rho(t, x) \, \mathrm{d} x = \int_{\mathbb{R}^2} \rho(0, x) \, \mathrm{d} x = 1.
$$
This  implies that $\{\rho(t, \cdot)\}_{t \geq 0}$ is a family of probability 
density functions on $\R^2$. 
By the superposition principle \cite{RXZ20}, if \eqref{2VNS} admits a density solution $\rho$, 
the following McKean--Vlasov SDE 
 admits a solution $X_t$:
\begin{align}\label{SDE0}
\mathrm{d} X_t = \left( \int_{\mathbb{R}^2} K_2(X_t - y) \rho(t, y) \, \mathrm{d} y \right) \mathrm{d} t + \mathrm{d} L^{(\alpha)}_t.
\end{align}
Note that the kernel $K_2$ is singular 
at the origin.

When $\alpha = 2$ and the initial vorticity is a finite Radon measure, the existence of solutions to the PDE \eqref{2VNS} was established by Giga, Miyakawa, and Osada in \cite{Gi-Mi-Os}, and uniqueness was proved by Gallagher and Gallay in \cite{Ga-Ga}. Their proofs rely heavily on heat kernel estimates for the operator $\Delta + b \cdot \nabla$ when $b$ is divergence-free and has a specific form (see \cite{O87}).
Recently,   the existence of weak solutions to \eqref{2VNS} with measure-valued initial data was obtained 
in \cite{Zh21} by directly constructing weak solutions to \eqref{SDE0}. Some regularity results were also 
established 
there. One of the goals of this paper is to derive analogous results for the fractional vorticity equation \eqref{2VNS} with $\alpha \in (1,2)$ by studying the associated SDE \eqref{SDE0}.

\medskip

\subsection{Main results}
We consider the following nonlocal equation on $\R^d$
with lower regularity on $b$ and $f$:
\begin{align}\label{PDE0}
\p_t u=\sL^\alpha_t u+b\cdot \nabla u+f,
\end{align}
where 
$$
\sL^\alpha_t f(x):={\rm p.v.}\int_{\mR^d}(f(x+y)-f(x)) K(t,y)\dif y,
$$
and  
$ K(t,y):\mR_+\times\mR^d\to(0,\infty)$ 
is a symmetric measurable function   in $y$ and satisfies
\begin{align}\label{KK0}
 0< \kappa_0|y|^{-d-\alpha}\leq K(t,y)\leq \kappa_1|y|^{-d-\alpha}.
\end{align}
Note that when $\kappa_0=\kappa_1$, up to a multiplicative constant, $\sL^\alpha_t$ is the standard fractional Laplacian 
$\Delta^{\alpha/2}:=-(-\Delta)^{\alpha/2}$.  
By symmetry of $K(t, y)$ in $y$, it is easy to see that  
$$
\sL^\alpha_t f(x)=\frac12\int_{\mR^d}(f(x+y)+f(x-y)-2f(x)) K(t,y)\dif y,
$$ 
and for $f\in C^2_b(\R^d) \cap L^2(\R^d)$,
\begin{align}\label{SY1}
\int_{\mR^d}\sL^\alpha_tf(x) g(x)\dif x=\int_{\mR^d}f(x)\sL^\alpha_t g(x)\dif x=-\cE_t^{(\alpha)}(f,g),
\end{align}
where
\begin{align}\label{SY91}
\cE_t^{(\alpha)}(f,g)
:= \frac12 \int_{\mR^d}\int_{\mR^d}(f(x)-f(y))(g(x)-g(y))K(t,x-y)\dif x\dif y .
\end{align}
Suppose that $g$ has support in $B_r$. Then for any $R\geq r$, we have
\begin{align}\label{AA12}
\begin{split}
\cE_t^{(\alpha)}(f,g)
&= \frac12 \int_{B_R}\int_{B_R}(f(x)-f(y))(g(x)-g(y))K(t,x-y)\dif x\dif y \\
&\quad+ \int_{B_R}\left(\int_{B^c_R}(f(x)-f(y)) K(t,x-y)\dif y \right)g(x)\dif x.
\end{split}
\end{align}
We assume that
\begin{align}\label{C2} 
  | b | +|  \nabla b | \in\mL^1_{\rm loc} (\R^d)
\quad \hbox{and} \quad  f\mbox{ is a distribution}.
\end{align}
We first introduce the following notion of super- and sub- solutions.

\bd \label{D:1.1}
A bounded measurable function  
$u(t, x)$ on $\mR \times \mR^d$
is called a super (resp. sub)-solution
of PDE \eqref{PDE0}
if for any nonnegative $\varphi\in C^\infty_c(\mR^{1+d})$,
\begin{align}\label{Def0}
\begin{split}
-\int_{\mR^{1+d}} u\p_t\varphi  \ \geq \hbox{ {\rm (resp. $\leq$)} } 
\int_{\mR^{1+d}} u\sL_t^\alpha\varphi
-\int_{\mR^{1+d}}(b\cdot\nabla \varphi+ \varphi \, \div b )u+\int_{\mR^{1+d}} f\varphi.
\end{split}
\end{align}
If $u$ is both a super and sub-solution, then it is called a solution.
\ed

\br\rm
Under \eqref{C2} and 
the boundedness assumption on $u$,
each term in the right hand side of \eqref{Def0} is well defined. 
Indeed, for $\varphi\in C^\infty_c(\mR^{1+d})$, we  have 
for each $t\in \R$, 
\begin{align*}
\|u\sL_t^\alpha\varphi\|_{\mL^1}
&\leq  \kappa_1\|u\|_\infty  \left\| \int_{\{|y|\leq 1\}}
\frac{\varphi(\cdot+y)+\varphi(\cdot-y)-2\varphi }{|y|^{d+\alpha}}\dif y\right\|_{\mL^1}\\
& +  \kappa_1\|u\|_\infty   \int_{\{ |y|>1\}}
\frac{\|\varphi(\cdot+y)+\varphi(\cdot-y)-2\varphi\|_{\mL^1}}{|y|^{d+\alpha}}\dif y
\leq  C<\infty. 
\end{align*}
To make sense of $\int_{\mR^{1+d}}\varphi \, \div b\,  u$, it suffices to assume 
$|  \div b | \in\mL^1_{\rm loc} (\R^d)$. 
However, we imposed a  stronger condition
 $|  \nabla b | \in\mL^1_{\rm loc} (\R^d)$ on $u$
in Definition \ref{D:1.1} because we later
 need a commutator estimate \eqref{Mo2} 
involving $u$
  (see Lemma \ref{Le213} below).
\er

\br\rm
If $u$ is a supper  (resp. sub)-solution, then so is $u+\kappa$  for any $\kappa\in\mR$.
\er

To state our main results, we need some notations. For  $r>0$ and $t_0\in\mR$,   define
\begin{align}\label{CY1}
B_r:=\{x\in\mR^d: |x|\leq r\}, \quad
Q^+_r(t_0):=[t_0,t_0+r]\times B_r, \quad  Q^-_r(t_0):=[t_0-r,t_0]\times B_r,
\end{align}
and
\begin{align}\label{CY11}
Q_r(t_0):=Q^+_r(t_0)\cup Q^-_r(t_0)=[t_0-r,t_0+r]\times B_r,
\quad  Q_r:=Q_r(0).
\end{align}
For a function $g:\mR^{1+d}\to\mR$, define 
\begin{align}\label{Tail}
{\rm Tail}(g;r)(t):=\int_{B^c_r}\frac{|g(t,y)|}{|y|^{d+\alpha}}\dif y.
\end{align}
For $\alpha\in(0,2)$ and $\beta\in[0, \alpha/2)$, we introduce the following index set for later use:
\begin{align}\label{Ind}
\mI^\beta_\alpha:=\Big\{(q,\bbp)\in
(1, \infty) \times (1, \infty)^d: \ 
\tfrac\alpha{q}+|\tfrac{1}{\bbp}|<\alpha-\beta\Big\}.
\end{align}
See \eqref{AM9}-\eqref{e:2.3} and \eqref{e:2.7} below for the notations of
$\| f\|_{\mL^\bbp}$, $\frac1{\bbp}$, $|\bbp|$, $\bH^\beta_{\bbp}$ and $\| u\|_{\mL^{q}_t\mL^{\bbp}_x}$.
We consider the following condition on $b$. 

\begin{enumerate}[{\bf (H$_b$)}]
\item  
In addition to \eqref{C2}, we assume that for some $(q_1,\bbp_1), (q_2,\bbp_2)\in\mI^0_\alpha$,
$$
\kappa_b:=\|b\b1_{Q_2}\|_{\mL^{q_1}_t\mL^{\bbp_1}_x}+\|\div b\b1_{Q_2}\|_{\mL^{q_2}_t\mL^{\bbp_2}_x}<\infty.
$$
\end{enumerate}

We will use the following parameter set:
$$
\Theta:=(\kappa_0,\kappa_1,\kappa_b,\bbp_1,q_1,\bbp_2,q_2,d,\alpha).
$$
Our first main result is the following local upper bound estimate of sub-solution to PDE \eqref{PDE0}.

\bt\label{Le32}
Let $\alpha \in (0, 2)$,
$\beta\in[0, \alpha/2)$ and $(q_0,\bbp_0)\in\mI_\alpha^\beta$.
Under \eqref{KK0} and {\bf (H$_b$)}, for any $p>0$, there is a constant
$C=C(\Theta,\beta,q_0,\bbp_0, p)>0$ such that for any sub-solution $u$ of PDE \eqref{PDE0} and $t\in\mR$,
$$
\|u^+\b1_{Q_1}\cI_t\|_{\mL^\infty}+\|\Delta^{\alpha/4}(u^+\chi_1)\cI_t\|_{\mL^2}
\lesssim_C \|u^+\b1_{Q_2}\cI_t\|_{\mL^p}
+\|{\rm Tail}(u^+; 1)\cI_t\|_{\mL^2_t}+\|\chi_2f\cI_t\|_{\mL^{ q_0}_t(\bH^{-\beta}_{\bbp_0})},
$$
where $\cI_t:=\b1_{(-\infty,t]}$ and for $j=1,2$, $\chi_j\in C^\infty_c(Q_{2j})$, are nonnegative and equal to $1$ on $Q_{j}$.
\et

\br\rm
When $\alpha=2$, this type of maximum estimate was established in \cite{ZZ18}. It is important to note that $f$ is allowed to be a distribution.
\er

Next we consider the following condition. 

\begin{enumerate}[{\bf (H$'_b$)}]
\item  In addition to \eqref{C2}, we assume that for some 
$(q_1,\bbp_1), (q_2,\bbp_2)\in\mI^0_\alpha$,
$$
\kappa_b:=\|b\b1_{Q_4}\|_{\mL^2}+\|\div b\b1_{Q_4}\|_{\mL^2}+
\|b\b1_{Q_4}\|_{\mL^{q_1}_t\mL^{\bbp_1}_x}+\|\div b\b1_{Q_4}\|_{\mL^{q_2}_t\mL^{\bbp_2}_x}<\infty.
$$
\end{enumerate}

Our second main result is the following weak Harnack inequality for super-solutions of PDE \eqref{PDE0}.

\bt\label{Th41}
Let $\alpha \in (0, 2)$ and  
$(q_0,\bbp_0)\in\mI_\alpha^0$.
Under  \eqref{KK0} and {\bf (H$'_b$)}, there are $p_0>0$ and $C=C(\Theta,
q_0,\bbp_0)>0$ 
such that  for any super-solution 
$u$ of \eqref{PDE0} in $Q_4 $  with $u\b1_{Q_4 }\geq 0$,
\begin{align}\label{WH1}
\|u\b1_{Q_{3/2}^+(-2)}\|_{\mL^{p_0}}\lesssim_C {\rm inf}_{Q_{3/2}^-(2)}u+\|f\b1_{Q_2}\|_{\mL^{q_0}_t\mL^{\bbp_0}_x}+\|{\rm Tail}(u^-,1)\b1_{[-2,2]}\|_{\mL^1_t}.
\end{align}
Moreover,  for any solution 
$u$ of \eqref{PDE0} in $Q_4 $  with $u\b1_{Q_4 }\geq 0$, 
the following Harnack inequality holds:
\begin{align}\label{WH2}
{\rm sup}_{Q_1^+(-2)}u\lesssim_C {\rm inf}_{Q_1^-(2)}u
+\|f\b1_{Q_2}\|_{\mL^{q_0}_t\mL^{\bbp_0}_x}
+\|u\b1_{[-4,4]}\|_{\mL^\infty}.
\end{align}
\et

When $b$ is zero, $f$ is bounded measurable, and 
$\sL^\alpha_t$
is replaced with  
a symmetric nonlocal operator whose energy on each ball is uniformly comparable 
to that of $\Delta^{\alpha/2}$,
Felsinger and Kassmann \cite{FK13} studied
weak Harnack inequality and H\"older regularity for caloric functions.
 Their results were later
 extended to more general symmetric nonlocal operators in \cite{St19}. In a recent work \cite{KW25}, Kassmann and Weidner showed a full Harnack inequality version with no tail term 
for globally nonnegative caloric functions.
 However, all the above results do not 
 cover the case of singular drift $b$  nor the case  of distributional source term  $f$.

\subsection{Organization of the paper}

The remainder of this paper is structured as follows. 
Section 2 introduces the necessary function spaces, with particular emphasis on the De Giorgi and Moser classes, which play crucial roles in {\red our} subsequent analysis. 
In Section 3, we establish an upper bound estimate for subsolutions of PDE \eqref{PDE0} by verifying their membership in a suitable De Giorgi class. 
Section 4 is devoted to proving a weak Harnack inequality for supersolutions of PDE \eqref{PDE0} by demonstrating that they belong to an appropriate Moser class. 
Section 5 addresses the global well-posedness and $L^\infty$-estimates for PDE \eqref{PDE0} under supercritical conditions, as well as H\"older regularity in the critical case when $\alpha \in (0,1]$. 
In Section 6, we apply the results from Section 5 to establish well-posedness for stochastic quasi-geostrophic equations. 
Section 7 presents existence results for weak solutions to a class of supercritical distribution-dependent SDEs, including the two-dimensional fractional Navier-Stokes equations. 
Finally, Section 8 is concerned with
 the existence and uniqueness of generalized martingale solutions for critical linear SDEs.

We conclude this introduction by summarizing some notations used throughout the paper.
\begin{itemize}
    \item The symbol $:=$ denotes a definition. For $c > 0$, the notation $f \lesssim_c g$ means $f \leq c g$.
    \item $\langle \cdot, \cdot \rangle$ represents the dual pairing between distributions and test functions.
\item For $a, b\in \R$, $a \vee b :=\max \{a, b\}$, $a\wedge b:= \min \{a, b\}$ and $a^+:=a\vee 0$, 
$a^-:=(-a)\vee 0$. 
 \end{itemize}

\section{Preliminary}

\subsection{Function spaces}
Let $N\in\mN$. For multi-index $\bbp=(p_1,\cdots,p_N)\in(0,\infty]^N$, we define
\begin{align}\label{AM9}
\|f\|_{\mL^\bbp}:=\Bigg(\int_\mR\Bigg(\int_\mR\cdots\left(\int_\mR |f(z_1,\cdots,z_N)|^{p_N}\dif z_N\right)^{\frac{p_{N-1}}{p_N}}\cdots
\dif z_2\Bigg)^{\frac{p_1}{p_2}}\dif z_1\Bigg)^{\frac{1}{p_1}},
\end{align}
with the standard modification when $p_i=\infty$, for some $i=1,\cdots, N$,
which is also called the mixed $L^\bbp$-space (see \cite{BP61}).
When $\bbp=(p,\cdots,p)\in(0,\infty]^N$, 
$\mL^\bbp$ is the usual $L^p (\R^N)$ space and 
we shall simply write it as $  \mL^p$.
 (
 Note that  when there is some $1\leq i\leq N$ so that
  $p_i \in (0, 1)$ and $p\in (0, 1)$, $\mL^\bbp$ and $L^p$ are not Banach spaces
as in these cases the norms defined above do not satisfy the triangle inequality.
However,  for most use later, these spaces are for indices larger than or equal 1.)
For multi-indices $\bbp, \bbq\in(0,\infty]^N$, we shall use the following notations:
$$
\bbp + \bbq := (p_1+q_1, \cdots, p_N+q_N), 
\ \  \bbp\cdot\bbq:=\sum_{i=1}^Np_iq_i,
$$
and 
$$
\bbp>\bbq\ \ (\hbox{\rm resp.  } \bbp\geq\bbq;\ \bbp=\bbq)\Longleftrightarrow p_i>q_i
 \ (\hbox{\rm resp.  }  p_i\geq q_i;\ p_i=q_i) \ \hbox{ for all } i=1,\cdots,N.
$$
Moreover,  we use bold number to denote constant vector in $\mR^N$, for example, 
$$
\boldsymbol{1}=(1,\cdots,1),\ \ \boldsymbol{2}=(2,\cdots,2),
$$
and also denote
\begin{equation}\label{e:2.2}
|\bbp | :=\bbp\cdot \boldsymbol{1},  \quad 
\tfrac1{\bbp}:=\big(\tfrac1{p_1},\cdots,\tfrac1{p_N}\big).
\end{equation} 
For any multi-indices $\bbp,\bbq,\bbr\in(0,\infty]^N$ with $\tfrac1\bbp+\tfrac1\bbr=\tfrac1\bbq$, the following H\"older's inequality holds
$$
\|fg\|_{\mL^\bbq}\leq \|f\|_{\mL^\bbp}\|g\|_{\mL^\bbr}.
$$
For any multi-indices $\bbp,\bbq,\bbr\in[1,\infty]^N$ with $\tfrac1\bbp+\tfrac1\bbr=\boldsymbol{1}+\tfrac1\bbq$, the following Young's inequality holds
$$
\|f*g\|_{\mL^\bbq}\leq \|f\|_{\mL^\bbp}\|g\|_{\mL^\bbr}.
$$

For $\beta\in\mR$ and $\bbp\in(1,\infty)^N$, we define the following Bessel potential space
\begin{equation}\label{e:2.3}
\bH^\beta_{\bbp}:=(I-\Delta)^{-\beta/ 2}(\mL^\bbp).
\end{equation}
The dual space of $\bH^\beta_{\bbp}$ is 
$\bH^{-\beta}_{\bbq}$, where $\frac1\bbp+\frac1\bbq= \bf 1$.
For $\beta\in(0,2)$ and $\bbp\in(1,\infty)^N$, by H\"ormander's multipliers theorem (see \cite[Corollary 3.1]{CGN}), an equivalent norm of $\bH_{\bbp}^\beta$ is given by
\begin{align}\label
{e:2.4} 
\|f\|_{\bH^\beta_{\bbp}}
  := \|f\|_{\mL^\bbp}+\|\Delta^{\beta/2}f\|_{\mL^\bbp}.
\end{align}
  We have the following interpolation inequality of  Gagliado-Nirenberge's type.
\bl
Let $\bbp,\bbq,\bbr\in(1,\infty)^N$ and 
$s,s_0,s_1\in\mR$, $\theta\in[0,1]$. Suppose that
\begin{align}\label{DD01}
\tfrac1{\bbp}\leq\tfrac{1-\theta}{\bbq}+\tfrac{\theta}{\bbr} 
\quad \hbox{and} \quad 
 s-|\tfrac{1}{\bbp}|=(1-\theta)\big(s_0-|\tfrac{1}{\bbq}|\big)+\theta\big(s_1-|\tfrac{1}{\bbr}|\big).
\end{align}
Then there is a constant $C=C(\bbp,\bbq,\bbr,s,s_0,s_1,\theta)>0$ such that
\begin{align}\label{Sob}
\|f\|_{\bH^s_{\bbp}}\lesssim_C\|f\|^{1-\theta}_{\bH^{s_0}_{\bbq}}\|f\|^\theta_{\bH^{s_1}_{\bbr}}.
\end{align}
\el

Let $\mR^N=\mR^{1+d}$ be the space-time space. 
For a function $f(t, x)$ defined on $\mR^{1+d}$ and $\nu \in (0, \infty]$ and $\bbp\in (0, \infty]^{d}$,
we define
\begin{equation} \label{e:2.7} 
\| f\|_{\mL^\nu_t\mL^\bbp_x} := 
\begin{cases}
 \left( \int_{\mR} \left( \|  f(t, \cdot)\|_{\mL^\bbp (\mR^d)} \right)^\nu dt \right)^{1/\nu}
 \quad &\hbox{if } 0<\nu <\infty ,\\
 {\rm esssup}_{t\in \mR}  \|  f(t, \cdot)\|_{\mL^\bbp (\mR^d)}  &\hbox{if } \nu=\infty.
 \end{cases} 
\end{equation}
For $\alpha\in(0,2)$, we introduce the following energy space for later use:
\begin{align}\label{VA}
\sV_\alpha:=\Big\{f(t, x): \,  \|f\|_{\sV_\alpha}:=\|f\|_{\mL^\infty_t\mL^2_x}+\|\Delta^{ \alpha /4}f\|_{\mL^2}<\infty\Big\}.
\end{align}

We have the following    Sobolev's type inequality.

\bl\label{Le22}
Let $\alpha\in(0,2)$. For $(\nu,\bbr)\in[2,\infty)^{1+d}$ with $\frac\alpha\nu+|\frac1 \bbr|=\tfrac{d}2$, it holds that
$$
\|u\|_{\mL^\nu_t\mL^\bbr_x}\lesssim_C \|u\|_{\mL^\infty_t\mL^2_x}+\|u\|_{\mL^2_t\mL^\bbp_x}\lesssim_C \|u\|_{\sV_\alpha},
$$
where $\bbp\in[2,\infty)^d$ is defined by $\frac2\bbp+\frac{\nu-2}{\bb2}=\frac \nu \bbr$ and $C=C(d,\alpha,\nu,\bbr)>0$.
\el
\begin{proof}
Let $\theta=\frac{2}\nu$. Then $\frac\theta\bbp+\frac{1-\theta}{\bb2}=\frac 1 \bbr$. By H\"older's inequality we have
$$
\|u\|_{\mL^\bbr_x}\leq \|u\|_{\mL^2_x}^{1-\theta}\|u\|_{\mL^\bbp_x}^\theta.
$$
Hence, 
$$
\|u\|_{\mL^\nu_t\mL^\bbr_x}\lesssim_C \|u\|^{1-\theta}_{\mL^\infty_t\mL^2_x}\|u\|^{\theta}_{\mL^2_t\mL^\bbp_x}.
$$
which in turn gives the first inequality by Young's inequality. Since $|\frac1\bbp|=\frac{d-\alpha}2$ by the assumption,  we have by \eqref{Sob} 
with $s=\theta=0$ and $s_0=\alpha/2, q=2$,
\begin{align}\label{DS1}
\|f\|_{\mL^\bbp_x}\lesssim\|f\|_{\bH^{\alpha/2}_2}\asymp  \|f\|_{\mL^2_x}+\|\Delta^{\alpha/4}f\|_{\mL^2_x}.
\end{align}
The second inequality follows.
\end{proof}

We also have the following Sobolev inequality in balls. 

\bl
For any $\alpha\in(0,1)$ and $\bbp\in[2,\infty)^d$ with 
$|\frac1\bbp|=\frac {d-\alpha} 2$,  
there is a constant 
$C=C(d,\bbp,\alpha)>0$ such that for any $0<r<R<\infty$
and $\eta\in C^\infty_c(B_r)$,
\begin{align}\label{SS1}
\|f\eta\|_{\mL^\bbp}\lesssim_C 
\left( 1+ \frac{R^d}{(R-r)^{d+\alpha}} \right)
 \|f\eta\|_{\mL^2}+\left(\int_{B_R}\int_{B_R}\frac{((f\eta)(x)-( f\eta)(y))^2}{|x-y|^{d+\alpha}}\dif x\dif y\right)^{1/2}.
\end{align}
\el
\begin{proof}
Let $\eta\in C^\infty_c(B_r)$. By \eqref{DS1} we have
\begin{align*}
\|f\eta\|_{\mL^\bbp}\lesssim\|f\eta\|_{\mL^2}+\|\Delta^{\alpha/4}(f\eta)\|_{\mL^2}.
\end{align*}
By \eqref{SY1} and \eqref{AA12} we have
\begin{align*}
\|\Delta^{ \alpha/ 4}(f\eta)\|_{\mL^2}^2
&=\frac12
\int_{\mR^d}\int_{\mR^d}\frac{((f\eta)(x)-(f\eta)(y))^2}{|x-y|^{d+\alpha}}\dif x\dif y\\
&=\frac12
\int_{B_R}\int_{B_R}\frac{((f\eta)(x)-(f\eta)(y))^2}{|x-y|^{d+\alpha}}\dif x\dif y\\
&\quad+\int_{B_R}(f\eta)^2(x)\left(\int_{B^c_R}\frac{\dif y}{|x-y|^{d+\alpha}}\right)\dif x.
\end{align*}
Noting that for $|x|\leq r$ and $|y|\geq R$, 
$$
|x-y|^{-d-\alpha}\leq (1-\tfrac r R)^{-d-\alpha} |y|^{-d-\alpha},
$$
we have
$$
\int_{B^c_R}\frac{\dif y}{|x-y|^{d+\alpha}}\dif x\leq (1-\tfrac r R)^{-d-\alpha}\int_{|y|\geq R} \frac{\dif y}{|y|^{d+\alpha}}
=\frac{R^d}{(R-r)^{d+\alpha}}\int_{|y|\geq 1}\frac{\dif y}{|y|^{d+\alpha}}.
$$
The desired estimate now follows.
 \end{proof}

\subsection{De Giorgi and Moser classes}

Throughout this subsection,  we shall fix an open  multi-index set
$\sI\subset(1,\infty)^{N}$  and $Q:=(Q_\tau)_{\tau\in[1,2]}$  an increasing   family of   bounded open sets in $\mR^N$
with
\begin{align}\label{DD81} 
 \cap_{\sigma>\tau}Q_\sigma=\bar Q_{\tau}.
\end{align}
For a set $Q\subset\mR^N$, we use $\b1_Q$ to denote the indicator function of $Q$.

The following De Giorgi class associated with $\sI$ and $Q$ is introduced in \cite{Zh2}. 

\bd  [De Giorgi class] \label{Def6}
We say  a function $u\in L^1(Q_2)$ is 
 in the De Giorgi class $\cD\cG^+_\sI(Q)$ 
if there are $\bbp_i\in\sI$,  $i=1,\cdots,m$, some $j\in\{1,\cdots, m\}$
and $\gamma, \cA\geq 0$ such that for any $\bbp\in\sI$, there is a constant
 $C_\bbp   = C_\bbp (Q)>0$ such that
for any $1\leq\tau<\sigma\leq 2$ and $\kappa\geq 0$,
\begin{align}\label{DX-1}
(\sigma-\tau)^{\gamma}\|\b1_{Q_{\tau}}(u-\kappa)^+\|_{\mL^\bbp}\lesssim_{C_\bbp}
\sum_{i=1}^{j}\|\b1_{Q_{\sigma}}(u-\kappa)^+\|_{\mL^{\bbp_i}}+\cA\sum_{i=j+1}^{m}\|\b1_{\{u>\kappa\}\cap Q_{\sigma}}\|_{\mL^{\bbp_i}}.
\end{align}
\ed

The following result is essentially due to De Giorgi and proven in  \cite[Theorem 2.10]{Zh2}.
\bt\label{TH20}
For any $u\in\cD\cG^+_\sI(Q)$ and $p>0$, there is a constant $C_0>0$  
depending only on $Q$, $p$, $N$ and the quantities  in the definition of $\cD\cG^+_\sI(Q)$ such that
\begin{align}\label{JK2} 
\|u^+\b1_{Q_1}\|_{\mL^\infty}\lesssim_{C_0}\|u^+\b1_{Q_2}\|_{\mL^p}+\cA.
\end{align}
\et

The following lemma can be found in \cite[Lemma 2.2.6]{SC1}, which will  be used to show the weak Harnack inequality. 

\bl\label{Le28}
Let $\mu$ be a probability measure over $Q_2$ and $C_0,\gamma>0$, $q\in(0,\infty]$. 
Let $h$ be a non-negative measurable function on $Q_2$. Suppose that for all $p\in(0,1\wedge \tfrac q2)$ and $1\leq\tau<\sigma\leq 2$,
\begin{align}\label{CON2}
\|h\b1_{Q_\tau}\|_{\mL^q(\mu)}\leq \left(C_0(\sigma-\tau)^{-\gamma}\right)^{1/p}\|h\b1_{Q_\sigma}\|_{\mL^p(\mu)},
\end{align}
and for all $\lambda>0$,
\begin{align}\label{CON1}
  \mu \left( \big\{(t,x)\in Q_2: \log h(t,x)>\lambda\big\} \right) \leq C_0/\lambda.
 \end{align}
Then for any $\tau\in[1,2)$, there is a constant $C=C(C_0,\tau,\gamma,q)>0$ such that
$$
\|h\b1_{Q_\tau}\|_{\mL^q(\mu)}\leq C.
$$
\el

To check condition \eqref{CON2}, we introduce the following Moser class for later use. 

\bd [Moser  class]\label{Def60}
Let $q_0\in(0,\infty]$.
We say a nonnegative function $u\in L^\infty(Q_2)$ is 
 in the Moser class $\cM^{q_0}_\sI(Q)$
if there are $\bbp_i\in\sI$, $i=1,\cdots,m$, $\gamma\geq 0$ and $\beta\geq 0$ 
such that for any $\bbp\in\sI$, there is a constant $C_\bbp>0$ such that
for any $1\leq\tau<\sigma\leq 2$ and $q\in(0,q_0)$,
\begin{align}\label{DX-01}
(\sigma-\tau)^{\gamma}\|\b1_{Q_{\tau}}u^q\|_{\mL^\bbp}\lesssim_{C_\bbp}(q^\beta+1)\sum_{i=1}^m\|\b1_{Q_{\sigma}}u^q\|_{\mL^{\bbp_i}}.
\end{align}
\ed
We have the following reverse H\"older inequality for $u\in \cM^{q_0}_\sI(Q)$.
 
\bt\label{Th21}
\begin{enumerate}[\rm (i)]
\item For any  $u\in\cM^\infty_\sI(Q)$, there are constants $\gamma_1, C_1>0$ 
depending only on  
the quantities in Definition \ref{Def60}
 such that 
for all $1\leq\tau<\sigma\leq 2$ and $p\in(0,1]$,
\begin{align}\label{JK22}
\|u\b1_{Q_\tau}\|_{\mL^\infty}\leq (C_1(\sigma-\tau)^{-\gamma_1})^{\frac1p}\|u\b1_{Q_\sigma}\|_{\mL^p}.
\end{align}

\item For any $u\in\cM^{q_0}_\sI(Q)$ with $0< q_0<\infty$,
there are positive constants $p_0, \gamma_2, C_2>0$
depending only on   
the quantities in Definition \ref{Def60}
such that 
for all $1\leq\tau<\sigma\leq 2$ and $p\in(0,1\wedge\frac{p_0}2)$,
\begin{align}\label{JK23}
\|u\b1_{Q_\tau}\|_{\mL^{p_0}}\leq(C_2(\sigma-\tau)^{-\gamma_2})^{\frac1p}\|u\b1_{Q_\sigma}\|_{\mL^p}.
\end{align}
\end{enumerate}
\et

\begin{proof}
Suppose $u\in\cM^{q_0}_\sI(Q)$ with $q_0\in (0, \infty]$. 
Let $\bbp_i\in\sI$ be as in Definition \ref{Def60}.
Since $\sI$ is open, one can find $\theta>1$ 
close to $1$ such that $\theta\bbp_i\in\sI$ for each $i=1,\cdots,m$.
Thus by \eqref{DX-01}, there are constants $C_0>0$  and $\beta \geq 0$
such that for all $1\leq\tau<\sigma\leq 2$ and $q\in(0,q_0)$,
\begin{align*}
(\sigma-\tau)^{\gamma}\sum_{i=1}^m\|\b1_{Q_{\tau}}u^q\|_{\mL^{\theta\bbp_i}}\lesssim_{C_0}(q^\beta+1)\sum_{i=1}^m\|\b1_{Q_{\sigma}}u^q\|_{\mL^{\bbp_i}}.
\end{align*}
Since for $a_1,a_2,\cdots,a_m\geq 0$,
$$
a_1\vee a_2\vee\cdots\vee a_m\leq a_1+a_2+\cdots+a_m\leq m(a_1\vee a_2\vee\cdots\vee a_m),
$$ 
the above inequality can be rewritten as
\begin{align*}
(\sigma-\tau)^{\gamma}\vee_{i=1}^m\|\b1_{Q_{\tau}}u\|^q_{\mL^{q\theta\bbp_i}}\leq mC_0 (q^\beta+1)\vee_{i=1}^m\|\b1_{Q_{\sigma}}u\|^q_{\mL^{q\bbp_i}}.
\end{align*}
In particular, if we set
\begin{align}\label{AA4}
\sM(\tau,q):=\vee_{i=1}^m\|\b1_{Q_{\tau}}u\|_{\mL^{q\bbp_i}},
\end{align}
then for all $1\leq\tau<\sigma\leq 2$ and $q\in(0,q_0)$,
\begin{align}\label{A1A4}
\sM(\tau,\theta q)\leq ((\sigma-\tau)^{-\gamma}mC_0 (q^\beta+1))^{1/q}\sM(\sigma,q).
\end{align}

(i) Suppose $q_0=\infty$. Fix $q\in(0,1)$, whose value will be determined below.
Let $\tau_n=\tau+2^{-n}(\sigma-\tau)$ and $q_n:=\theta^n q$. Iterating \eqref{A1A4} with 
$\tau=\tau_n$, $\sigma=\tau_{n-1}$ and $q=q_{n-1}$, 
we obtain
\begin{align*}
\sM(\tau,q_n)&\leq\sM(\tau_n,q_n)\leq
\Big((\tau_{n-1}-\tau_n)^{-\gamma}mC_0(q_{n-1}^\beta+1)\Big)^{1/q_{n-1}}\sM(\tau_{n-1},q_{n-1})\\
&\leq\cdots\leq\prod_{j=1}^{n}\Big((\tau_{j-1}-\tau_j)^{-\gamma}mC_0(q^\beta_{j-1}+1)\Big)^{1/q_{j-1}}\sM(\tau_0,q),
\end{align*}
where the first inequality is due to $Q_\tau\subset Q_{\tau_n}$.
Noting that
$$
\sM(\tau,\infty)=\|u\b1_{Q_\tau}\|_{\mL^\infty},\ \ \sM(\tau_0,q)=\vee_{i=1}^m\|\b1_{Q_{\sigma}}u\|_{\mL^{q\bbp_i}}
$$
and by $q\in(0,1)$ and $\theta>1$,
\begin{align*}
&\prod_{j=1}^\infty\Big((\tau_{j-1}-\tau_j)^{-\gamma}mC_0 (q^\beta_{j-1}+1)\Big)^{1/q_{j-1}}
=\prod_{j=1}^\infty\Big(2^{\gamma j}(\sigma-\tau)^{-\gamma}mC_0 (\theta^{\beta(j-1)}q^\beta+1)\Big)^{\theta^{1-j}/q}\\
&\qquad\leq\prod_{j=1}^\infty\Big(2^{\gamma j}(\sigma-\tau)^{-\gamma}2mC_0 \theta^{\beta j}\Big)^{\theta^{1-j}/q}
=\Big(\e^{\sum_{j=1}^\infty\theta^{-j}(j\ln (2^\gamma\theta^\beta)+\ln ((\sigma-\tau)^{-\gamma}2mC_0))}\Big)^{\theta/q}\\
&\qquad=\Big(\e^{\frac\theta{(\theta-1)^2}\ln (2^\gamma\theta^\beta)+\frac1{\theta-1}\ln ((\sigma-\tau)^{-\gamma}2mC_0)}\Big)^{\theta/q}
=\Big(C_1(\sigma-\tau)^{-\theta\gamma/(\theta-1)}\Big)^{1/q},
\end{align*}
where $C_1=C(\theta,\gamma,\beta,m, C_0)>0$, 
we obtain that for any $q\in(0,1)$,
\begin{align}\label{DP9}
\|u\b1_{Q_\tau}\|_{\mL^\infty}\leq \Big(C_1(\sigma-\tau)^{-\theta\gamma/(\theta-1)}\Big)^{1/q}\vee_{i=1}^m\|\b1_{Q_{\sigma}}u\|_{\mL^{q\bbp_i}}.
\end{align}
Let
$$
\delta:=\max_{i,j}p_{ij}\vee 1
\quad  \hbox{and} \quad 
\bbp_i=(p_{i1},\cdots, p_{iN}).
$$
For $p\in(0,1)$ and $q:=p/\delta\in(0,1)$, by H\"older's inequality, we have
\begin{align*}
\|\b1_{Q_{\sigma}}u\|_{\mL^{q\bbp_i}}\leq C\|\b1_{Q_{\sigma}}u\|_{\mL^p}
\end{align*}
where $C=C(\bbp_i,Q_2, \bbp)>0$.
Substituting this into \eqref{DP9}, we obtain \eqref{JK22} with $\gamma_1=\theta\gamma\delta/(\theta-1)$.

\medskip

(ii) Suppose $q_0\in (0, \infty)$. Without loss of generality, we assume $q_0=1$.
Let $\tau_n=\tau+2^{-n}(\sigma-\tau)$. Iterating \eqref{A1A4} with $q=\theta^{-1}$, we obtain
\begin{align} \label{e:2.17}
\sM(\tau,1)&\leq\sM(\tau_n,1)\leq
\Big((\tau_{n-1}-\tau_n)^{-\gamma}mC_0\Big)^{\theta}\sM(\tau_{n-1},\theta^{-1}) \nonumber \\
&\leq\cdots  \leq  \prod_{j=1}^{n}
\Big((\tau_{n-j}-\tau_{n-j+1})^{-\gamma}mC_0 \Big)^{\theta^j}
\sM(\tau_0,\theta^{-n})  .
\end{align}
Noting that
$$
\sum_{j=1}^n\theta^j=\frac{\theta(\theta^n-1)}{\theta-1},
\qquad
 \sum_{j=1}^n(n-j)\theta^j=\frac{\theta^{n+1}-n\theta^2+(n-1)\theta}{(\theta-1)^2},
$$
we have
\begin{align}
& \prod_{j=1}^n\Big((\tau_{n-j}-\tau_{n-j+1})^{-\gamma}mC_0  \Big)^{\theta^j}
=\prod_{j=1}^n\Big(2^{\gamma(n-j+1)}(\sigma-\tau)^{-\gamma}mC_0 \Big)^{\theta^j } \nonumber \\
&\quad =  \exp \left( \sum_{j=1}^n  \theta^j \left((n-j+1)\gamma\ln 2+\ln ((\sigma-\tau)^{-\gamma}2 mC_0) ) \right)  \right)  \nonumber \\
&\quad\leq \e^{C_3(\theta^n-1)\ln (C_4(\sigma-\tau)^{-\gamma})}=(C^{C_3}_4(\sigma-\tau)^{-C_3\gamma})^{\theta^n-1} \label{e:2.18} 
\end{align}
for some $C_3, C_4 >0$. 
 Moreover, if we let 
$$
p_0:=\min_{i,j}p_{ij},\quad  p_1:=\max_{i,j}p_{ij}  \quad \hbox{and} \quad \bbp_i=(p_{i1},\cdots, p_{iN}),
$$
then by H\"older's inequality, we have for some $C>0$,
$$
\|u\b1_{Q_\tau}\|_{\mL^{p_0}}\lesssim_C\vee_{i=1}^m\|u\b1_{Q_{\tau}}\|_{\mL^{\bbp_i}}=\sM(\tau,1),
$$
and
$$
\sM(\tau_0,\theta^{-n})=\vee_{i=1}^m\|u\b1_{Q_{\sigma}}\|_{\mL^{\bbp_i/\theta^n}}\lesssim_C\|u\b1_{Q_{\sigma}}\|_{\mL^{p_1/\theta^n}}.
$$
Combining the above calculations, we have from \eqref{e:2.17}-\eqref{e:2.18} that for some $C_5 >0$,
 $$
\|u\b1_{Q_\tau}\|_{\mL^{p_0}}\leq
(C^{C_3}_5(\sigma-\tau)^{-C_3\gamma})^{\theta^n-1} \|u\b1_{Q_{\sigma}}\|_{\mL^{p_1/\theta^n}}.
$$
This yields \eqref{JK23} by suitable choices of $n$.
 \end{proof}

The following lemma is essentially contained in the proof of 
\cite[p. 126, Lemma 6.21]{Li96}, which 
can be used to verifying condition \eqref{CON1}.

\bl\label{Le24}
Let $(E,\cE,\mu)$ be a measure space and $f:[t_0,t_1]\times E\to \mR $  a $\cB([t_0,t_1])\times\cE$-measurable function. Let $V(t):=\int_E f(t,x)\mu(\dif x)$.
Suppose that $V(t)<\infty$ and there are $c_0, C_0>0$ such that for 
any $t<t'$ in $[t_0,t_1]$, 
$$
c_0\int^{t'}_t\!\!\!\int_E|f(s,x)-V(s)|^2\mu(\dif x)\dif s\leq V(t')- V(t)+C_0.
$$
Then we have for all $\lambda>0$,
\begin{align}\label{GF0}
(\dif t\otimes\mu)\Big\{(t,x)\in[t_0,t_1]\times E: \, f(t,x)-V(t_1 ) 
> \lambda +C_0 
\Big\}\leq 1/(c_0\lambda) ,
\end{align}
and
\begin{align}\label{GF1}
  (\dif t\otimes\mu)\Big\{
 (t,x)\in[t_0,t_1]\times E: \,  f(t, x) -V(t_0) < -( \lambda +C_0) 
  \Big\} \leq 1/ (c_0\lambda) .
\end{align}
\el

\begin{proof}
Let $g(t,x):=f(t,x)-V(t_1)-C_0$ and $P(t):=V(t_1)-V(t)+C_0$. By the assumption we have for all $t\in[t_0,t_1]$,
\begin{align}\label{DS0}
c_0\int^{t_1}_t\!\!\!\int_E|g(s,x)+P(s)|^2\mu(\dif x)\dif s\leq P(t).
\end{align}
For $\lambda>0$ and $t\in[t_0,t_1]$, if we define
$$
\Gamma_\lambda(t):=\{x: g(t,x)>\lambda\},
$$
then since $P(t)\geq 0$ by \eqref{DS0}, we have
$$
g(t,x)+P(t)\geq \lambda+P(t)>0,\ \ x\in\Gamma_\lambda(t).
$$
Thus,
\begin{align*}
\int_E|g(t,x)+P(t)|^2\mu(\dif x)\geq \int_{\Gamma_\lambda(t)}|g(t,x)+P(t)|^2\mu(\dif x)
\geq \mu(\Gamma_\lambda(t))|\lambda+P(t)|^2,
\end{align*}
which together with \eqref{DS0} implies that for  all $t\in[t_0,t_1]$, 
$$
c_0\int^{t_1}_t\mu(\Gamma_\lambda(s))|\lambda+P(s)|^2\dif s\leq P(t).
$$
Let $q(t)$ solve the following ODE:
\begin{align}\label{ODE}
c_0\int^{t_1}_t\mu(\Gamma_\lambda(s))|\lambda+q(s)|^2\dif s=q(t).
\end{align}
Since $q(t_1)=0$ and $\dot q(t)\leq 0$, we have $q(t)\in[0,\infty]$ for $t\in[t_0,t_1]$.
In particular,
$$
P(t)-q(t)\geq\int^{t_1}_t C(s)(P(s)-q(s))\dif s,
$$
where 
$$
C(s)=c_0\mu(\Gamma_\lambda(s))(2\lambda+ P(s)+q(s))\geq 0.
$$
Since $P(t_1)-q(t_1)>0$, solving the above integral inequality yields
$$
P(t)\geq q(t)\geq 0,\ \ t\in[t_0,t_1],
$$
which means that $q(t)$ does not explore.
Thus one can solve ODE \eqref{ODE} to obtain
$$
c_0\int^{t_1}_t\mu(\Gamma_\lambda(s))\dif s=-\int^{t_1}_t\frac{\dif q(s)}{(\lambda+q(s))^2}=\frac1\lambda-\frac1{\lambda+q(t)}
\leq\frac{1}\lambda,
$$
which gives \eqref{GF0}. If we let $\bar g(t,x):=f(t,x)-V(t_0)+C_0$ and $\bar P(t):=V(t)-V(t_0)+C_0$, then we can similarly show \eqref{GF1}.
 \end{proof}

\section{Local upper bound of sub-solutions: Proof of Theorem \ref{Le32}}

Throughout this section,  we assume {\bf (H$_b$)} holds.
The aim of this section is to show Theorem \ref{Le32}, where the key point 
is to verify that any sub-solution of  \eqref{PDE0} lies in the De-Giorgi class introduced in Definition \ref{Def6}.

Let $\rho\in C^\infty_c(Q_1)$ be a symmetric probability density function in $\mR^{1+d}$. For $\eps>0$, we shall use the following mollifiers:
\be \label{e:3.3}
\rho_\eps(t,x)=\eps^{-1-d}\rho(\eps^{-1}t,\eps^{-1}x).
\ee
For a locally integrable function $u$ in $\mR^{1+d}$, we always write
$$
u_\eps(t,x)=(u*\rho_\eps)(t,x)=\int_{\mR^{1+d}}u(s,y)\rho_\eps(t-s,x-y)\dif y\dif s.
$$
Let $u$ be a  super (resp.  sub)-solution of PDE \eqref{PDE0}.  By taking $\varphi=\rho_\eps(t-\cdot,x-\cdot)$ as a test function in \eqref{Def0}, 
we obtain
\begin{align}\label{Mo1}
\begin{split}
\p_t u_\eps\geq \hbox{ (resp. $\leq$)}&\ \sL^\alpha_tu_\eps+(b\cdot\nabla u)*\rho_\eps+f_\eps
=\Delta^{\alpha/2}u_\eps+b\cdot\nabla u_\eps+[\rho_\eps, b\cdot\nabla]u+f_\eps,
\end{split}
\end{align}
where we have used the notation 
 \begin{align}\label{Mo3}
[\rho_\eps, b\cdot\nabla]u:=(b\cdot\nabla u)*\rho_\eps-b\cdot\nabla u_\eps.
\end{align}
Since $b,\nabla b\in\mL^1_{\rm loc}$ and $u\in\mL^\infty_{\rm loc}$, by DiPerna-Lions' commutator estimate (see \cite{Di-Li}), we have
\begin{align}\label{Mo2}
\lim_{\eps\to 0}\| \b1_{Q_R}  [\rho_\eps, b\cdot\nabla]u \|_{\mL^1}=0
\quad \hbox{for any }   R>0.
\end{align}

We first show the following technical result.
\bl\label{Le213}
Suppose that $u$ is a sub-solution of \eqref{PDE0}. Then for any $\kappa\in\mR$, 
$(u-\kappa)^+$ is a sub-solution of \eqref{PDE0} with 
$f\b1_{\{u>\kappa\}}$ in place of $f$.
\el

\begin{proof}
Let $\beta:\mR\to\mR$ be a non-decreasing and convex $C^2$-function. By \eqref{Mo1} and the chain rule, we have
\begin{align*}
\p_t \beta(u_\eps)&\leq (\sL^\alpha_tu_\eps+b\cdot\nabla u_\eps+[\rho_\eps, b\cdot\nabla]u+f_\eps)\beta'(u_\eps).
\end{align*}
Note that
\begin{align*}
\left( \sL^\alpha_t u_\eps \right)\beta'(u_\eps)
 -\sL^\alpha_t \beta(u_\eps)&=\lim_{\delta\downarrow 0}\int_{|x-y|\geq\delta}H_\eps(x,y)K(t,x-y)\dif y,
\end{align*}
where
$$
H_\eps(x,y):=(u_\eps(y)-u_\eps(x))\beta'(u_\eps(x))-(\beta(u_\eps(y))-\beta(u_\eps(x))).
$$
By Taylor's expansion, we have
\begin{align*}
H_\eps(x,y)&=(u_\eps(y)-u_\eps(x))\left[\beta'(u_\eps(x))-\int^1_0\beta'(s (u_\eps(y)-u_\eps(x))+u_\eps(x))\dif s\right]\\
&=-(u_\eps(y)-u_\eps(x))^2\int^1_0\!\!\!\int^1_0s\beta''(ts(u_\eps(y)-u_\eps(x))+u_\eps(x))\dif s\dif t.
\end{align*}
Since $\beta''\geq0$, $H_\eps(x,y) \leq 0$ and so we have 
$$
\p_t \beta(u_\eps)\leq \sL^\alpha_t\beta(u_\eps)+b\cdot\nabla\beta(u_\eps)+
\beta'(u_\eps) \,   [\rho_\eps, b\cdot\nabla]u
+f_\eps\beta'(u_\eps).
$$
By taking limits and using \eqref{Mo2}, it is easy to see that  
for any nonnegative $\varphi\in C^\infty_c(\mR^{1+d})$,
\begin{align}\label{DM0}
-\int_{\mR^{1+d}} \beta(u)\p_t\varphi\leq\int_{\mR^{1+d}} \beta(u)\sL^\alpha_t\varphi
-\int_{\mR^{1+d}}(b\cdot\nabla \varphi+\varphi \,  \div b)
\beta(u)+\int_{\mR^{1+d}} f\beta'(u)\varphi.
\end{align}
Finally, for $\kappa\in\mR$ and $\eps>0$, define
$$
\beta_\eps(r)=\int^r_{-\infty}\int^t_{-\infty}\beta''_\eps(s)\dif s\dif t,
$$
where
$$
\beta''_\eps(s)=\frac{s-\kappa}{\eps^2}\b1_{[\kappa,\kappa+\eps)}(s)+\frac{\kappa+2\eps-s}{\eps^2}\b1_{[\kappa+\eps,\kappa+2\eps]}(s).
$$
Clearly,  $\beta_\eps$ is an increasing convex $C^2$-function with
$$
\lim_{\eps\downarrow 0}\beta_\eps(r)=(r-\kappa)^+,\quad  \lim_{\eps\downarrow 0}\beta'_\eps(r)=\b1_{r>\kappa}.
$$
Using $\beta_\eps$ in place of the $\beta$ in \eqref{DM0} and taking limits $\eps\downarrow 0$, we conclude the proof.
\end{proof}

Now we show the following energy estimate.
\bl\label{Le12}
Let $u$ be a nonnegative sub-solution of PDE \eqref{PDE0}.
For any $\eta\in C^\infty_c(Q_r;[0,1])$, where $r\in[1,2]$,
there is a constant $C=C(d,\alpha)>0$ such that for any  $t\in\mR$ and $\theta\in(1,2)$,
\begin{align}\label{ES5}
\begin{split}
\|u\eta\cI_t\|^2_{\sV_\alpha}
&\lesssim_C\|u^2(\p_s\eta^2-b\cdot\nabla\eta^2-
\eta^2 \,  \div b ) 
\cI_t\|_{\mL^1}+(\|\nabla\eta\|^2_\infty+1)\|u\b1_{Q_{\theta r}}\cI_t\|^2_{\mL^2}\\
&\quad +(\theta-1)^{-2(d+\alpha)}\|\b1_{ \{ u\eta\not=0\}}\cI_t \|^2_{\mL^2}\|{\rm Tail}(u;r)\cI_t\|_{\mL^2_t}^2+\|\<\!\!\!\<f\eta,u\eta\>\!\!\!\>\cI_t\|_{\mL^1_t},
\end{split}
\end{align}
 where $\|\cdot\|_{\sV_\alpha}$ is defined by \eqref{VA},
 and $\<\!\!\!\<\cdot,\cdot\>\!\!\!\>$ stands for the usual dual pair between the distribution  and test functions.
\el

\begin{proof}
Multiplying both sides of \eqref{Mo1} by $u_\eps\eta^2$ and integrating in $\mR^d$, by \eqref{SY1} we obtain
\begin{align}\label{SA1}
\frac12\int_{\mR^{d}}\p_t(u_\eps\eta)^2-\frac12\int_{\mR^{d}}u^2_\eps\p_t\eta^2
&\leq-\cE_t^{(\alpha)}(u_\eps, u_\eps\eta^2)+\int_{\mR^{d}}(b\cdot\nabla u_\eps+[\rho_\eps, b\cdot\nabla]u+f_\eps)u_\eps\eta^2,
\end{align}
where $\cE_t^{(\alpha)}(u_\eps, u_\eps\eta^2)$ is defined by \eqref{SY91}, $[\rho_\eps, b\cdot\nabla]u$ is defined by \eqref{Mo3}.
Noting that for any $a,b,\beta,\gamma\in\mR$,
\begin{align*}
(a-b)(a \beta^2-b\gamma^2)&=(a\beta-b\gamma)^2-ab(\beta-\gamma)^2,
\end{align*}
by definition \eqref{SY91}, we have
\begin{align*}
-\cE_t^{(\alpha)}(u_\eps, u_\eps\eta^2)&
=\frac12
\int_{\mR^d}\int_{\mR^d}(u_\eps(y)-u_\eps(x))((u_\eps\eta^2)(x)-(u_\eps\eta^2)(y))K(t,x-y)\dif x\dif y\\
&=-\frac12
\int_{\mR^d}\int_{\mR^d}(u_\eps(x)\eta(x)-u_\eps(y)\eta(y))^2K(t,x-y)\dif x\dif y\\
&\quad+\frac12
\int_{\mR^d}\int_{\mR^d}u_\eps(x)u_\eps(y)(\eta(x)-\eta(y))^2K(t,x-y)\dif x\dif y\\
&\stackrel{\eqref{KK0}}{\leq}-\kappa_0\|\Delta^{\alpha/4}(u_\eps\eta)\|^2_{\mL^2_x}+I.
\end{align*}
Here we dropped the time variable $t$ in $u_\eps$ and $\eta$.
Since $\eta$ has support in $Q_r\subset Q_{\theta r}$, where $\theta\in(1,2)$, we have
\begin{align*}
I &=\frac12
\int_{B_{\theta r}}\int_{B_{\theta  r}}u_\eps(x)u_\eps(y)(\eta(x)-\eta(y))^2K(t,x-y)\dif x\dif y\\
&\quad+
\int_{B_{\theta r}}\int_{B^c_{\theta r}}u_\eps(x)u_\eps(y)(\eta(x)-\eta(y))^2K(t,x-y)\dif x\dif y
=: I_1+I_2.
\end{align*}
For $I_1$, by $ab\leq(a^2+b^2)/2$ and the symmetry of $K(t,y)=K(t,-y)$, we have
\begin{align*}
I_1&\leq
\frac14
\int_{B_{\theta r}}\int_{B_{\theta r}}(|u_\eps(x)|^2+|u_\eps(y)|^2)(\eta(x)-\eta(y))^2K(t,x-y)\dif x\dif y\\
&=\frac12
\int_{B_{\theta r}}\int_{B_{\theta r}}|u_\eps(x)|^2(\eta(x)-\eta(y))^2K(t,x-y)\dif x\dif y\\
&\stackrel{\eqref{KK0}}{\leq}\frac{\kappa_1}2
\left(\sup_x\int_{\mR^d}\frac{(\eta(x)-\eta(y))^2}{|x-y|^{d+\alpha}}\dif y\right)\|u_\eps\b1_{Q_{\theta r}}\|^2_{\mL^2_x}\\
&\lesssim(\|\nabla\eta\|^2_\infty+\|\eta\|^2_\infty)\|u_\eps\b1_{Q_{\theta r}}\|^2_{\mL^2_x}.
\end{align*}
For $I_2$, noting that $\theta\in(1,2)$ and for $|x|\leq r$ and $|y|\geq \theta r$, 
$$
|x-y|^{-d-\alpha}\leq (1-\tfrac1\theta)^{-d-\alpha} |y|^{-d-\alpha},
$$
since $\eta\in C^\infty_c(Q_r;[0,1])$ and $u_\eps\geq 0$, we have
\begin{align*}
I_2&\stackrel{\eqref{KK0}}{\leq}
\kappa_1\int_{B_r}\int_{B^c_{\theta r}}\frac{u_\eps(x)u_\eps(y)\eta(x)^2}{|x-y|^{d+\alpha}}\dif x\dif y\\
&\leq \kappa_1
(1-\tfrac1\theta)^{-d-\alpha}\int_{B_r}(u_\eps\eta)(x)\dif x\int_{B^c_{\theta r}}\frac{u_\eps(y)}{|y|^{d+\alpha}}\dif y\\
&=\kappa_1
\theta^{d+\alpha}(\theta-1)^{-d-\alpha}\|u_\eps\eta\|_{\mL^1_x}{\rm Tail}(u_\eps, \theta r)\\
&\leq \kappa_1 2^{d+\alpha}
(\theta-1)^{-d-\alpha}\|u_\eps\eta\|_{\mL^1_x}{\rm Tail}(u_\eps, r).
\end{align*}
Combining the above calculations, we obtain that for some $C=C(d,\alpha)>0$,
\begin{align*}
\cE_t^{(\alpha)}(u_\eps, u_\eps\eta^2)&\geq \|\Delta^{\alpha/4}(u_\eps\eta)\|^2_{\mL^2_x}
-C(\|\nabla\eta\|^2_\infty+1)\|u_\eps\b1_{Q_{\theta r}}\|^2_{\mL^2_x}\\ 
&\quad-\kappa_12^{d+\alpha}
(\theta-1)^{-d-\alpha}\|u_\eps\eta\|_{\mL^1_x}{\rm Tail}(u_\eps; r).
\end{align*}
We clearly have
\begin{align*}
\int_{\mR^{d}}(b\cdot\nabla u_\eps)(u_\eps\eta^2)
&=\frac12\int_{\mR^{d}}(b\cdot\nabla u^2_\eps)\eta^2
=-\frac12\int_{\mR^{d}}u^2_\eps(b\cdot\nabla\eta^2+
\eta^2 \,  \div b ) . 
\end{align*}
Substituting the above two estimates into \eqref{SA1} and then integrating both sides  from $-\infty$ to $t$, 
we obtain that for some $C=C(d,\alpha)>0$,
\begin{align*}
&\frac12\|(u_\eps\eta)(t)\|_{\mL^2_x}^2+\int^t_{-\infty}\|\Delta^{\alpha/4}(u_\eps\eta)\|^2_{\mL^2_x}\\
&\leq \frac12\int^t_{-\infty}\!\int_{\mR^{d}}u^2_\eps\Big(\p_s\eta^2-b\cdot\nabla\eta^2
- \eta^2\, \div b 
\Big)+\int^t_{-\infty}\!\int_{\mR^{d}}([\rho_\eps, b\cdot\nabla]u+f_\eps) u_\eps\eta^2\\
&\quad+C(\|\nabla\eta\|^2_\infty+1)\int^t_{-\infty}\|u_\eps\b1_{Q_{\theta r}}\|^2_{\mL^2_x}+C(\theta-1)^{-d-\alpha}\int^t_{-\infty}\|u_\eps\eta\|_{\mL^1_x}{\rm Tail}(u_\eps, r).
\end{align*}
Taking limits $\eps\to 0$ and by Fatou's lemma and \eqref{Mo2}, we obtain that for some $C=C(d,\alpha)>0$,
\begin{align*}
&\frac12\|(u\eta)(t)\|_{\mL^2_x}^2+\int^t_{-\infty}\|\Delta^{\alpha/4}(u\eta)\|^2_{\mL^2_x}
\leq\frac12 \int^t_{-\infty}\|u^2(\p_s\eta^2-b\cdot\nabla\eta^2- \eta^2 \, \div b 
)\|_{\mL^1_x}+\int^t_{-\infty}\<\!\!\!\<f\eta,u\eta\>\!\!\!\>\\
&\qquad\quad+C(\|\nabla\eta\|^2_\infty+1)\int^t_{-\infty}\|u\b1_{Q_{\theta r}}\|^2_{\mL^2_x}+C(\theta-1)^{-d-\alpha}\int^t_{-\infty}\|u\eta\|_{\mL^1_x}{\rm Tail}(u, r).
\end{align*}
Taking
  supremum in the time variable from $-\infty$ to $t$ on  both sides, 
we get by H\"older's inequality,
\begin{align*}
&\frac12\|u\eta\cI_t\|_{\mL^\infty_t\mL^2_x}^2+\|\Delta^{\alpha/4}(u\eta)\cI_t\|^2_{\mL^2_t\mL^2_x}
\leq\frac12\|u^2(\p_s\eta^2-b\cdot\nabla\eta^2-
  \eta^2 \, \div b 
)\cI_t\|_{\mL^1}+\|\<\!\!\!\<f\eta,u\eta\>\!\!\!\>\cI_t\|_{\mL^1_t}\\
&\quad+C(\|\nabla\eta\|^2_\infty+1)\|u\b1_{Q_{\theta r}}\cI_t\|^2_{\mL^2}+C(\theta-1)^{-d-\alpha}
\|u\eta\cI_t\|_{\mL^\infty_t\mL^2_x}\|\b1_{\{u\eta\not=0\}}\cI_t\|_{\mL^2_t\mL^2_x}\|{\rm Tail}(u_\eps, r)\|_{\mL^2_t},
\end{align*}
which in turn implies \eqref{ES5} by Young's inequality.
\end{proof}

\bl\label{Le31}
Let $\beta\in[0, \alpha/ 2)$ and  
$(q_0,\bbp_0)\in\mI_\alpha^\beta$.
For any $\eps\in(0,1)$, there exists a constant $C_\eps=C(\eps,\alpha,\beta,d,\bbp_0,q_0)>0$ such that
\begin{align}\label{38}
\|\<\!\!\!\<f\eta,u\eta\>\!\!\!\>\cI_t\|_{\mL^1_t}\leq 
\eps\|u\eta\cI_t\|^2_{\sV_\alpha}+C_\eps\|f\eta\cI_t\|^2_{\mL^{ q_0}_t(\bH^{-\beta}_{\bbp_0})}\|\b1_{\{u\eta\not=0\}}\cI_t\|^2_{\mL^{\nu_0}_t\mL^{\bbr_0}_x},
\end{align}
where $(\nu_0,\bbr_0)\in(2,\infty)^{1+d}$ satisfies $\tfrac\alpha{\nu_0}+|\tfrac{1}{\bbr_0}|>\tfrac d2$.
\el
\begin{proof}
Let $ (\nu_1,\bbr_1)\in(1,\infty)^{1+d}$ be defined by
$$
\tfrac1{ q_0}+\tfrac1{\nu_1}=1,\ \ \tfrac1{\bbp_0}+\tfrac1{\bbr_1}={\bf 1}.
$$
Since the dual space of $\bH^\beta_{\bbr_1}$ is $\bH^{-\beta}_{\bbp_0}$, by H\"older's inequality in the time variable, we have
\begin{align}\label{DD6}
\|\<\!\!\!\<f\eta,u\eta\>\!\!\!\>\cI_t\|_{\mL^1_t}\leq
\|f\eta\cI_t\|_{\mL^{ q_0}_t(\bH^{-\beta}_{\bbp_0})}\|u\eta\cI_t\|_{\mL^{ \nu_1}_t(\bH^{\beta}_{\bbr_1})}.
\end{align}
Since $\tfrac\alpha{q_0}+|\tfrac{1}{\bbp_0}|<\alpha-\beta$, one can choose $\theta\in(\frac{2\beta}\alpha,1)$ so that
\begin{align}\label{DD91}
\tfrac\alpha{q_0}+|\tfrac{1}{\bbp_0}|<\alpha-\tfrac{\theta\alpha}2.
\end{align}
Let $\nu_2\in(0,\infty)$ and $\bbr_2\in(1,\infty)^d$ be defined by
$$
\tfrac1{ \nu_1}=\tfrac{1-\theta}{ \nu_2}+\tfrac{\theta}{2},\ \ \tfrac1{\bbr_1}=\tfrac{1-\theta}{\bbr_2}+\tfrac{\theta}{\bb2}
$$
so that
$$
|\tfrac{1}{\bbr_1}|+(1-\theta)(0-|\tfrac{1}{\bbr_2}|)+\theta(\tfrac\alpha 2-\tfrac{d}{2})=\tfrac{\theta\alpha}2>\beta.
$$
By Gagliado-Nirenberge's inequality \eqref{Sob}, we have
\begin{align*}
\|u\eta\cI_t\|_{\mL^{ \nu_1}_t(\bH^{\beta}_{\bbr_1})}
\lesssim \|u\eta\cI_t\|^{1-\theta}_{\mL^{ \nu_2}_t\mL^{\bbr_2}_x}\|u\eta\cI_t\|^{\theta}_{\mL^2_t(\bH^{\alpha/2}_{2})}.
\end{align*}
Let $\nu_0:=\frac{2-\theta}{1-\theta}\, \nu_2$ and $\bbr_0:=\frac{2-\theta}{1-\theta}\, \bbr_2$. By definition and \eqref{DD91},
it is easy to see that
$$
(\nu_0,\bbr_0)\in(2,\infty)^{1+d}
\quad \hbox{and} \quad 
\tfrac\alpha{\nu_0}+|\tfrac{1}{\bbr_0}|>\tfrac d2.
$$
Thus, by H\"older's inequality and Lemma \ref{Le22}, we have
\begin{align*}
\|u\eta\cI_t\|^{1-\theta}_{\mL^{ \nu_2}_t\mL^{\bbr_2}_x}
&\leq \|\b1_{\{u\eta\not=0\}}\cI_t\|^{1-\theta}_{\mL^{ \nu_0(1-\theta)}_t\mL^{\bbr_0(1-\theta)}_x}
\|u\eta\cI_t\|^{1-\theta}_{\mL^{\nu_0}_t\mL^{\bbr_0}_x}
\lesssim\|\b1_{\{u\eta\not=0\}}\cI_t\|_{\mL^{ \nu_0}_t\mL^{\bbr_0}_x}
\|u\eta\cI_t\|^{1-\theta}_{\sV_\alpha}.
\end{align*}
Hence,
$$
\|u\eta\cI_t\|_{\mL^{ \nu_1}_t(\bH^{\beta}_{\bbr_1})}
\lesssim \|\b1_{\{u\eta\not=0\}}\cI_t\|_{\mL^{\nu_0}_t\mL^{\bbr_0}_x}\|u\eta\cI_t\|_{\sV_\alpha}.
$$
Substituting this into \eqref{DD6} and by Young's inequality, we obtain \eqref{38}.
\end{proof}

Now we are in a position to give the 

\begin{proof}[Proof of Theorem \ref{Le32}]
Let $\sI\subset\mR^{1+d}$ be an open index subset defined by
\begin{align}\label{II}
\sI:=\Big\{(\nu,\bbr)\in(2,\infty)^{1+d}: \tfrac\alpha{\nu}+|\tfrac{1}\bbr|>\tfrac d2\Big\}.
\end{align}
For fixed $t\in[-2,2]$, we define
$$
\wt Q^t_\tau:=Q_\tau\cap((-\infty,t)\times\mR^{d}).
$$
By Theorem \ref{TH20}, it suffices to show that
$$
u\in\cD\cG^+_{\sI}(\wt Q^t_\cdot)  \quad \mbox{with $\cA=\|\chi_2f\cI_t\|_{\mL^{ q_0}_t(\bH^{-\beta}_{\bbp_0})}+\|{\rm Tail}(u^+;1)\cI_t\|_{\mL^2_t}$.}
$$ 
More precisely, we want to show that for any 
$(\nu,\bbr)\in\sI$, there are $\gamma,C>0$ only depending on $\kappa_0, \alpha, \beta, d$, $\nu,\bbr$ and $q_i,\bbp_i, i=0,1,2$ such that
for any $1\leq\tau<\sigma\leq 2$, $|t|\leq2$ and $\kappa\geq 0$,
\begin{align}\label{DJ180}
\begin{split}
&(\sigma-\tau)^{\gamma}\|\b1_{Q_\tau}(u-\kappa)^+\cI_t\|_{\mL^\nu_t\mL^\bbr_x}
\lesssim_C\|\b1_{Q_\sigma}(u-\kappa)^+\cI_t\|_{\mL^2}+\sum_{i=1,2}\|\b1_{Q_\sigma}(u-\kappa)^+\cI_t\|_{\mL^{\nu_i}_t\mL^{\bbr_i}_x}\\
&\qquad
+\|\b1_{\{u>\kappa\}\cap Q_\sigma}\cI_t \|_{\mL^2}\|{\rm Tail}(u^+;1)\cI_t\|_{\mL^2_t}
+\|\chi_2f\cI_t\|_{\mL^{ q_0}_t(\bH^{-\beta}_{\bbp_0})}\|\b1_{\{u>\kappa\}\cap Q_\sigma}\cI_t\|_{\mL^{\nu_0}_t\mL^{\bbr_0}_x},
\end{split}
\end{align}
where $(\nu_i,\bbr_i)\in\sI$, $i=0,1,2$ are determined by $( q_i,\bbp_i)$.
By Lemma \ref{Le213}, we may assume $\kappa=0$ and $u$ is a nonnegative sub-solution. 

Let $1\leq \tau<\sigma\leq 2$ and 
$$\eta\in C^\infty_c(Q_{(\tau+\sigma)/2};[0,1]),\ \ \eta|_{Q_\tau}\equiv 1$$   
satisfy that for some universal constant $C>0$,
\begin{align}\label{DX877}
\|\p_t\eta\|_\infty+\|\nabla\eta\|_\infty\leq C(\sigma-\tau)^{-1}.
\end{align}
Fix $t\in[-2,2]$. By \eqref{ES5} with $\theta:=2\sigma/(\tau+\sigma)\in(1,2)$ and \eqref{DX877},  we have
\begin{align}\label{ES55}
\begin{split}
\|u\eta\cI_t\|^2_{\sV_\alpha}
&\lesssim\|u^2(b\cdot\nabla\eta^2+ \eta^2 \, \div b 
)\cI_t\|_{\mL^1}+(\sigma-\tau)^{-2}\|u\b1_{Q_\sigma}\cI_t\|^2_{\mL^2}\\
&\quad +(\sigma-\tau)^{-2(d+\alpha)}\|\b1_{\{u\eta\not=0\}}\cI_t \|^2_{\mL^2}\|{\rm Tail}(u, 1)\cI_t\|_{\mL^2_t}^2+\|\<\!\!\!\<f\eta,u\eta\>\!\!\!\>\cI_t\|_{\mL^1_t}.
\end{split}
\end{align}
By H\"older's inequality and Young's inequality, we have
\begin{align*}
\|u^2(b\cdot\nabla\eta^2+ \eta^2 \, \div b
)\cI_t\|_{\mL^1}
&\leq2\|u^2|b|\b1_{Q_\sigma}\cI_t\|_{\mL^1}\|\nabla\eta\|_\infty+\|u^2|\div b|\b1_{Q_\sigma}\cI_t\|_{\mL^1}\\
&\leq 2\|u\b1_{Q_\sigma}\cI_t\|^2_{\mL^{\nu_1}_t\mL^{\bbr_1}_x}\|b\b1_{Q_\sigma}\|_{\mL^{q_1}_t\mL^{\bbp_1}_x}\|\nabla \eta\|_\infty\\
&\quad+\|u\b1_{Q_\sigma}\cI_t\|^2_{\mL^{\nu_2}_t\mL^{\bbr_2}_x}\|\div b\b1_{Q_\sigma}\|_{\mL^{q_2}_t\mL^{\bbp_2}_x},
\end{align*}
where $(\nu_i,\bbr_i)\in\sI$, $i=1,2$ are defined by
$$
\tfrac1{q_i}+\tfrac2{\nu_i}=1,\ \ \tfrac1{\bbp_i}+\tfrac2{\bbr_i}={\bf 1},\ i=1,2.
$$
Moreover, by Lemma \ref{Le31} we have
$$
\|\<\!\!\!\<f\eta,u\eta\>\!\!\!\>\cI_t\|_{\mL^1_t}\leq
\eps\|u\eta\cI_t\|^2_{\sV_\alpha}+C_\eps\|f\eta\cI_t\|^2_{\mL^{ q_0}_t(\bH^{-\beta}_{\bbp_0})}\|\b1_{\{u\eta\not=0\}}\cI_t\|^2_{\mL^{\nu_0}_t\mL^{\bbr_0}_x}.
$$
Combining the above estimates and letting $\eps$ be small enough, we obtain
\begin{align*}
\|u\eta\cI_t\|^2_{\sV_\alpha}
& \lesssim  (\sigma-\tau)^{-2}\|u\b1_{Q_\sigma}\cI_t\|^2_{\mL^2}+(\sigma-\tau)^{-2}\sum_{i=1,2}\|u\b1_{Q_\sigma}\cI_t\|^2_{\mL^{\nu_i}_t\mL^{\bbr_i}_x}\\
&\quad+(\sigma-\tau)^{-2(d+\alpha)}\|\b1_{\{u\eta\not=0\}}\|^2_{\mL^2} \|{\rm Tail}(u,1)\cI_t\|_{\mL^2_t}^2\\
&\quad+\|f\eta\cI_t\|^2_{\mL^{ q_0}_t(\bH^{-\beta}_{\bbp_0})}\|\b1_{\{u\eta\not=0\}}\cI_t\|^2_{\mL^{\nu_0}_t\mL^{\bbr_0}_x}.
\end{align*}
Since for any $(\nu,\bbr)\in\sI$,  by Lemma \ref{Le22}, we have
$$
\|u\b1_{Q_\tau}\cI_t\|_{\mL^\nu_t\mL^\bbr_x}\leq\|u\eta\cI_t\|_{\mL^\nu_t\mL^\bbr_x}\leq C\|u\eta\cI_t\|_{\sV_\alpha}.
$$
Also, since $\eta\chi_2=\eta$, we have
$$
\|f\eta\cI_t\|_{\mL^{ q_0}_t(\bH^{-\beta}_{\bbp_0})}\leq C\|f\chi_2\cI_t\|_{\mL^{ q_0}_t(\bH^{-\beta}_{\bbp_0})}\|\eta\|_{C^1_b}
\leq C(\sigma-\tau)^{-1}\|f\chi_2\cI_t\|_{\mL^{ q_0}_t(\bH^{-\beta}_{\bbp_0})}.
$$
Thus we get \eqref{DJ180}, and by Theorem \ref{TH20}, for any $p>0$ there is a $C>0$ such that
$$
\|u\b1_{Q_1}\cI_t\|_{\mL^\infty}\lesssim_C \|u\b1_{Q_2}\cI_t\|_{\mL^p}+\|{\rm Tail}(u^+; 1)\cI_t\|_{\mL^2_t}+\|\chi_2f\cI_t\|_{\mL^{ q_0}_t(\bH^{-\beta}_{\bbp_0})}.
$$
Moreover, by applying Lemma \ref{Le12} with $\eta=\chi_1$, as above and using the above $L^\infty$-estimate, we obtain the desired estimate.
The proof is complete.
\end{proof}

\section{Weak Harnack inequality for super-solutions}

Throughout this section we assume {\bf (H$'_b$)}.
The aim of this section is to show Theorem \ref{Th41}.
In the following proof, we always assume that $u$ is  a super-solution 
of \eqref{PDE0} in $Q_4 $  with $u\b1_{Q_4 }\geq 0$.      
Suppose that 
\begin{align}\label{KK1}
\kappa\geq \|f\b1_{Q_2}\|_{\mL^{q_0}_t\mL^{\bbp_0}_x}+\|{\rm Tail}(u^-,1)\b1_{[-2,2]}\|_{\mL^1_t}.
\end{align}
Recalling $u_\eps=u*\rho_\eps$, we define
\begin{align}\label{AA1}
w_\eps:=u_\eps+\kappa
\quad \hbox{and} \quad   w:=u+\kappa.
\end{align}
Clearly, $\lim_{\eps \to 0} w_\eps =w$.
By \eqref{Mo1}, we have
\begin{align}\label{DA4}
\p_t w_\eps \geq \sL^\alpha_tw_\eps+(b\cdot\nabla w)*\rho_\eps+f_\eps
=\sL^\alpha_tw_\eps+b\cdot\nabla w_\eps+[\rho_\eps, b\cdot\nabla]u+f_\eps,
\end{align}
where $[\rho_\eps, b\cdot\nabla]u$ is defined by \eqref{Mo3}. For fixed $q>0$, let 
$$
\bar q:=(q+1)\vee 2.
$$
Let $\eta:\mR^{1+d}\to [0,1]$ be a $C^1$-function with support contained in $Q_\sigma$, where $\sigma\in(0,4]$.
Multiplying both sides of \eqref{DA4} by $w_\eps^{-q}\eta^{\bar q}$ and integrating over $\mR^d$, by \eqref{SY1} and \eqref{AA12}, 
we have for $t\in (-4, 4)$ and $\eps >0$ sufficiently small, 
\begin{align}\label{DZ1}
\begin{split}
\int_{\mR^d} \p_t w_\eps w_\eps^{-q}\eta^{\bar q}\dif x  &
\geq -\cE_t^{(\alpha)}(w_\eps,w_\eps^{-q}\eta^{\bar q})
+\int_{\mR^d} (b\cdot\nabla w_\eps+[\rho_\eps, b\cdot\nabla]u+f_\eps)w_\eps^{-q}\eta^{\bar q}\dif x  \\ 
&=\frac12
\int_{B_\sigma}\!\int_{B_\sigma}(w_\eps(y)-w_\eps(x))((w_\eps^{-q}\eta^{\bar q})(x)-(w_\eps^{-q}\eta^{\bar q})(y))K(t,x-y)\dif x\dif y \\
&\quad+
\int_{B_\sigma}\left(\int_{B^c_\sigma}(w_\eps(y)-w_\eps(x))K(t,x-y)\dif y\right)(w_\eps^{-q}\eta^{\bar q})(x)\dif x\\
&\quad+\int_{\mR^d} (b\cdot\nabla w_\eps)w_\eps^{-q}\eta^{\bar q}
\dif x +\int_{\mR^d} ([\rho_\eps, b\cdot\nabla]u+f_\eps)w_\eps^{-q}\eta^{\bar q} \dif x \\
 &=:\cJ^\eps_1 (t) +\cJ^\eps_2 (t) +\cJ^\eps_3 (t) +\cJ^\eps_4 (t).
 \end{split}
\end{align}
Below, we shall consider three cases for the calculations about \eqref{DZ1}: 
$q=1$, $q\in (1, \infty)$ and $q\in (0, 1)$.

\medskip

\begin{itemize}
\item In Lemma \ref{Le42} below, we take $q=1$ to derive an estimate of $\log w$.
\item In Lemma \ref{Le45} below, we take $q\in(1,\infty)$ to get an upper bound estimate of $w^{-1}$.
\item In  Lemma \ref{Le46} below, we take $q\in(0,1)$ to derive an integrability estimate of $w^{1-q}$.
\end{itemize}
From these estimates, we can establish Theorem \ref{Th41}.

\medskip

Before presenting 
the proofs, we first recall some well-known facts.
We need the following  
 algebraic  inequalities which are proved  in \cite[Lemma 3.3]{FK13}.

\medskip

\bl\label{Le26}
\begin{enumerate}[\rm (i)]
\item 
For all $a,b>0$ and $\alpha,\beta\geq 0$,
\begin{align}\label{EQ9}
(b-a)(\alpha^2a^{-1}-\beta^2b^{-1})\geq\Big(\log\tfrac{a}{\alpha}-\log\tfrac b\beta\Big)^2\alpha\beta-(\alpha-\beta)^2.
\end{align}
\item Let $q>1$ and $\ell(q):=\max\{4,\frac{6q-5}{2}\}$. It holds that for all $a,b>0$ and $\alpha,\beta\geq 0$,
\begin{align*}
(b-a) \left(\alpha^{q+1}a^{-q}-\beta^{q+1}b^{-q} \right)
&\geq\tfrac{1}{q-1}\left[ \Big(\tfrac{\beta}{b}\Big)^{\frac{q-1}2}-\Big(\tfrac{\alpha}{a}\Big)^{\frac{q-1}2}\right]^2
-\ell(q)(\beta-\alpha)^2\left[ \Big(\tfrac{\beta}{b}\Big)^{q-1}+\Big(\tfrac{\alpha}{a}\Big)^{q-1}\right].
\end{align*}
\item Let $q\in(0,1)$, $\ell_1(q):=\frac{2q}{3(1-q)}$ and $\ell_2(q):=\frac{4q}{1-q}+\frac{9}{q}$. It holds that for all $a,b>0$ and $\alpha,\beta\geq 0$,
$$
(b-a) \left(\alpha^2a^{-q}-\beta^2b^{-q} \right)
\geq \ell_1(q)\left( \beta b^{\frac{1-q}2}-\alpha a^{\frac{1-q}2}\right)^2
-\ell_2(q)(\beta-\alpha)^2\left( b^{1-q}+a^{1-q}\right).
$$
\end{enumerate}
\el

\medskip

\bl\label{Le42}
Under \eqref{KK1}, it holds that for some $a=a(w(0,\cdot))\in\mR$ and all $\lambda>0$,
\begin{align}\label{GF00}
(\dif t\otimes\dif x)\big\{(t,x)\in Q^-_2(0): \log w(t,x)-a-C_0>\lambda\big\}\leq C_1/\lambda ,
\end{align}
and
\begin{align}\label{GF01}
(\dif t\otimes\dif x)\big\{(t,x)\in Q^+_2(0): a-C_0-\log w(t,x)>\lambda\big\}\leq C_1/\lambda ,
\end{align}
where $C_0>0$ depends only on $d,\alpha$ and $\|\div b\b1_{Q_4}\|_{\mL^2}$, $\|b\b1_{Q_4}\|_{\mL^2}$ and $C_1=C_1(d,\alpha)>0$.
\el
\begin{proof}
Let $\eta(x):=((3-|x|)\wedge 1)\vee0$. Note that
\begin{align}\label{KA33}
\eta|_{B_2}=1,\quad \eta|_{B^c_3}=0.
\end{align}
Define
$$
v_\eps(t,x):=\log\big(\tfrac{w_\eps(t,x)}{\eta(x)}\big) \b1_{\{ | x|<3\} },
\quad  V_\eps(t):=\int_{\mR^d}v_\eps(t,x)\mu(\dif x),
$$
where $\mu$ is a probability measure given by
$$
\mu(\dif x):=\eta(x)^2\dif x/\|\eta\|^2_{\mL^2_x}.
$$
Taking $q=1$ and $\sigma=4$ in \eqref{DZ1}, we estimate each term of \eqref{DZ1} as follows.

\begin{enumerate}[$\bullet$]
\item Since $\eta$ does not depend on $t$, we have by the chain rule,  
\begin{align}\label{KA2}
\int_{\mR^d} \p_t w_\eps w_\eps^{-1}\eta^{2}=\int_{\mR^d} \p_t \log w_\eps  \eta^2=\int_{\mR^d} \p_t \log \big(\tfrac{w_\eps}{\eta}\big)  \eta^2
=\p_t\int_{\mR^d}v_\eps  \eta^2=\|\eta\|^2_{\mL^2_x}\p_tV_\eps.
\end{align}

\item For 
$\cJ^\eps_1 (t)$,
by \eqref{EQ9}, \eqref{KK0} and the weighted Poincare's inequality \cite[Proposition 2.3]{FK13},
\begin{align}
\cJ^\eps_1 (t) 
&=\frac12\int_{B_4}\!\int_{B_4}(w_\eps(y)-w_\eps(x))((w_\eps^{-1}\eta^2)(x)-(w_\eps^{-1}\eta^2)(y))K(t,x-y)\dif x\dif y\no\\
&\geq\frac12 \int_{B_4}\!\int_{B_4}
\left(\log\frac{w_\eps(x)}{\eta(x)}-\log\frac{w_\eps(y)}{\eta(y)}\right)^2\eta(x)\eta(y)K(t,x-y)\dif x\dif y
 \no\\
&\quad-\frac12\int_{B_4}\!\int_{B_4}(\eta(x)-\eta(y))^2K(t,x-y)\dif x\dif y\no\\
&\geq c_0\|\eta\|^2_{\mL^2}
\int_{B_4}
\left(v_\eps(x)-V_\eps\right)^2\mu(\dif x)-C_2,\label{KA3}
\end{align}
where $c_0, C_2>0$ only depend on $d,\alpha, \kappa_0,\kappa_1$. 

\item For 
$\cJ^\eps_2 (t)$, 
noting that for $|x|\leq 3$ and $|y|\geq 4$,
$$
|x-y|^{-d-\alpha}\leq 4^{d+\alpha}|y|^{-d-\alpha},
$$
and $\eta(x)=0$ for $|x|>3$,  we have for some $C_3=C_3(d,\alpha)>0$,
\begin{align*}
 - \cJ^\eps_2 (t) 
 &=\int_{B_{3}}\left(\int_{B^c_4}(w_\eps(x)-w_\eps(y))K(t,x-y)\dif y\right)(w_\eps^{-1}\eta^2)(x)\dif x\\
&\stackrel{\eqref{KK0}}{\leq} \int_{B_{3}}\int_{B^c_4}\left(\eta^2(x)+\frac{1}\kappa w^-_\eps(y)\right)\frac{\kappa_14^{d+\alpha}\dif y}{|y|^{d+\alpha}}\dif x\\
&\leq C_3\Big(1+{\rm Tail}(w^-_\eps,4)/\kappa\Big).
\end{align*}

\item For $\cJ^\eps_3 (t)$,   we have
by the chain rule and integration by parts,  
\begin{align*}
 \cJ^\eps_3 (t)
 &=\int_{\mR^d} (b\cdot\nabla \log w_\eps)\eta^2
=\int_{\mR^d} (b\cdot\nabla(v_\eps+ \log\eta))\eta^2\\
&=-\int_{\mR^d} v_\eps\div (b\eta^2)+\int_{\mR^d} (b\cdot\nabla\log \eta)\eta^2\\
&=\int_{\mR^d} (V_\eps-v_\eps)\div (b\eta^2)+\int_{\mR^d} (b\cdot\nabla \eta)\eta\\
&=\int_{B_4} (V_\eps-v_\eps)\eta(
\eta \, \div b 
+2b\cdot\nabla\eta)+\int_{B_4} (b\cdot\nabla \eta)\eta.
\end{align*}
By Young's inequality,  for some $C_4=C_4(d,\alpha)\geq 1$,
\begin{align}
\left|\cJ^\eps_3 (t) \right|
&\leq\tfrac{c_0}{2}\int_{B_4} (v_\eps-V_\eps)^2\eta^2+C_4\int_{B_4}(|\div b|^2+|b|^2+|b|)\no\\
&=\tfrac{c_0\|\eta\|^2_{\mL^2_x}}{2}\int_{B_4} (v_\eps(x)-V_\eps)^2\mu(\dif x)+C_4\int_{B_4}(|\div b|^2+|b|^2+|b|).\label{KA4}
\end{align}
\end{enumerate}
For $-2\leq t_0<t_1\leq 2$, substituting \eqref{KA2}-\eqref{KA4} into \eqref{DZ1}
with $q=1$ and integrating from $t_0$ to $t_1$, we get
\begin{align*}
&\tfrac{c_0\|\eta\|^2_{\mL^2_x}}{2}\int^{t_1}_{t_0}\!\!\!\int_{B_4}\left(v_\eps(s,x)-V_\eps(s)\right)^2\mu(\dif x)\dif s
\leq\|\eta\|^2_{\mL^2_x}\(V_\eps(t_1)-V_\eps(t_0)\)\\
&\qquad\qquad+C_3\int^{t_1}_{t_0}\Big(1+\frac{{\rm Tail}(w^-_\eps,4)}{\kappa}\Big)+C_4\int^{t_1}_{t_0}\!\!\!\int_{B_4}(|\div b|^2+|b|^2+|b|)\\
&\qquad\qquad\qquad+\int^{t_1}_{t_0}\!\!\!\int_{B_4}\frac{|[\rho_\eps, b\cdot\nabla]u+f_\eps|}{\kappa}\eta^2.
\end{align*}
Letting $\eps\to 0$ and by Fatou's lemma and \eqref{Mo2}, we arrive at
\begin{align*}
&\tfrac{c_0\|\eta\|^2_{\mL^2_x}}{2}\int^{t_1}_{t_0}\!\!\!\int_{B_4}\left(v(s,x)-V(s)\right)^2\mu(\dif x)\dif s
\leq\|\eta\|^2_{\mL^2_x}\(V(t_1)-V(t_0)\)\\
&\qquad\qquad+C_3\int^{t_1}_{t_0}\Big(1+\frac{{\rm Tail}(w^-,4)}{\kappa}\Big)+C_4\int^{t_1}_{t_0}\!\!\!\int_{B_4}\Big(|\div b|^2+|b|^2+|b|+\frac{|f|} \kappa\Big),
\end{align*}
where $V(t):=\int_{\mR^d}v(t,x)\mu(\dif x)=\int_{B_3}\log\big(\tfrac{w(t,x)}{\eta(x)}\big)\mu(\dif x)$. Hence,
$$
\tfrac{c_0}{2}\int^{t_1}_{t_0}\!\!\!\int_{B_4}\left(v(s,x)-V(s)\right)^2\mu(\dif x)\dif s\leq V(t_1)-V(t_0)+C_0,
$$
where by \eqref{KK1},
$$
C_0:=\frac{C_3}{\|\eta\|^2_{\mL^2_x}}\int^{2}_{-2}\Big(1+\frac{{\rm Tail}(w^-,4)}\kappa\Big)
+\frac{C_4}{\|\eta\|^2_{\mL^2_x}}\int^{2}_{-2}\!\int_{B_4}\Big(|\div b|^2+|b|^2+|b|+\frac{|f|} \kappa\Big)<\infty.
$$
By Lemma \ref{Le24} and \eqref{KA33}, we obtain \eqref{GF00} and \eqref{GF01} with $a=V(0)$.
\end{proof}

\bl\label{Le45}
There are constants $\gamma_0, C_0>0$ such that for all $1\leq\tau<\sigma\leq 2$ and $p\in(0,1]$,
$$
\big\|w^{-1}\b1_{Q^-_\tau(2)}\big\|_{\mL^\infty}\leq (C_0(\sigma-\tau)^{-\gamma_0})^{\frac1p}\big\|w^{-1}\b1_{Q^-_\sigma(2)}\big\|_{\mL^p}.
$$
\el
\begin{proof}
Let $\sI\subset\mR^{1+d}$ be as in \eqref{II}. By (i) of Theorem \ref{Th21}, it suffices to show that
$$
w^{-1}\in\cM^\infty_\sI(Q^-_\cdot(2)).
$$
More precisely, we want to show that there are $(\nu_i,\bbr_i)\in\sI$, $i=0,1,2,3$ such that 
for any $(\nu,\bbr)\in\sI$, there is a constant $C>0$ such that
 for  any $1\leq\tau<\sigma\leq 2$ and  $q>0$,
$$
(\sigma-\tau)^{(d+\alpha)/2}\|w^{-q}\b1_{Q^-_{\tau}(2)}\|_{\mL^\nu_t\mL^\bbr_x}
\lesssim_{C} (2q+1)\sum_{i=0}^3\|w^{-q}\b1_{Q^-_{\sigma}(2)}\|_{\mL^{\nu_i}_t\mL^{\bbr_i}_x},
$$
or equivalently,  for any $q>1$,
\begin{align}\label{DX-02}
(\sigma-\tau)^{(d+\alpha)/2}\|w^{(1-q)/2}\b1_{Q^-_{\tau}(2)}\|_{\mL^\nu_t\mL^\bbr_x}
\lesssim_{C} q\sum_{i=0}^3\|w^{(1-q)/2}\b1_{Q^-_{\sigma}(2)}\|_{\mL^{\nu_i}_t\mL^{\bbr_i}_x}.
\end{align}

Let $1\leq \tau<\sigma\leq 2$ and $\eta\in C^\infty_c(Q_{(\tau+\sigma)/2}(2);[0,1])$   with $\eta|_{Q_\tau(2)}\equiv 1$ and satisfy \eqref{DX877}.
Recalling \eqref{AA1},   we define for $q>1$, 
$$
v_\eps:=w_\eps^{(1-q)/2}.
$$
We estimate each term in \eqref{DZ1} with $q>1$ and $\bar q=q+1$ as  follows. 
\begin{enumerate}[$\bullet$]
\item By the chain rule and integration by parts, we have
\begin{align*}
\int_{\mR^d} \p_t w_\eps w_\eps^{-q}\eta^{q+1}
=\frac{1}{1-q}\int_{\mR^d}\p_t v_\eps^2\eta^{q+1}
=\frac{1}{1-q}\int_{\mR^d}\Big(\p_t( v_\eps^2\eta^{q+1})-v_\eps^2\p_t\eta^{q+1}\Big).
\end{align*}

\item For 
$\cJ^\eps_1 (t)$, 
 since $\eta$ has support in $Q_\sigma(2)$, by (ii) of Lemma \ref{Le26} we have
\begin{align*}
\cJ^\eps_1 (t) 
&=\frac12\int_{B_\sigma}\!\int_{B_\sigma}(w_\eps(y)-w_\eps(x))((w_\eps^{-q}\eta^{q+1})(x)-(w_\eps^{-q}\eta^{q+1})(y))K(t,x-y)\dif x\dif y\\
&\geq \frac1{2(q-1)}\int_{B_\sigma}\!\int_{B_\sigma}\left(\Big(\frac{\eta(x)}{w_\eps(x)}\Big)^{\frac{q-1}{2}}-\Big(\frac{\eta(y)}{w_\eps(y)}\Big)^{\frac{q-1}{2}}\right)^2K(t,x-y)\dif x\dif y\\
&-\frac{\ell(q)}2\int_{B_\sigma}\!\int_{B_\sigma}(\eta(x)-\eta(y))^2\left(\Big(\frac{\eta(x)}{w_\eps(x)}\Big)^{q-1}+\Big(\frac{\eta(y)}{w_\eps(y)}\Big)^{q-1}\right)K(t,x-y)\dif x\dif y\\
&=:\cJ^\eps_{11} (t) +\cJ^\eps_{12} (t).
\end{align*}
For $\cJ^\eps_{11} (t)$, for any $\bbp\in(2,\infty)^d$ with $|\frac1\bbp|=\frac{d-\alpha}2$, by Sobolev's embedding \eqref{SS1}, we have
\begin{align*}
\cJ^\eps_{11} (t) 
&\stackrel{\eqref{KK0}}{\geq}\frac{\kappa_0}{2(q-1)}\int_{B_\sigma}\!\int_{B_\sigma}\left[\big(v_\eps\eta^{\frac{q-1}2}\big)(x)-\big(v_\eps\eta^{\frac{q-1}2}\big)(y)\right]^2\frac{\dif x\dif y}{|x-y|^{d+\alpha}}\\
&\geq \tfrac {c_0}{q-1}\|v_\eps\eta^{\frac{q-1}2}\|^2_{\mL^{\bbp}_x}-\tfrac {C_0(\sigma-\tau)^{-d-\alpha}}{q-1}\|v_\eps\eta^{\frac{q-1}2}\|^2_{\mL^2_x}.
\end{align*}
For $\cJ^\eps_{12} (t) $, since $|\eta|\leq 1$ and $q>1$, 
by $\ell(q)=\max\{4,\frac{6q-5}{2}\}$ and \eqref{DX877},
we have
\begin{align*}
|\cJ^\eps_{12} (t) |&\stackrel{\eqref{KK0}}{\leq} \kappa_1\ell(q)\|\nabla\eta\|_\infty^2\int_{B_\sigma}\!\int_{B_\sigma}|x-y|^2\Big(\frac{\eta(x)}{w_\eps(x)}\Big)^{q-1}\frac{\dif x\dif y}{|x-y|^{d+\alpha}}\\
&\lesssim q (\sigma-\tau)^{-2}\|v_\eps\eta^{\frac{q-1}2}\|_{\mL^2_x}^2\int_{B_{4\sigma}}|y|^{2-d-\alpha}\dif y.
\end{align*}

\item For $\cJ^\eps_2 (t) $, noting that for $|x|\leq(\tau+\sigma)/2$ and $|y|\geq\sigma$,
$$
|x-y|^{-d-\alpha}\leq \tfrac{(2\sigma)^{d+\alpha}}{(\sigma-\tau)^{d+\alpha}}|y|^{-d-\alpha},
$$
due to $\kappa\leq w_\eps$ on $Q_4$ and  
supp($\eta$)$\subset Q_{(\tau+\sigma)/2}(2)\subset Q_\sigma(2)$, we have
\begin{align}
-\cJ^\eps_2 (t)
&=\int_{B_{(\tau+\sigma)/2}}\left(\int_{B^c_\sigma}(1-w_\eps(y)/w_\eps(x))K(t,x-y)\dif y\right)(v_\eps^2\eta^{q+1})(x)\dif x\no\\
&\stackrel{\eqref{KK0}}{\leq} \kappa_1\int_{B_{(\tau+\sigma)/2}}\left(\int_{B^c_\sigma}\frac{1+w^-_\eps(y)/\kappa}{|x-y|^{d+\alpha}}\dif y\right)(v_\eps^2\eta^{q+1})(x)\dif x\no\\
&\leq \tfrac{\kappa_1(2\sigma)^{d+\alpha}
\b1_{[0,4]} (t) }
{(\sigma-\tau)^{d+\alpha}}
\int_{B_{(\tau+\sigma)/2}}\left(\int_{B^c_\sigma}
\frac{1+w^-_\eps(y)/\kappa}{|y|^{d+\alpha}}\dif y\right)v^2_\eps(x)\dif x\no\\
&\lesssim
\tfrac{\b1_{[0,4]} (t) }
{(\sigma-\tau)^{d+\alpha}}\Big(1+\tfrac{{\rm Tail}(w^-_\eps,\sigma)}\kappa\Big)\|v_\eps\b1_{B_\sigma}\|^2_{\mL^2_x}.\label{KK3}
\end{align}
\item For $\cJ^\eps_3 (t) $, integrating it over $[0,2]$ with respect to the time variable 
and by the integration by parts and H\"older's inequality, we have
for $\frac2{\nu_i}+\frac{1}{q_i}=1$ and $\frac2{\bbr_i}+\frac{1}{\bbp_i}={\bf 1}$, $i=1,2$,
\begin{align}
\int^2_0|\cJ^\eps_3(s)|\dif s
&=\tfrac{1}{q-1}\int^2_0\left|\int_{\mR^d} w_\eps^{1-q} 
\left(b\cdot\nabla\eta^{q+1}+
\eta^{q+1} \, \div b 
\right)\right|\dif s\no\\
&\leq\tfrac{1}{q-1}\int^2_0\!\!\!\int_{\mR^d} v_\eps^2((q+1)|b\cdot\nabla\eta|\eta^q+|\div b|\eta^{q+1})\dif s\no\\
\begin{split}\label{ZD3}
&\leq\tfrac{q+1}{q-1}\|v_\eps\b1_{Q^-_\sigma(2)}\|_{\mL^{\nu_1}_t\mL^{\bbr_1}_x}^2\|b\b1_{Q^-_\sigma(2)}\|_{\mL^{q_1}_t\mL^{\bbp_1}_x}\|\nabla\eta\|_\infty\\
&\quad+\tfrac{1}{q-1}\|v_\eps\b1_{Q^-_\sigma(2)}\|_{\mL^{\nu_2}_t\mL^{\bbr_2}_x}^2\|\div b\b1_{Q^-_\sigma(2)}\|_{\mL^{q_2}_t\mL^{\bbp_2}_x}.
\end{split}
\end{align}

\item For $\cJ^\eps_4 (t) $, integrating it over $[0,2]$ with respect to the time variable 
and thanks to 
$w_\eps\geq\kappa$ on $Q_4$
and by H\"older's inequality, we have
for $\frac2{\nu_0}+\frac{1}{q_0}=1$ and $\frac2{\bbr_0}+\frac{1}{\bbp_0}={\bf 1}$,
\begin{align*}
\int^2_0\left|\cJ^\eps_4(s)\right|\dif s
&\leq\int_{Q^-_\sigma(2)} (v_\eps^2(|f|/\kappa)+\kappa^{-q}[\rho_\eps, b\cdot\nabla]u)\\
&\leq 
\|v_\eps\b1_{Q^-_\sigma(2)}\|_{\mL^{\nu_0}_t\mL^{\bbr_0}_x}^2\|f\b1_{Q^-_\sigma(2)}\|_{\mL^{q_0}_t\mL^{\bbp_0}_x}/\kappa
+\kappa^{-q}\|[\rho_\eps, b\cdot\nabla]u\b1_{Q^-_\sigma(2)}\|_{\mL^1}.
\end{align*}
\end{enumerate}
Substituting the above estimates into \eqref{DZ1} and integrating from $0$ to $t$, where $t\in[0,2]$, 
and multiplying both sides by $q-1$, by \eqref{KK1}, we obtain that for any $\bbp\in(2,\infty)^d$ with $|\frac1\bbp|=\frac{d-\alpha}2$, 
\begin{align*}
&\|(v_\eps\eta^{\frac{q+1}2})(t)\|^2_{\mL^2_x}+\|v_\eps\eta^{\frac{q-1}2}\|^2_{\mL^2_t\mL_x^{\bbp}}\\
&\quad\lesssim_C\Big(\|\p_t\eta^{q+1}\|_\infty+1+ \tfrac{(q-1)q}
{(\sigma-\tau)^2}\Big)\|v_\eps\b1_{Q^-_\sigma(2)}\|_{\mL^2}^2\\
&\qquad+\tfrac{q-1}{(\sigma-\tau)^{d+\alpha}}\Big(\tfrac{\|{\rm Tail}(w^-_\eps,\sigma)\b1_{[0,2]}\|_{\mL^1_t}}\kappa+1\Big)\|v_\eps\b1_{Q^-_\sigma(2)}\|^2_{\mL^2}\\
&\qquad+\tfrac{q+1}{\sigma-\tau}\|v_\eps\b1_{Q^-_\sigma(2)}\|_{\mL^{\nu_1}_t\mL^{\bbr_1}_x}^2+\|v_\eps\b1_{Q^-_\sigma(2)}\|_{\mL^{\nu_2}_t\mL^{\bbr_2}_x}^2\\
&\qquad+(q-1)\|v_\eps\b1_{Q^-_\sigma(2)}\|_{\mL^{\nu_0}_t\mL^{\bbr_0}_x}^2+\kappa^{-q}\|\b1_{Q^-_\sigma(2)}[\rho_\eps, b\cdot\nabla]u\|_{\mL^1},
\end{align*}
where the constant $C$ depends only on $d,\alpha$ and $\kappa_0$.
Let $v=w^{(1-q)/2}$.
Taking limit $\eps\to 0$ and by \eqref{Mo2} and then taking supremum in $t\in[0,2]$, we arrive at
\begin{align*}
&\|v\b1_{Q^-_\tau(2)}\|^2_{\mL^\infty_t\mL^2_x}+\|v\b1_{Q^-_\tau(2)}\|^2_{\mL^2_t\mL_x^\bbp}\\
&\quad\lesssim_C\Big(\tfrac{q+1}{\sigma-\tau}+1+
\tfrac{(q-1)q}{(\sigma-\tau)^2}\Big)\|v\b1_{Q^-_\sigma(2)}\|_{\mL^2}^2\\
&\qquad+\tfrac{q-1}{(\sigma-\tau)^{d+\alpha}}\Big(\tfrac{\|{\rm Tail}(w^-,\sigma)\b1_{[0,2]}\|_{\mL^1_t}}\kappa+1\Big)\|v\b1_{Q^-_\sigma(2)}\|^2_{\mL^2}\\
&\qquad+\tfrac{q+1}{\sigma-\tau}\|v\b1_{Q^-_\sigma(2)}\|_{\mL^{\nu_1}_t\mL^{\bbr_1}_x}^2+\|v\b1_{Q^-_\sigma(2)}\|_{\mL^{\nu_2}_t\mL^{\bbr_2}_x}^2\\
&\qquad+(q-1)\|v\b1_{Q^-_\sigma(2)}\|_{\mL^{\nu_0}_t\mL^{\bbr_0}_x}^2\\
&\quad\lesssim_C\tfrac{q^2}{(\sigma-\tau)^{d+\alpha}}\|v\b1_{Q^-_\sigma(2)}\|^2_{\mL^2}+\tfrac{q+1}{\sigma-\tau}\|v\b1_{Q^-_\sigma(2)}\|_{\mL^{\nu_1}_t\mL^{\bbr_1}_x}^2\\
&\qquad+\|v\b1_{Q^-_\sigma(2)}\|_{\mL^{\nu_2}_t\mL^{\bbr_2}_x}^2+(q-1)\|v\b1_{Q^-_\sigma(2)}\|_{\mL^{\nu_0}_t\mL^{\bbr_0}_x}^2,
\end{align*}
where in the second step we used \eqref{KK1}.
By Lemma \ref{Le22}, we obtain \eqref{DX-02}. The proof is complete.
\end{proof}

\bl\label{Le46}There are constants $p_0,\gamma_1, C_1>0$ such that for all $1\leq\tau<\sigma\leq 2$ and $p\in(0,1]$,
$$
\big\|w\b1_{Q^+_\tau(-2)}\big\|_{\mL^{p_0}}\leq (C_1(\sigma-\tau)^{-\gamma_1})^{\frac1p}\big\|w\b1_{Q^+_\tau(-2)}\big\|_{\mL^p}.
$$
\el
\begin{proof}
Let $\sI\subset\mR^{1+d}$ be as in \eqref{II}. By (ii) of Theorem \ref{Th21}, it suffices to show that
$$
w\in\cM^{1/4}_\sI(Q^+_\cdot(-2)).
$$
More precisely, we want to show that there are $(\nu_i,\bbr_i)\in\sI$, $i=0,1,2,3$ such that
for any $(\nu,\bbr)\in\sI$, there is a constant $C>0$ such that for any
 $1\leq\tau<\sigma\leq 2$ and $q\in(0,\frac14)$,
$$
(\sigma-\tau)^{(d+\alpha)/2}\|\b1_{Q^+_\tau(-2)}w^q\|_{\mL^\nu_t\mL^\bbr_x}
\lesssim_C \sum_{i=0}^3\|\b1_{Q^+_\tau(-2)}w^q\|_{\mL^{\nu_i}_t\mL^{\bbr_i}_x},
$$
 or equivalently, for any  $q\in(\frac12,1)$,
\begin{align}\label{DX-03}
(\sigma-\tau)^{(d+\alpha)/2}\|\b1_{Q^+_\tau(-2)}w^{(1-q)/2}\|_{\mL^\nu_t\mL^\bbr_x}
\lesssim_C\sum_{i=0}^3\|\b1_{Q^+_\tau(-2)}w^{(1-q)/2}\|_{\mL^{\nu_i}_t\mL^{\bbr_i}_x}.
\end{align}

Let $1\leq \tau<\sigma\leq 2$ and $\eta\in C^\infty_c(Q_{(\tau+\sigma)/2}(-2);[0,1])$   with $\eta|_{Q_\tau(-2)}\equiv 1$ and satisfy \eqref{DX877}.
Recalling \eqref{AA1}, for 
$q\in(\frac12,1)$,
we define
$$
v_\eps:=w_\eps^{(1-q)/2}.
$$
We estimate each term in \eqref{DZ1} with 
$q\in(\frac12,1)$
and $\bar q=2$ as follows. 
\begin{enumerate}[$\bullet$]
\item By the chain rule and integration by parts, we have
\begin{align*}
\int_{\mR^d}\p_t w_\eps w_\eps^{-q}\eta^2
=\tfrac{1}{1-q}\int_{\mR^d}\p_t v_\eps^2\eta^2
=\tfrac{1}{1-q}\int_{\mR^d}\p_t( v_\eps^2\eta^2)-\tfrac{1}{1-q}\int_{\mR^d} v_\eps^2\p_t\eta^2.
\end{align*}

\item For $\cJ^\eps_1 (t)$,  by (iii) of Lemma \ref{Le26} we have
\begin{align*}
\cJ^\eps_1 (t) &=\frac12\int_{B_\sigma}\!\int_{B_\sigma} 
(w_\eps(y)-w_\eps(x))((w_\eps^{-q}\eta^2)(x)-(w_\eps^{-q}\eta^2)(y))K(t,x-y)\dif x\dif y\\
&\geq \frac{\ell_1(q)}2\int_{B_\sigma}\!\int_{B_\sigma}
\big[\eta(x)w_\eps(x)^{\frac{1-q}{2}}-\eta(y)w_\eps(y)^{\frac{1-q}{2}}\big]^2K(t,x-y)\dif x\dif y\\
&-\frac{\ell_2(q)}2\int_{B_\sigma}\!\int_{B_\sigma}
(\eta(x)-\eta(y))^2\left[w_\eps(x)^{1-q}+w_\eps(y)^{1-q}\right]K(t,x-y)\dif x\dif y\\
&=: \cJ^\eps_{11} (t) +\cJ^\eps_{12} (t).
\end{align*}
For $\cJ^\eps_{11} (t)$,  for any $\bbp\in(2,\infty)^d$ with $|\frac1\bbp|=\frac{d-\alpha}2$
 and $q\in(\frac12,1)$,
by Sobolev's embedding \eqref{SS1} and $\ell_1(q)=\frac{2q}{3(1-q)}$, 
we have
\begin{align*}
\cJ^\eps_{11} (t) 
&\stackrel{\eqref{KK0}}{\geq}\frac{\kappa_0\ell_1(q)}2\int_{B_\sigma}\!\int_{B_\sigma}
\frac{\left[\eta(x)v_\eps(x)-\eta(y)v_\eps(y)\right]^2}{|x-y|^{d+\alpha}}\dif x\dif y\\
&\geq \tfrac{1}{1-q}\left(c_0\|\eta v_\eps\|^2_{\mL^{\bbp}_x}-C_0(\sigma-\tau)^{-d-\alpha}\|\eta v_\eps\|^2_{\mL^2_x}\right).
\end{align*}
For $\cJ^\eps_{12} (t) $, since supp($\eta$)$\subset Q_{(\tau+\sigma)/2}\subset Q_\sigma$, 
by $\ell_2(q)=\frac{4q}{1-q}+\frac{9}{q}$  and \eqref{DX877},
we have
\begin{align*}
|\cJ^\eps_{12} (t) |
&\stackrel{\eqref{KK0}}{\leq} 2\kappa_1\ell_2(q)\b1_{[-4,0]}\|\nabla\eta\|_\infty^2\int_{B_\sigma}\!\int_{B_\sigma}
\frac{|x-y|^2v^2_\eps(x)}{|x-y|^{d+\alpha}}\dif x\dif y\\
&\lesssim \tfrac{1}{(1-q)(\sigma-\tau)^2}\b1_{[-4,0]}\|v_\eps\b1_{Q_\sigma}\|_{\mL^2_x}^2\int_{B_{4\sigma}}|y|^{2-d-\alpha}\dif y.
\end{align*}

\item For $\cJ^\eps_2 (t)$, as in \eqref{KK3} we have
\begin{align*}
-\cJ^\eps_2 (t)
\lesssim\tfrac{\b1_{[-4,0]}}{(\sigma-\tau)^{d+\alpha}}\Big(1+\tfrac{{\rm Tail}(w^-_\eps,\sigma)}\kappa\Big)\|v_\eps\b1_{B_\sigma}\|^2_{\mL^2_x}.
\end{align*}

\item For $\cJ^\eps_2 (t) $, as in \eqref{ZD3}, integrating it over $[-2,0]$ with respect to the time variable 
and by the integration by parts and H\"older's inequality, we have
for $\frac2{\nu_i}+\frac{1}{q_i}=1$ and $\frac2{\bbr_i}+\frac{1}{\bbp_i}=1$ with $i=1,2$,
\begin{align*}
\int^0_{-2}|\cJ^\eps_3(s)|\dif s
&=\tfrac{1}{1-q}\int^0_{-2}\left|\int_{\mR^d} w_\eps^{1-q}(b\cdot\nabla\eta^2+
\eta^2 \, \div b )\right|\dif s\no\\
&\leq\tfrac{1}{1-q}\int^0_{-2}\!\int_{\mR^d} v_\eps^2(2|b\cdot\nabla\eta|\eta+|\div b| \, \eta^2)\dif s\no\\
&\leq\tfrac{1}{1-q}\|v_\eps\b1_{Q^+_\sigma(-2)}\|_{\mL^{\nu_1}_t\mL^{\bbr_1}_x}^2\|b\b1_{Q^+_\sigma(-2)}\|_{\mL^{q_1}_t\mL^{\bbp_1}_x}\|\nabla\eta\|_\infty\\
&\quad+\tfrac{1}{1-q}\|v_\eps\b1_{Q^+_\sigma(-2)}\|_{\mL^{\nu_2}_t\mL^{\bbr_2}_x}^2\|\div b\b1_{Q^+_\sigma(-2)}\|_{\mL^{q_2}_t\mL^{\bbp_2}_x}.
\end{align*}

\item For $\cJ^\eps_4 (t)$, integrating it over $[-2,0]$ with respect to the time variable 
and thanks to  
$w_\eps\geq\kappa$ on $Q_4$
and by H\"older's inequality, we have
for $\frac2{\nu_0}+\frac{1}{q_0}=1$ and $\frac2{\bbr_0}+\frac{1}{\bbp_0}=1$,
\begin{align*}
\int^0_{-2}\left|\cJ^\eps_4(s)\right|\dif s\leq 
\|v_\eps\b1_{Q^+_\sigma(-2)}\|_{\mL^{\nu_0}_t\mL^{\bbr_0}_x}^2\|f\b1_{Q^+_\sigma(-2)}\|_{\mL^{q_0}_t\mL^{\bbp_0}_x}/\kappa
+\kappa^{-q}\|[\rho_\eps, b\cdot\nabla]u\b1_{Q^+_\sigma(-2)}\|_{\mL^1}.
\end{align*}
\end{enumerate}
Substituting the above estimates into \eqref{DZ1} and integrating from $t$ to $0$, where $t\in[-2,0]$, 
and multiplying both sides by $1-q$, by \eqref{KK1}, we obtain that  for any $\bbp\in(2,\infty)^d$ with $|\frac1\bbp|=\frac{d-\alpha}2$, 
\begin{align*}
& \|(v_\eps\eta)(t)\|^2_{\mL^2_x}+\int^0_{t}\|(v_\eps\eta)(s)\|^2_{\mL_x^\bbp}\dif s  \\
&\quad \lesssim 
 \Big(|\p_t\eta^2\|_\infty+1+\tfrac{\ell_2(q)}{(\sigma-\tau)^2}\Big)\|v_\eps\b1_{Q^+_\sigma(-2)}\|_{\mL^2}^2 \\
&\qquad +\tfrac{1}{(\sigma-\tau)^{d+\alpha}}\Big(\tfrac{\|{\rm Tail}(w^-_\eps,\sigma)\b1_{[-2,0]}\|_{\mL^1_t}}\kappa+1\Big)\|v_\eps\b1_{Q^+_\sigma(-2)}\|^2_{\mL^2}\\
&\qquad +\tfrac{1}{\sigma-\tau}\|v_\eps\b1_{Q^+_\sigma(-2)}\|_{\mL^{\nu_1}_t\mL^{\bbr_1}_x}^2+\|v_\eps\b1_{Q^+_\sigma(-2)}\|_{\mL^{\nu_2}_t\mL^{\bbr_2}_x}^2  \\
&\qquad  +\|v_\eps\b1_{Q^+_\sigma(-2)}\|_{\mL^{\nu_0}_t\mL^{\bbr_0}_x}^2+\kappa^{-q}\|\b1_{Q^+_\sigma(-2)}[\rho_\eps, b\cdot\nabla]u\|_{\mL^1},
\end{align*}
where the constant $C$ depends only on $d,\alpha$ and $\kappa_0$.
Taking limit $\eps\to 0$ and by \eqref{Mo2} and then taking supremum in $t\in[0,2]$, we arrive at
\begin{align*}
&\|v\b1_{Q^+_\tau(-2)}\|^2_{\mL^\infty_t\mL^2_x}+\|v\b1_{Q^+_\tau(-2)}\|^2_{\mL^2_t\mL_x^\bbp} \\
&\quad\lesssim \Big(\tfrac{1}{\sigma-\tau}+1+\tfrac{\ell_2(q)}{(\sigma-\tau)^2}\Big)\|v\b1_{Q^+_\sigma(-2)}\|_{\mL^2}^2\\
&\qquad  +\tfrac{2}{(\sigma-\tau)^{d+\alpha}}\Big(\tfrac{\|{\rm Tail}(w^-,\sigma)\b1_{[-2,0]}\|_{\mL^1_t}}\kappa+1\Big)\|v\b1_{Q^+_\sigma(-2)}\|^2_{\mL^2}\\
&\qquad+\tfrac{2}{\sigma-\tau}\|v\b1_{Q^+_\sigma(-2)}\|_{\mL^{\nu_1}_t\mL^{\bbr_1}_x}^2+\|v\b1_{Q^+_\sigma(-2)}\|_{\mL^{\nu_2}_t\mL^{\bbr_2}_x}^2
+\|v\b1_{Q^+_\sigma(-2)}\|_{\mL^{\nu_0}_t\mL^{\bbr_0}_x}^2\\
&\quad\lesssim  \tfrac{1}{(\sigma-\tau)^{d+\alpha}}\|v\b1_{Q^+_\sigma(-2)}\|^2_{\mL^2}+\tfrac{1}{\sigma-\tau}\|v\b1_{Q^+_\sigma(-2)}\|_{\mL^{\nu_1}_t\mL^{\bbr_1}_x}^2 \\
&\qquad +\|v\b1_{Q^+_\sigma(-2)}\|_{\mL^{\nu_2}_t\mL^{\bbr_2}_x}^2+\|v\b1_{Q^+_\sigma(-2)}\|_{\mL^{\nu_0}_t\mL^{\bbr_0}_x}^2,
\end{align*}
where in the second step we have used $\ell_2(q)=\frac{4q}{1-q}+\frac{9}{q}$ and \eqref{KK1}.
Now by Lemma \ref{Le22}, we obtain \eqref{DX-03}. This completes the proof.  
\end{proof}

 We are now in a position to give the 

\begin{proof}[Proof of Theorem \ref{Th41}]
Let $C_0$ and $a$ be as in Lemma \ref{Le42} and define $h(t,x):=\e^{a-C_0}/w(t,x)$. Then by \eqref{GF01}, we have for some $C_1=C_1(d,\alpha)>0$,
$$
(\dif t\otimes\dif x)\big\{(t,x)\in Q^+_2(0): \log h(t,x)>\lambda\big\}\leq C_1\lambda^{-1},\ \ \lambda>0.
$$
Moreover, by Lemma \ref{Le45}, we have for any 
$p\in(0,1)$ and $1\leq\tau<\sigma\leq 2$,
$$
\big\|h\b1_{Q^-_\tau(2)}\big\|_{\mL^\infty}\leq (C(\sigma-\tau)^{-\gamma_0})^{1/p}\big\|h\b1_{Q^-_\sigma(2)}\big\|_{\mL^p}.
$$
Since $Q^-_2(2)=Q_2^+(0)$, we can use Lemma \ref{Le28} to deduce that for some $c_0\in(0,1)$,
\begin{align}\label{UU1}
\|h\b1_{Q^-_{3/2}(2)}\|_{\mL^\infty}\leq c_0^{-1}\Rightarrow c_0\e^{a-C_0}\leq {\rm inf}_{Q^-_{3/2}(2)}w.
\end{align}
On the other hand, if we let $\bar h(t,x):=w(t,x)\e^{-a-C_0}$.
Then by \eqref{GF00}, we have
$$
(\dif t\otimes\dif x)\big\{(t,x)\in Q^-_2(0): \log \bar h(t,x)>\lambda\big\}\leq C_1\lambda^{-1}.
$$
Moreover, by Lemma \ref{Le46}, we have for some $p_0>0$ and 
all $p\in(0,1)$ and $1\leq\tau<\sigma\leq 2$,
$$
\big\|\bar h\b1_{Q^+_\tau(-2)}\big\|_{\mL^{p_0}}\leq(C(\sigma-\tau)^{-\gamma_1})^{\frac1p}\big\|\bar h\b1_{Q^+_\sigma(-2)}\big\|_{\mL^p}.
$$
Since $Q^+_2(-2)=Q_2^-(0)$,  we can use Lemma \ref{Le28} to deduce that
$$
\|\bar h\b1_{Q_{3/2}^+(-2)}\|_{\mL^{p_0}}\leq C\Rightarrow \|w\b1_{Q_{3/2}^+(-2)}\|_{\mL^{p_0}}\leq C\e^{a+C_0},
$$
which together with \eqref{UU1} yields \eqref{WH1}.
Inequality  \eqref{WH2}  follows
 from  \eqref{WH1}, Theorem \ref{Le32} 
and a shifting the time variable and suitable scaling 
argument.
\end{proof}

\bc\label{Co46}
Let $(q_0,\bbp_0)\in\mI_\alpha^0$.
Under {\bf (H$'_b$)},
there is some $\theta\in(0,1)$ depending only on $\kappa_0,q_i,\bbp_i, d,\alpha$ such that
for any $f$ with $\|f\b1_{Q_4}\|_{\mL^{q_0}_t\mL^{\bbp_0}_x}<\infty$ 
and any solution $u$ of \eqref{PDE0} in $Q_6$ with $u\b1_{Q_4}\geq 0$,
\begin{align}\label{Osc1}
\underset{Q_{1/2}}{\osc} u\leq \theta\,\underset{Q_6}{\osc}\, u+\|f\b1_{Q_4}\|_{\mL^{q_0}_t\mL^{\bbp_0}_x}
+\|u\b1_{[-4,4]}\|_{\mL^\infty}.
\end{align}
\ec

\begin{proof}
Following the proof of \cite[Theorem 5.4.7]{SC1},
by \eqref{WH2}, if $u\b1_{Q_4}\geq 0$, then
\begin{align}\label{A19}
\fint_{Q^+_1(-2)}u\leq\sup_{Q^+_1(-2)} u\lesssim\inf_{Q^-_1(2)} u+\sA,
\end{align}
where
$$
\sA:=\|f\b1_{Q_2}\|_{\mL^{q_0}_t\mL^{\bbp_0}_x}+\|u\b1_{[-4,4]}\|_{\mL^\infty_t}.
$$
Let $M_u, m_u$ be the supremum and infimum of $u$ in $Q_4$. Similarly, let $M^+_u$ and $m^+_u$ be the  supremum and infimum of $u$ in $Q_1^-(2)$.
By applying \eqref{A19} to the non-negative solutions $M_u-u$ and $u-m_u$, we obtain that for some $C_0>1$,
$$
M_u-\fint_{Q^+_1(-2)}u\lesssim_{C_0}\inf_{Q^-_1(2)}(M_u-u)+\sA/2
$$
and
$$
\fint_{Q^+_1(-2)}u-m_u\lesssim_{C_0}\inf_{Q^-_1(2)}(u-m_u)+\sA/2.
$$
Hence,
$$
M_u-m_u\lesssim_{C_0} M_u-M^+_u+m^+_u-m_u+\sA,
$$
which implies that
$$
C_0(M^+_u-m^+_u)\leq (C_0-1)(M_u-m_u)+C_0\sA.
$$
Dividing both sides by $C_0$, we obtain 
$$
\underset{Q_1^-(2)}{\osc} u\leq \theta\,\underset{Q_4}{\osc}\, u+\|f\b1_{Q_2}\|_{\mL^{q_0}_t\mL^{\bbp_0}_x}
+\|u\b1_{[-4,4]}\|_{\mL^1_t}.
$$
Using the above estimate for $u(t+3/2,x)$ and by the monotonicity of $\osc_Q u$ in $Q$, we get the desired estimate.
\end{proof}

\br\rm
If $\alpha\in[1,2)$ and $(q_1,\bbp_1)$ in  {\bf (H$'_b$)} satisfies $\frac{\alpha}{q_1}+|\frac1{\bbp_1}|=\alpha-1$
(the critical case), 
then 
by \eqref{Osc1} and a standard scaling 
argument, 
one can show that
for some $\gamma\in(0,1)$ and $C>0$,
$$
\sup_{s\not=t\in[0,T]}\sup_{x\not=y\in B_4}\frac{|u(t,x)-u(s,y)|}{(|t-s|+|x-y|)^\gamma}\lesssim_C 
\|f\b1_{Q_4}\|_{\mL^{q_0}_t\mL^{\bbp_0}_x}
+\|u\b1_{[-4,4]}\|_{\mL^\infty}.
$$
\er

\section{Global $L^\infty$-estimate and H\"older regularity}

To derive  the global $L^\infty$ and H\"older estimates, we introduce some localized spaces as in \cite{ZZ18}. 
Let $\chi\in C^\infty_c(\mR^d)$ be a smooth function with $\chi(x)=1$ for $|x|\leq 1$ and $\chi(x)=0$ for $|x|>2$.
For $r>0$ and $z\in\mR^d$, define
\begin{align}\label{CHI}
\chi^z_r(x):=\chi((x-z)/r).
\end{align}
Fix $r>0$. For $\beta\in\mR$ and $\bbp\in(1,\infty)^d$, we introduce the following localized $\bH^\beta_\bbp$-space:
\begin{align}\label{Ck}
\widetilde \bH^\beta_\bbp:=\Big\{f: \nor f\nor_{\bH^\beta_\bbp}:=\sup_z\|f\chi^z_r\|_{\bH^\beta_\bbp}
=\sup_z\|f(\cdot-z)\chi^0_r\|_{\bH^\beta_\bbp}<\infty\Big\},
\end{align}
and  the localized space-time function space 
\begin{align}\label{GG1}
\wt\mL^q_T(\widetilde \bH^\beta_\bbp):=\Big\{f: \nor f\nor_{\wt\mL^q_T(\widetilde \bH^\beta_\bbp)}:=
\sup_{z\in\mR^d}\| f\chi^z_r\|_{\mL^q_T(\bH^\beta_\bbp)}<\infty\Big\},
\end{align}
where for $T>0$ and a Banach space $\mB$,
$$
\mL^q_T(\mB):=L^q([0,T];\mB).
$$
The following properties of $\wt\mL^q_T(\widetilde \bH^\beta_\bbp)$  
follows easily from 
its
definition (see \cite[Proposition 4.1]{ZZ18}).

\begin{enumerate}[{\rm(i)}]
\item For $r\not=r'>0$, there is a constant $C=C(d,\alpha,r,r',p,q)\geq 1$ such that
\begin{align}\label{GT1}
C^{-1}\sup_{z}\|f\chi^{z}_{r'}\|_{\mL^q_T(\bH^\beta_\bbp)}\leq \sup_{z}\|f\chi^{z}_r\|_{\mL^q_T(\bH^\beta_\bbp)}
\leq C \sup_{z}\|f\chi^{z}_{r'}\|_{\mL^q_T(\bH^\beta_\bbp)}.
\end{align}
In other words, the definition of 
the space $\wt\mL^q_T(\widetilde \bH^\beta_\bbp)$ does not depend on the choice of $r$.
\item For any $f\in\wt\mL^q_T(\widetilde \bH^\beta_\bbp)$, it holds that 
$$
f_\eps(t,x):=(f*\rho_\eps)(t,x)\in C^\infty_b([0,T]\times\mR^d),
$$
 and for some $C=C(T,d,\beta,\bbp,q)>0$,
\begin{align}\label{GT2}
\nor f_\eps\nor_{\wt\mL^q_T(\widetilde \bH^\beta_\bbp)}\leq C\nor f\nor_{\wt\mL^q_T(\widetilde \bH^\beta_\bbp)},\ \forall \eps\in(0,1),
\end{align}
and for any $\varphi\in C^\infty_c(\mR^{d})$,
\begin{align}\label{GT3}
\lim_{\eps\to0}\|(f_\eps-f)\varphi\|_{\mL^q_T(\bH^\beta_\bbp)}=0.
\end{align}
Here  $\rho_\eps$ is the mollifier defined in \eqref{e:3.3}.
\end{enumerate}

\subsection{Global well-posedness}

 First of all,  we present  
the following existence and uniqueness result of global weak solutions for PDE \eqref{PDE0} with initial value $u(t)|_{t\leq 0}=0$. 
In the following we set $\cI_t:=\b1_{[0,t]}$ and assume
$$
b(t,\cdot)|_{t\leq 0}=0 \quad \hbox{and} \quad f(t,\cdot)|_{t\leq 0}=0. 
$$

\bt\label{Th51}
Let $\beta\in[0,\alpha/2)$ and $T>0$. Suppose that for some $(q_1,\bbp_1)$ and $ (q_2,\bbp_2)\in\mI^0_\alpha$,
\begin{align}\label{CON11} 
\kappa_0:=\nor b\nor_{\wt\mL^{q_1}_T\wt\mL^{\bbp_1}_x}+\nor\div b\nor_{\wt\mL^{q_2}_T\wt\mL^{\bbp_2}_x}<\infty.
\end{align}
Then for any $f\in \wt\mL^{q_0}_T(\wt\bH^{-\beta}_{\bbp_0})$, where $(q_0,\bbp_0)\in\mI_\alpha^\beta$,  
there is a unique global weak solution $u$ to PDE \eqref{PDE0}  with initial value $u(t)|_{t\leq 0}=0$ such that for some $C=C(T,\kappa_0,q_i,\bbp_i)>0$
and all $t\in[0,T]$,
\begin{align}\label{EA80}
\|u\cI_t\|_{\mL^\infty_T\mL^\infty_x}+\nor \Delta^{\alpha/4}u\cI_t\nor_{\widetilde \mL^2_T\wt\mL^2_x}\lesssim_C
\nor f\cI_t\nor_{\wt\mL^{q_0}_T(\wt\bH^{-\beta}_{\bbp_0})}.
\end{align}
\et

\begin{proof}
We divide the proof into two steps.

{\sc Step 1}. We first show \eqref{EA80}.  Fix $x_0\in\mR^{d}$. Define
$$
u_{x_0}(t,x):=u(t,x-x_0),\quad  b_{x_0}(t,x):=b(t,x-x_0),\quad  f_{x_0}(t,x):=f(t,x-x_0).
$$
It is easy to check  that $u_{x_0}$ is a weak solution of PDE \eqref{PDE0} but with $(b_{x_0}, f_{x_0})$
in place of $(b, f)$ there.  
Applying Theorem \ref{Le32} to $u^+_{x_0}$ and  $u^-_{x_0}$ with $p=1$, respectively,   
there is a constant $C=C(\kappa_0,q_i,\bbp_i)>0$ such that for all $t\in[0,T]$ and $x_0\in\mR^{d}$,
\begin{align}\label{DD50}
&\|u_{x_0}\b1_{Q_1}\cI_t\|_{\mL^\infty}+\|\Delta^{\alpha/4}u_{x_0}\chi^0_1\cI_t\|_{\mL^2}\no\\
&\quad\lesssim_C \|u_{x_0}\b1_{Q_2}\cI_t\|_{\mL^1}+\|{\rm Tail}(u_{x_0}; 1)\cI_t\|_{\mL^2_t}
+\|\chi^0_2f_{x_0}\cI_t\|_{\mL^{ q_0}_t(\bH^{-\beta}_{\bbp_0})}\no\\
&\quad\lesssim_C\int^t_0\|u(s)\|_{\mL^\infty}\dif s+\left(\int^t_0\|u(s)\|^2_{\mL^\infty}\dif s\right)^{1/2}
+\nor f\cI_t\nor_{\wt\mL^{ q_0}_T(\wt\bH^{-\beta}_{\bbp_0})}.
\end{align}
Taking supremum in $x_0\in\mR^{d}$ and square for both sides, we obtain that for all $t\in[0,T]$,
\begin{align*}
\|u(t)\|^2_{\mL^\infty}\lesssim_C\int^t_0\|u(s)\|^2_\infty\dif s
+\nor f\cI_t\nor^2_{\wt\mL^{ q_0}_T(\wt\bH^{-\beta}_{\bbp_0})},
\end{align*}
 which  by Gronwall's inequality  yields that
$$
\sup_{s\in[0,t]}\|u(s)\|_{\mL^\infty}\lesssim\nor f\cI_t\nor_{\wt\mL^{ q_0}_T(\wt\bH^{-\beta}_{\bbp_0})}.
$$
 Thus, we obtain \eqref{EA80} by \eqref{DD50}, which in particular gives  the uniqueness of weak solutions.
 
 \medskip
 
{\sc Step 2}. In this step we show the existence of weak solutions. Define
$$
 b_\eps:=b*\rho_\eps \quad \hbox{and} \quad   f_\eps:=f*\rho_\eps.
$$
By \eqref{CON11}, we have for each $\eps\in(0,1)$,
$b_\eps, f_\eps\in C^\infty_b([0,T]\times\mR^{d})
$
and 
\begin{align}\label{DS3}
\sup_{\eps \in (0, 1]} 
\Big(\nor b_\eps\nor_{\wt\mL^{q_1}_T\wt\mL^{\bbp_1}_x}+\nor \div b_\eps\nor_{\wt\mL^{q_2}_T\wt\mL^{\bbp_2}_x}\Big)<\infty . 
\end{align}
It is well known that there is a unique smooth solution $u_\eps$ 
to  the following PDE:
\begin{align}\label{PDE10}
\p_t u_\eps=\Delta^{\alpha/2} u_\eps+b_\eps\cdot \nabla u_\eps+f_\eps
\quad \hbox{with } u_\eps(0)=0.
\end{align}
In particular, for any $\varphi\in C^\infty_c([0,T]\times\mR^d)$,
\begin{align}\label{PDE11}
\begin{split}
\int_{\mR^{1+d}}u_\eps\p_t\varphi&=-\int_{\mR^{1+d}}\<u_\eps,\Delta^{\alpha/2}\varphi\>
-\int_{\mR^{1+d}}(b_\eps\cdot\nabla \varphi+
\varphi \,  \div b_\eps 
)u_\eps+\int_{\mR^{1+d}}f_\eps\varphi.
\end{split}
\end{align}
Moreover, by \eqref{DS3} and \eqref{EA80},  
\begin{align}\label{SD1}
\sup_{\eps\in(0,1]}(\|u_\eps\|_{\mL^\infty_T\mL^\infty_x}+\nor \Delta^{\alpha/4}u_\eps\nor_{\widetilde \mL^2_T\wt\mL^2_x})<\infty.
\end{align}
Hence, there is a sequence $\eps_k\to 0$ and $u\in \mL^\infty$ such that
for each $\varphi\in C^\infty_c([0,T]\times\mR^d)$,
$$
\lim_{k\to \infty}  \int_{\mR^{1+d}}u_{\eps_k}\varphi  =  \int_{\mR^{1+d}}u\varphi.
$$
Taking limits on both sides of \eqref{PDE11} along $\eps_k$,
we conclude that  $u$ is a weak solution of PDE \eqref{PDE0}.  This completes the proof.
\end{proof}

\subsection{H\"older regularity for $\alpha\in(0,1]$}
In this subsection we establish a priori
 H\"older regularity in the following critical case:
\begin{align}\label{CC9}
\|b\|_{\rm crit}:=\b1_{\alpha\in(0,1)}\|b\|_{\mL^\infty_t(\dot\cC^{1-\alpha})}+\b1_{\alpha=1}\|b\|_{\mL^\infty_t({\rm BMO})}<\infty,
\end{align}
where for $\beta\in(0,1)$,
$$
\|f\|_{\dot\cC^{\beta}}:=\sup_{ x\not=y\in\mR^d} \frac{|f(x)-f(y)|}{|x-y|^\beta},
$$
and BMO  is the usual BMO space, that is., for $f_r (x):=\fint_{B_r} f(x+y)\dif y$,
$$
\|f\|_{\rm BMO}:=\sup_{r>0,x\in\mR^d}\left(\fint_{B_r}|f(x+y)-f_r (x)|\dif y\right).
$$
Note that for any $p\geq 1$ (cf. \cite[Corollary 7.1.9]{G19}),
\begin{align}\label{59}
\|f\|_{\rm BMO}&\asymp\sup_{r>0,x\in\mR^d}\left(\fint_{B_r}|f(x+y)-  f_r (x)|^p\dif y\right)^{1/p}.
\end{align}

The following is the main result of this section.

\bt\label{Th52}
Let $\alpha\in(0,1]$ and $T>0$. Suppose that 
for some $\bbp\in(1,\infty)^d$ with $|\frac1{\bbp}|<\alpha$,
\begin{align}\label{CN1}
\kappa_0:=\|\div b\|_{\mL^\infty_T\mL^\infty_x}+\nor b\nor_{\wt\mL^\infty_T\wt\mL^\bbp_x}+\|b\|_{\rm crit}<\infty.
\end{align}
For any $f\in \wt\mL^{q_0}_T(\wt\mL^{\bbp_0}_x)$, where $(q_0,\bbp_0)\in\mI^0_\alpha$,
 there are constants $\gamma\in(0,1)$ and $C>0$ depending only on $T,\kappa_0$, $q_0, \bbp_0,d,\alpha$ such that
$$
\sup_{s\not=t\in[0,T]}\sup_{x\not=y}\frac{|u(t,x)-u(s,y)|}{(|t-s|+|x-y|)^\gamma}\leq C\nor f\nor_{\wt\mL^{q_0}_T\wt\mL^{\bbp_0}_x}.
$$
\et

\br\rm
In the case of $\alpha=1$, for $p\geq 1$, by \eqref{59}, there is a constant $C=C(d,p)>0$ such that for each $z\in\mR^d$,
\begin{align*}
\|b(\cdot,\cdot-z)\b1_{B_1}\|_{\mL^\infty_T\mL^p_x}
\leq  C\|b\|_{\mL^\infty_T({\rm BMO})}
+\|b(\cdot,\cdot-z)\b1_{B_1}\|_{\mL^\infty_T\mL^1_x}.
\end{align*}
Taking supremum in $z\in\mR^d$, we get
$$
\nor b\nor_{\wt\mL^\infty_T\wt\mL^p_x}\leq C\|b\|_{\mL^\infty_T({\rm BMO})}
+\nor b\nor_{\wt\mL^\infty_T\wt\mL^1_x}.
$$
In this case, we can replace \eqref{CN1} by 
$$
\kappa_0:=\|\div b\|_{\mL^\infty_T\mL^\infty_x}+\nor b\nor_{\wt\mL^\infty_T\wt\mL^1_x}+\|b\|_{\rm crit}<\infty.
$$
\er

For the proof of Theorem \ref{Th52}, we  will
 use the Cafferalli-Vaseur method \cite{CV10}. We need the following technical lemma (cf. \cite[Proposition 7.1.5]{G19}).
\bl
For $\beta\in(0,1]$, 
there is a constant $C=C(d)>0$ such that for any $r>0$ and $x_1,x_2\in\mR^d$ with $|x_1-x_2|\leq 1$,
\begin{align}\label{DF1}
\left|\fint_{B_{r}(x_1)} f-\fint_{B_{r}(x_2)} f\right|
\lesssim_C
\left\{
\begin{aligned}
&\|f\|_{\rm   BMO }
\log\(\e+\tfrac{|x_1-x_2|}{r}\),\\
&\|f\|_{\dot\cC^\beta}|x_1-x_2|^\beta.
\end{aligned}
\right.
\end{align}
\el

Fix $\lambda>0$. Let $\xi_\lambda(t)$ be a solution of the following ODE in $\mR^d$:
$$
\dot \xi_\lambda(t)=\lambda^{\alpha-1}\fint_{B_2}b(\lambda^\alpha t,\lambda(y-\xi_\lambda(t)))\dif y
=\lambda^{\alpha-1}\fint_{B_{2\lambda}(\lambda\xi_\lambda(t))}b(\lambda^\alpha t,\cdot),
$$
with the initial value $\xi_\lambda(0)=0$.
Since $x\mapsto \lambda^{\alpha-1}\fint_{B_{2\lambda}(\lambda x)}b(\lambda^\alpha t,\cdot)$ is continuous and 
has linear growth by \eqref{DF1}, 
by Peano's theorem, there is at least one solution to the above ODE.
Let $u$ be a weak solution of PDE \eqref{PDE0}. Define
$$
u_\lambda(t,x):=u(\lambda^\alpha t,\lambda(x-\xi_\lambda(t))),\quad 
f_\lambda(t,x):=\lambda^\alpha f(\lambda^\alpha t, \lambda(x-\xi_\lambda(t)))
$$
and
$$
b_\lambda(t,x):=\lambda^{\alpha-1}b(\lambda^\alpha t, \lambda(x-\xi_\lambda(t)))-\dot\xi_\lambda(t).
$$
Then by the chain rule, it is easy to see that
$$
\p_t  u_\lambda=\Delta^{\alpha / 2}  u_\lambda+ b_\lambda\cdot\nabla  u_\lambda+ f_\lambda.
$$
Clearly we have
\begin{align}\label{CCA1}
\int_{B_2} b_\lambda(t,x)\dif x\equiv0,\quad   \|b_\lambda\|_{\rm crit}=\|b\|_{\rm crit},
\end{align}
and there is a constant $C=C(d,q_0,\bbp_0)>0$ such that for all $\lambda\in(0,1)$,
\begin{align}
\label{CCA1b} 
\|f_\lambda\b1_{Q_2}\|_{\mL^{q_0}_t\mL^{\bbp_0}_x} \leq C\lambda^{\alpha-\frac{\alpha}{q_0}-|\frac{1}{\bbp_0}|}\| f\|_{\mL^{q_0}_t\mL^{\bbp_0}_x}.
\end{align}

\bl
For any $p\geq 1$, there is a constant $C=C(p,d,\alpha)>0$ such that for any $\lambda>0$,
\begin{align}\label{CCA0}
\|b_\lambda\b1_{Q_2}\|_{\mL^\infty_t\mL^{p}_x}\leq C\, \|b\|_{\rm crit}.
\end{align}
\el
\begin{proof}
By definition we have
\begin{align*}
\|b_\lambda(t)\b1_{B_2}\|^p_{\mL^{p}_x}&=\int_{B_2}\left|\lambda^{\alpha-1}b(\lambda^\alpha t, \lambda(x-\xi_\lambda(t)))
-\dot\xi_\lambda(t)\right|^p\dif x\\
&\leq \lambda^{(\alpha-1)p}\sup_z\int_{B_2}\left|b(\lambda^\alpha t, \lambda x-z)-\fint_{B_2} b(\lambda^\alpha t,\lambda y-z)\dif y\right|^p\dif x\\
&= \lambda^{(\alpha-1)p-d}\sup_z\int_{B_{2\lambda}}\left|b(\lambda^\alpha t, x-z)-\fint_{B_{2\lambda}} b(\lambda^\alpha t,y-z)\dif y\right|^p\dif x.
\end{align*}
If $\alpha=1$, by \eqref{59}, we clearly have for some $C=C(d)>0$,
$$
\|b_\lambda(t)\b1_{B_2}\|^p_{\mL^{p}_x}\leq C\|b(\lambda^\alpha t)\|^p_{{\rm BMO}}.
$$
If $\alpha\in(0,1)$, noting that
$$
\left|b(\lambda^\alpha t, \lambda x-z)-\fint_{B_2} b(\lambda^\alpha t, \lambda y-z)\dif y\right|
\leq\|b(\lambda^\alpha t)\|_{\dot \cC^{1-\alpha}}\fint_{B_2}(\lambda|x-y|)^{1-\alpha}\dif y,
$$
as above, by definition we have
$$
\|b_\lambda(t)\b1_{B_2}\|^p_{\mL^{p}_x}\leq \|b(\lambda^\alpha t)\|^p_{\dot \cC^{1-\alpha}}\int_{B_2}\left(\fint_{B_2}|x-y|^{1-\alpha}\dif y\right)^p\dif x
\leq C\|b(\lambda^\alpha t)\|^p_{\dot \cC^{1-\alpha}}.
$$
The proof is complete.
\end{proof}

\bl\label{Le54}
Under \eqref{CC9}, there is a constant $C=C(d,\alpha,\|b\|_{\rm crit})>0$ such that for all $t\in[0,1]$ and $\lambda\in(0,1)$,
$$
|\xi_\lambda(t)|\lesssim_C t+
\lambda^{-1}\int^{\lambda^\alpha t}_0\left|\fint_{B_{2\lambda}}b(s, \cdot)\right|\dif s.
$$
\el
\begin{proof}
By the definition of $\xi_\lambda$ 
and \eqref{DF1}, there is a constant $C_0>0$ such that for all $\lambda\in(0,1)$
\begin{align*}
|\dot \xi_\lambda(t)|&\leq
\lambda^{\alpha-1}\left|\fint_{B_2}(b(\lambda^\alpha t,\lambda(y-\xi_\lambda(t)))-b(\lambda^\alpha t,\lambda y))\dif y\right|
+\lambda^{\alpha-1}\left|\fint_{B_2}b(\lambda^\alpha t,\lambda y)\dif y\right|\\
&\leq C_0\|b\|_{\rm crit}(1+|\xi_\lambda(t)|)+\lambda^{\alpha-1}\left|\fint_{B_2}b(\lambda^\alpha t,\lambda y)\dif y\right|.
\end{align*}
Hence,
\begin{align*}
|\xi_\lambda(t)|\leq C_0\|b\|_{\rm crit}\int^t_0(1+|\xi_\lambda(s)|)\dif s+\lambda^{\alpha-1}\int^t_0\left|\fint_{B_{2\lambda}}b(\lambda^\alpha s, y)\dif y\right|\dif s.
\end{align*}
The desired estimate now follows from the Gronwall's inequality. 
\end{proof}

Now we can give the 

\begin{proof}[Proof of Theorem \ref{Th52}]
In the following we set
$$
u(t,\cdot)|_{t\leq 0}=0,\ \ b(t,\cdot)|_{t\leq 0}=0,\ \ f(t,\cdot)|_{t\leq 0}=0.
$$
Since $\frac{\alpha}{q_0}+|\frac{1}{\bbp_0}|<\alpha$, one can choose $\bar q_0<q_0$ so that
$$
\tfrac{\alpha}{\bar q_0}+|\tfrac{1}{\bbp_0}|<\alpha.
$$
Thus by \eqref{CN1} and Theorem \ref{Th51} with $(\bar q_0,\bbp_0)\in\mI^0_\alpha$, we have for all $t\in[0,T]$,
\begin{align}\label{AAM3}
\|u\cI_t\|_{\mL^\infty_T\mL^\infty_x}\leq C\, \nor f\cI_t\nor_{\wt\mL^{\bar q_0}_T\wt\mL^{\bbp_0}_x}.
\end{align}
Let 
$$
w_0:=u,\ \ v_0:=b,\ \  g_0:=f. 
$$
Fix $\lambda\in(0,1)$, whose value will be determined below. For $k\in\mN$, we inductively define
$$
\dot \eta_k(t)=\lambda^{\alpha-1}\fint_{B_2} v_{k-1}(\lambda^\alpha t, \lambda(y-\eta_k(t)))\dif y
\quad \hbox{with } \eta_k(0)=0,
$$
and for $\Gamma_k(t,x):=(\lambda^\alpha t,\lambda(x-\eta_k(t)))\in\mR^{1+d}$,
$$
w_{k}:=w_{k-1}\circ\Gamma_k,\ \
v_k:=\lambda^{\alpha-1}(v_{k-1}\circ\Gamma_k)-\dot \eta_k,\ \ g_{k}:=\lambda^\alpha (g_{k-1}\circ\Gamma_k).
$$
From the construction and by induction and the chain rule, it is easy to see that
$$
\p_t w_k=\Delta^{\alpha/2}w_k+v_k\cdot\nabla w_k+g_k,
$$
and
\begin{align}
w_n(t,x)=w_0\big(\lambda^{n\alpha}t, \lambda^nx-a_n(t)\big),\ \ 
g_n(t,x)=\lambda^{n\alpha}g_0\big(\lambda^{n\alpha}t, \lambda^nx-a_n(t)\big),\label{ABC2}
\end{align}
where
\begin{align}\label{ABC1}
a_n(t)=\sum_{k=1}^{n}\lambda^{k}\eta_k(\lambda^{(n-k)\alpha}t).
\end{align}
By \eqref{CCA0} with $(v_k, v_{k-1})$ in the role of $(b_\lambda, b)$  and then  \eqref{CCA1}, 
 we have for any $p\geq 1$,
\begin{align}\label{AQ0}
\|v_k\b1_{Q_2}\|_{\mL^\infty_t\mL^{p}_x}\lesssim_C\|v_{k-1}\|_{\rm crit}=\|b_\lambda\|_{\rm crit}=\|b\|_{\rm crit}.
\end{align}
Hence, by Corollary \ref{Co46}, for some $\theta\in(0,1)$ and all $k\in\mN$,
\begin{align}\label{SQ1}
\underset{Q_{1/2}}{\osc}(w_k)\leq \theta\,\underset{Q_6}{\osc}(w_k)+\|w_k\cI_4\|_\infty+\|g_k\b1_{Q_4}\|_{\mL^{q_0}_t\mL^{\bbp_0}_x}.
\end{align}
Note that for any $p\geq 1$ and $\lambda\in(0,1)$,
\begin{align*}
\left|\fint_{B_{2\lambda}}v_{k-1}(t,y)\dif y\right|&\leq \left(\fint_{B_{2\lambda}}|v_{k-1}(t,y)|^p\dif y\right)^{1/p}\\
&=\lambda^{-d/p} \left(\frac{1}{|B_2|}\int_{B_{2\lambda}}|v_{k-1}(t,y)|^p\dif y\right)^{1/p}\\
&\leq\lambda^{-d/p} \left(\frac{1}{|B_2|}\int_{B_{2}}|v_{k-1}(t,y)|^p\dif y\right)^{1/p}\\
&\lesssim_C\lambda^{-d/p}\|v_{k-1}\|_{\rm crit}=\lambda^{-d/p}\|b\|_{\rm crit},
\end{align*}
where in the last inequality we have used that
$\int_{B_{2}}v_{k-1}(t,y)\dif y=0$ and \eqref{59} for $\alpha=1$, and the constant $C$ only depends on $d,\alpha,p$.
By Lemma \ref{Le54} and \eqref{AQ0}, for any $p\geq 1$, there is a $C=C(p,d,\alpha, \|b\|_{\rm crit})>0$ such that for all $k\in\mN$ and $\lambda, t\in(0,1)$,
\begin{align}
|\eta_k(t)|
&\lesssim_C t+\lambda^{-1}\int^{\lambda^{\alpha}t}_0\left|\fint_{B_{2\lambda}}v_{k-1}(s,y)\dif y\right|\dif s\no\\
&\lesssim_C t+\lambda^{\alpha-1-\frac d p}t\|v_{k-1}\|_{\rm crit}
\lesssim_C(1+\lambda^{\alpha-1-\frac d p})t.\label{AQ2}
\end{align}
In particular, for fixed $p>d/\alpha$, one can choose $\lambda_0$ small enough so that for all $\lambda\leq\lambda_0$ and $k\in\mN$,
$$
(t,x)\in Q_6\Rightarrow\Gamma_k(t,x)=(\lambda t^\alpha, \lambda(x-\eta_k(t)))\in Q_{1/2}.
$$
Thus, 
$$
\underset{Q_6}{\osc}(w_k)=\underset{Q_6}{\osc}(w_{k-1}\circ\Gamma_k)\leq \underset{Q_{1/2}}{\osc}(w_{k-1}),
$$ 
and by \eqref{SQ1},
\begin{align*}
\underset{Q_{1/2}}{\osc}(w_k)
&\leq  \theta\,\underset{Q_{1/2}}{\osc}(w_{k-1})+\|w_k\cI_4\|_\infty+\|g_k\b1_{Q_4}\|_{\mL^{q_0}_t\mL^{\bbp_0}_x}.
\end{align*}
Iterating the above inequality yields
\begin{align}\label{AQ5}
\underset{Q_{1/2}}{\osc}(w_n)\leq \theta^n\underset{Q_{1/2}}{\osc}(w_0)+\sum_{k=1}^n\theta^{n-k}
\Big(\|w_k\cI_4\|_\infty+\|g_k\b1_{Q_4}\|_{\mL^{q_0}_t\mL^{\bbp_0}_x}\Big).
\end{align}
By \eqref{ABC1} and \eqref{AQ2} we have for all $\lambda,t\in(0,1)$,
\begin{align*}
a_n(t)&\lesssim_C \sum_{k=1}^{n}\lambda^{k}(1+\lambda^{\alpha-1-\frac d p})\lambda^{(n-k)\alpha}t\\
&=(1+\lambda^{\alpha-1-\frac d p})\lambda^{n\alpha}\sum_{k=1}^{n}\lambda^{(1-\alpha)k}t\\
&\leq 2\lambda^{n\alpha-\frac dp}nt.
\end{align*}
From this we can deduce
that for some $n_0=n_0(\lambda)$ large enough and all $n\geq n_0$,
\begin{align}\label{AQ3}
Q_{(\lambda/2)^n}\subset\big\{(s,y): s=\lambda^{n\alpha}t, y=\lambda^n x-a_n(t), (t,x)\in Q_{1/2}\big\}.
\end{align}
Indeed, for $|s|\leq(\lambda/2)^n$ and $|y|\leq (\lambda/2)^n$, we have
\begin{align*}
|t|=\lambda^{-n\alpha}|s|\leq\tfrac{\lambda^{(1-\alpha)n}}{2^{n}}\leq \tfrac12,
\end{align*}
and for some $n_0=n_0(\lambda)$ large enough and all $n\geq n_0$,
\begin{align*}
|x|\leq \lambda^{-n}(|y|+|a_n(t)|)\leq\tfrac{1}{2^{n}}+\tfrac{C\lambda^{-\frac dp}n}{2^{n}}\leq\tfrac12.
\end{align*}
This establishes \eqref{AQ3}. Thus, by \eqref{ABC2}, we have
 for $n\geq n_0$,
\begin{align}\label{AQ6}
\underset{Q_{(\lambda/2)^{n}}}{\osc}(w_0)\leq \underset{Q_{1/2}}{\osc}(w_n).
\end{align}
Moreover, by \eqref{ABC2} and  H\"older's inequality, we have
\begin{align}
\|w_k\cI_4\|_\infty&=\|w_0(\lambda^{k\alpha}\cdot,\cdot)\cI_4\|_\infty=\|u\cI_{4\lambda^{k\alpha}}\|_\infty
\stackrel{\eqref{AAM3}}{\lesssim}\nor f\cI_{4\lambda^{k\alpha}}\nor_{\wt\mL^{\bar q_0}_T\wt\mL^{\bbp_0}_x}\no\\
&\lesssim \lambda^{k\alpha(\frac 1{\bar q_0}-\frac 1{q_0})}\nor f\cI_{4\lambda^{k\alpha}}\nor_{\wt\mL^{q_0}_T\wt\mL^{\bbp_0}_x}
\leq \lambda^{k(\frac \alpha{\bar q_0}-\frac \alpha{q_0})}\nor f\nor_{\wt\mL^{q_0}_T\wt\mL^{\bbp_0}_x},
\end{align}
and by the change of variable,
\begin{align}
\|g_k\b1_{Q_4}\|_{\mL^{q_0}_t\mL^{\bbp_0}_x}
&=\lambda^{k\alpha}
\left(\int_\mR\|f(\lambda^\alpha t,\lambda^k(\cdot-a_k(t)))\b1_{Q_4}\|^{q_0}_{\mL_x^{\bbp_0}}\dif t\right)^{1/q_0}\no\\
&\leq\lambda^{k(\alpha-\frac{\alpha}{q_0}-|\frac{1}{\bbp_0}|)}\sup_{z\in\mR^d}
\left(\int_0^{4\lambda^\alpha}\|f( t,\cdot-z)\b1_{B_{4\lambda^k}}\|^{q_0}_{\mL_x^{\bbp_0}}\dif t\right)^{1/q_0}\no\\
&\leq\lambda^{k(\alpha-\frac{\alpha}{q_0}-|\frac{1}{\bbp_0}|)}\nor f\nor_{\wt\mL^{q_0}_T\wt\mL^{\bbp_0}_x}.\label{AQ7}
\end{align}
Combining \eqref{AQ5}, \eqref{AQ6}-\eqref{AQ7}, we obtain  for $\delta:=(\alpha-\frac{\alpha}{q_0}-|\frac{1}{\bbp_0}|)\wedge (\frac \alpha{\bar q_0}-\frac \alpha{q_0})$
and any $n\geq n_0$,
\begin{align*}
\underset{Q_{(\lambda/2)^{n}}}{\osc}(w_0)
&\leq\theta^n\underset{Q_{1/2}}{\osc}(w_0)+C\sum_{k=1}^n\theta^{n-k}\lambda^{k\delta}\nor f\nor_{\wt\mL^{q_0}_T\wt\mL^{\bbp_0}_x}\\
&=\theta^n\underset{Q_{1/2}}{\osc}(w_0)+C\theta^n\sum_{k=1}^n(\lambda^\delta/\theta)^k\nor f\nor_{\wt\mL^{q_0}_T\wt\mL^{\bbp_0}_x}.
\end{align*}
In particular, if $\lambda\leq (\theta/2)^{1/\delta}\wedge\lambda_0$, then by \eqref{AAM3}, for some $C_0>0$ and any $n\geq n_0$,
$$
\underset{Q_{(\lambda/2)^{n}}}{\osc}(u)=\underset{Q_{(\lambda/2)^{n}}}{\osc}(w_0)
\leq\theta^n\left[\underset{Q_{1/2}}{\osc}(w_0)+\frac{C\lambda^\delta}{\theta-\lambda^\delta}
\nor f\nor_{\wt\mL^{q_0}_T\wt\mL^{\bbp_0}_x}\right]\leq C_0\theta^n\nor f\nor_{\wt\mL^{q_0}_T\wt\mL^{\bbp_0}_x}.
$$
From this and by shifting, it is standard to derive the desired H\"older continuity.
\end{proof}
 
The following stability result will be used to establish 
the well-posedness of generalized martingale problem in the critical case.

\bt\label{Th42}
Let $T>0$ and $(b_n,f_n)$ be a sequence of 
functions  satisfying 
\begin{enumerate}[\rm (i)]
\item $b_n$ satisfies \eqref{CN1} uniformly in $n$ and $\sup_n\nor f_n\nor_{\wt\mL^{q_0}_T\wt\mL^{\bbp_0}_x}<\infty$
for some $(q_0,\bbp_0)\in\mI^0_\alpha$.
\item  $(b_n,\div b_n, f_n)$ converges to $(b,\div b, f)$ in Lebesgue measure as $n\to\infty$. 
\end{enumerate}
Let $u_n$ and $u$ be the respective weak solutions of PDE \eqref{PDE0} corresponding to $(b_n,f_n)$ and $(b,f)$ with zero initial value.
Then for any bounded domain $Q\subset[0,T]\times\mR^d$,
\begin{align}\label{CC41}
\lim_{n\to\infty}\sup_{(t,z)\in Q}|u_{n}(t,x)-u(t,x)|=0.
\end{align}
\et
\begin{proof}
First of all, by (i) and \eqref{EA80} we have
\begin{align}\label{CC4}
\sup_n\|u_n\|_{\mL^\infty_T\mL^\infty_x}<\infty.
\end{align}
By  Theorem \ref{Th52}, there are constants $\gamma, C>0$ such that for all $s,t\in[0,T]$ and $x,y\in\mR^d$,
$$
\sup_{n\geq 1} |u_n(t,x)-u_n(s,y)|\leq C(|t-s|+|x-y|)^\gamma.
$$
In particular, by \eqref{CC4}  and Ascolli-Arzela's theorem, 
there is a subsequence $n_k$ and 
$$
\bar u\in C_b([0,T]\times\mR^{d})
$$ 
such that  for any bounded domain $Q\subset[0,T]\times\mR^d$,
$$
\lim_{k\to\infty}\sup_{(t,z)\in Q}|u_{n_k}(t,z)-\bar u(t,z)|=0.
$$
Since for any $\varphi\in C^\infty_c((0,T)\times\mR^d)$,
$$
-\int_{\mR^{1+d}} u_n\p_t\varphi =
\int_{\mR^{1+d}} u_n\Delta^{\alpha/2}\varphi
-\int_{\mR^{1+d}}(b_n\cdot\nabla \varphi+\div b_n\,\varphi)u_n+\int_{\mR^{1+d}} f_n\varphi.
$$
By (ii) and taking limits, one sees that $\bar u=u$ is the unique weak solution of PDE \eqref{PDE0} with zero initial value.
By a contradiction method, one sees that \eqref{CC41} holds.
\end{proof}

\section{Critical stochastic quasi-geostrophic equations with additive noises}

We consider the following stochastic quasi-geostropic equation in $\mR^2$:
\begin{align}\label{SPDE0}
\dif u=\big( \Delta^{1/2}u+\cR u\cdot\nabla u \big) \dif t+\sum_{k\in\mN}g_k\dif W^k_t,\ \ u(0)=u_0,
\end{align}
where  $W=(W^k)_{k\in\mN}$ is a sequence of independent   one-dimensional standard Brownian motions,  
 $g=(g_k)_{k\in\mN}$ is a sequence of space-time functions in some Sobolev spaces, and
$$
\cR u:=(-\p_2\Delta^{-1/2} u, \p_1\Delta^{-1/2}u). 
$$
It is well known that for any $p\in(1,\infty)$ (cf. \cite{G19}),
\begin{align}\label{RZ1}
\|\cR u\|_{L^p}\leq C\|u\|_{L^p},\ \ \|\cR u\|_{\rm BMO}\leq C\|u\|_{\mL^\infty},
\end{align}
and for any $\gamma\in(0,1)$,
\begin{align}\label{RZ2}
\|\cR u\|_{\cC^\gamma}\leq C\|u\|_{\cC^\gamma}.
\end{align}
We assume  that for some  $\beta\in(\frac12,1]$ and any $T>0$,
\begin{align}\label{REG1}
u_0\in \bH^1_2\cap \bH^1_4
\quad \hbox{and} \quad 
 g\in \mL^\infty_T(\bH^\beta_2(\ell^2)\cap \bH^\beta_4(\ell^2)).
\end{align}
Here we do not pursue the best possible
conditions on $u_0$ and $g$, but concentrate on how to use Theorem  \ref{Th52} to 
establish the pathwise uniqueness of SPDE \eqref{SPDE0}.

\bd\label{Def1}
Let ${\frak F}:=(\Omega,\sF,\bP;(\sF_t)_{t\geq 0})$ be a stochastic basis. Let $u:[0,T]\times\Omega\to \mL^2_x$ be a progressively measurable process
and $W = (W^1, W^2, \cdots )$ a sequence of i.i.d. Brownian motions. We say $({\frak F}, u,W)$  is 
a weak solution of SPDE \eqref{SPDE0} if
\begin{enumerate} [\rm (i)] 
\item For $\bP$-almost all $\omega$, it holds that for any $T>0$,
$$
u(\cdot,\omega)\in C([0,T]; L^2_w)\cap L^2(0,T; \bH^{1/2}_2),
$$
where $L^2_w$ is the usual $L^2$-space endowed with the weak topology.
\item  For any $\phi\in C^\infty_c( \mR^2)$, it holds that for all $t>0$, 
$$
\<u(t),\phi\>=\<u(0),\phi\>+\int^t_0\<u(s),\Delta^{1/2}\phi\>\dif s+\int^t_0\<u(s),\cR u(s)\cdot\nabla\phi\>\dif s+\sum_k\int^t_0\<g_k(s),\phi\>\dif W^k_s.
$$
\end{enumerate} 
If in addition, 
$u$ is progressively measurable with respect to the 
augmented minimum filtration $\{ \sF^W_t\}_{t\geq 0}$ generated by $W$, 
then we say $u$ a strong solution to  \eqref{SPDE0}. 
 \ed

The following is the main result of this section. 

\bt\label{Th62}
Under \eqref{REG1}, there is a unique strong solution for 
stochastic quasi-geostropic equation \eqref{SPDE0}  with regularities that for any $p,q\geq 2$ and $T>0$,
\begin{align}\label{REG4}
\bE\left[ \|u\|^q_{\mL^\infty_T\mL^p_x}\right]+
\bE\left[ \|u\|^2_{\mL^2_T\bH^{1/2}_2}\right] <\infty.
\end{align}
\et
Note that by Sobolev's embedding, under \eqref{REG1}, we have for any $p\geq 2$,
$$
u_0\in \mL^p_x,\ \ g\in L^\infty([0,T];\mL^2_x(\ell_2)\cap \mL^p_x(\ell^2)).
$$
Thus, the existence of a weak solution with estimate \eqref{REG4} has been shown in \cite[Theorem 3.7]{BM19}. 
By Yamada-Watanabe's theorem, our main task is to establish the pathwise uniqueness.
To this aim, we  first establish
  some a priori regularity estimates for any weak solution $u$. First of all, we study the following stochastic convolution:
$$
\mS_g(t):=\sum_{k\in\mN}\int^t_0P_{t-s}g_k(s)\dif W^k_s,
$$
where  $(P_t)_{t\geq 0}$ 
is the transition  semigroup of the isotropic Cauchy process, or equivalently, of $\Delta^{1/2}$  on $\R^2$.
 Note that for any $p\in(1,\infty)$ and $\beta\in[0,1]$,
\begin{align}\label{SM3}
\|\nabla P_t g\|_{\mL^p_x(\ell^2)}\leq Ct^{\beta-1}\|\Delta^{\beta/2}g\|_{\mL^p_x(\ell^2)},\ \ t>0.
\end{align}
We need the following estimate about stochastic convolution $\mS_g(t)$ due to \cite{DKZ87}. 
\bl \label{37}
Let $p,q\geq2$ and $\beta\in(\frac12+\frac1q,1]$. For any $T>0$ and $g\in\mL^q_T(\bH^\beta_p(\ell^2))$, there  is a continuous version $\wt\mS_g$ of
$[0,T]\ni t\mapsto \mS_g(t):=\int^t_0P_{t-s}g_k(s)\dif W^k_s\in \bH^1_p$ so that
for some $C=C(q,p,T,\beta)>0$,
\begin{align}\label{323}
\mE\left[ \sup_{t\in[0,T]}\left\|\wt\mS_g(t)\right\|^q_{\bH^1_p}\right] 
\lesssim_C\int^T_0\|g(s)\|^q_{\bH^\beta_p(\ell^2)}\dif s.
\end{align}
\el
\begin{proof}
We use the factorization method in \cite{DKZ87}. Noting that
$$
\int^t_s(t-r)^{\kappa-1}(r-s)^{-\kappa}\dif r=\frac{\pi}{\sin(\kappa\pi)},\ \ 0\leq s<t<\infty,\ \kappa\in(0,1),
$$
by stochastic Fubini's theorem, we have  for each $t\in(0,T)$,
$$
\mS_g(t)=\sum_{k\in\mN}\int^t_0P_{t-s}g_k(s)\dif W^k_s=\frac{\sin(\kappa \pi)}{\pi}\int^t_0(t-s)^{\kappa-1}P_{t-s}G(s)\dif s=:\wt\mS_g(t),\ \mP-a.s.
$$
where
$$
G(t):=\sum_{k\in\mN}\int^t_0(t-s)^{-\kappa}P_{t-s}g_k(s)\dif W^k_s.
$$
Let $\kappa\in(\beta-\frac12-\frac1q,\beta-\frac12)$. By Minkowskii's inequality and H\"older's inequality, we have
\begin{align}
\left\|\nabla\wt\mS_g(t)\right\|_{\mL^p_x}
&\lesssim\int^t_0(t-s)^{\kappa-1}\|\nabla P_{t-s}G(s)\|_{\mL^p_x}\dif s
\leq\int^t_0(t-s)^{\kappa-1}\|\nabla G(s)\|_{\mL^p_x}\dif s\no\\
&\leq\left(\int^t_0(t-s)^{(\kappa-1)q/(q-1)}\dif s\right)^{1-1/q}\left(\int^t_0\|\nabla G(s)\|^q_{\mL^p_x}\dif s\right)^{1/q}.\label{AB1}
\end{align}
From this a priori estimate, as in \cite[Lemma 1]{DKZ87}, one sees that $t\mapsto \nabla\wt\mS_g(t)$ is continuous in $\mL^p_x$
as long as the last integral is finite a.s. Indeed, by BDG's inequality and Fubini's theorem, we have
\begin{align*}
\int^T_0\mE \left\|\nabla G(t)\right\|_{\mL^p_x}^q\dif t
&\lesssim\int^T_0\left(\int^t_0(t-s)^{-2\kappa}\|\nabla P_{t-s}g(s)\|^2_{\mL^p_x(\ell_2)}\dif s\right)^{q/2}\dif t \\
&\lesssim\int^T_0 \left(\int^t_0(t-s)^{-2(1-\beta+\kappa)}\|\Delta^{\beta/2}g(s)\|^2_{\mL^p_x(\ell^2)}\dif s\right)^{q/2}\dif t \\
&\lesssim \int^T_0 \int^t_0(t-s)^{-2(1-\beta+\kappa)}\|\Delta^{\beta/2}g(s)\|^2_{\mL^p_x(\ell^2)}\dif s\dif t\\
&= \int^T_0 \left(\int^T_s(t-s)^{-2(1-\beta+\kappa)}\dif t\right)\|\Delta^{\beta/2}g(s)\|^2_{\mL^p_x(\ell^2)}\dif s\\
&\lesssim T^{2(\beta-\kappa)-1}\int^T_0\|\Delta^{\beta/2}g(s)\|^2_{\mL^p_x(\ell^2)}\dif s.
\end{align*}
This together with \eqref{AB1} yields 
$$
\mE\left(\sup_{t\in[0,T]}\left\|\nabla\wt\mS_g(t)\right\|^q_{\mL^p_x}\right)\lesssim_C\int^T_0\|\Delta^{\beta/2}g(s)\|^q_{\mL^p_x(\ell^2)}\dif s.
$$
Similarly, one can show
$$
\mE\left(\sup_{t\in[0,T]}\left\|\wt\mS_g(t)\right\|^q_{\mL^p_x}\right)\lesssim_C\int^T_0\|g(s)\|^q_{\mL^p_x(\ell^2)}\dif s.
$$
This completes the proof.
\end{proof}

\bl
Under \eqref{REG1}, for any weak solution $u$ of SQG \eqref{SPDE0},  $T>0$ and  $p\in(2,4)$, it holds that
\begin{align}\label{AQ9}
u\in\mL^p_T(\bH^1_p),\ \ a.s.
\end{align}
In particular, $u$ is a strong solution of SQG \eqref{SPDE0}  in PDE's sense.
\el

\begin{proof}
Define
$$
\wt u(t):=P_t u_0+\mS_g(t).
$$
By \eqref{323} and $\|P_tu_0\|_{\bH^1_4}\leq \|u_0\|_{\bH^1_4}$, for any $q\geq 2$ we have
\begin{align}
\mE\left[ \|\wt u\|_{\mL^\infty_T(\cC^{1/2}_x)}^q\right]
\lesssim\mE\left[ \|\wt u\|_{\mL^\infty_T(\bH^1_4)}^q\right]
\lesssim \|u_0\|_{\bH^1_4}+\int^T_0\|g(s)\|^q_{\bH^\beta_4(\ell^2)}\dif s<\infty.
\end{align}
Clearly, $\wt u$ solves the  SPDE:
$$
\dif \wt u= \Delta^{1/2}\wt u \dif t+g_k\dif W^k_t
\quad \hbox{with} \quad \wt u(0)=u_0.
$$
Let $\bar u:=u-\wt u$. It is easy to see that
$$
\p_t \bar u= \Delta^{1/2}\bar u+\cR u\cdot\nabla\bar u+\cR u\cdot\nabla \wt u
\quad \hbox{with} \quad \bar u(0)=0.
$$
Note that for any $p>1$,
$$
\div \cR u=0 \quad \hbox{and} \quad 
\|\cR u\|_{\mL^{\infty}_t\mL^{p}_x}\leq C\|u\|_{\mL^{\infty}_t\mL^{p}_x}.
$$
By \eqref{EA80} with $b=\cR u$ and $f=\cR u\cdot\nabla \wt u$, for $p\in(2,4)$, there is a constant
 $\gamma\in(0,1/4)$ such that
\begin{align*}
\|\bar u\|_{\mL^\infty_T(\cC^\gamma_x)}\lesssim\|\cR u\cdot\nabla \wt u\|_{\mL^\infty_T\mL^{p}_x}
&\lesssim \|\cR u\|_{\mL^\infty_T\mL^{4p/(4-p)}_x}\|\nabla \wt u\|_{\mL^\infty_T\mL^4_x}<\infty.
\end{align*}
This together with $\|\wt u\|_{\mL^\infty_T(\cC^{1/4}_x)}<\infty$ implies that
\begin{align*}
\|\cR u\|_{\mL^\infty_T(\cC^\gamma_x)}
&\lesssim \|u\|_{\mL^\infty_T(\cC^\gamma_x)}\leq
 \|\bar u\|_{\mL^\infty_T(\cC^\gamma_x)}+\|\wt u\|_{\mL^\infty_T(\cC^\gamma_x)}<\infty.
\end{align*}
Thus for $p\in(2,4)$, by \cite[Lemma 5.1]{Zh13} with $b=\cR u$ and $f=\cR u\cdot\nabla \wt u$,
\begin{align*}
\|\bar u\|_{\mL^p_T(\bH^1_p)}&\lesssim \|\cR u\cdot\nabla \wt u\|_{\mL^p_T\mL^p_x}
\leq T^{1/p}\|\cR u\cdot\nabla \wt u\|_{\mL^\infty_T\mL^p_x}<\infty.
\end{align*}
The  desired result now follows from 
$u=\bar u+\wt u$ and $\wt u\in\mL^\infty_T(\bH^1_p)$ a.s.
\end{proof}
Now we can give the
\begin{proof}[Proof of Theorem \ref{Th62}]
Let $u_1,u_2$ be two weak solutions of SQG \eqref{SPDE0} defined on the same probability space with the same initial value. 
Define $U=u_1-u_2$. Then $U(0)=0$ and
$$
\p_tU=\Delta^{ 1/2}U+\cR u_1\cdot\nabla U+\cR U\cdot\nabla u_2.
$$
Let $\|f\|_{\cC^\gamma}:=\|f\|_\infty+\|f\|_{\dot \cC^\gamma}$ and $T>0$.
For $p\in(3,4)$, by Theorems \ref{Th51} and \ref{Th52}, there are constants $\gamma\in(0,1)$ and $C>0$ 
such that for any $t\in[0,T]$,
\begin{align*}
\|U(t)\|^{p}_{\cC^\gamma}
&\lesssim_C \int^t_0\|\cR U(s)\cdot\nabla u_2(s)\|^{p}_{\mL^{p}_x}\dif s\\
&\lesssim_C \int^t_0\|\cR U(s)\|^{p}_{\mL^\infty_x}\|\nabla u_2(s)\|^{p}_{\mL^{p}_x}\dif s\\
&\lesssim_C \int^t_0\|\cR U(s)\|^{p}_{\cC^\gamma_x}\|\nabla u_2(s)\|^{p}_{\mL^{p}_x}\dif s\\
&\lesssim_C \int^t_0\|U(s)\|^{p}_{\cC^\gamma_x}\|\nabla u_2(s)\|^{p}_{\mL^{p}_x}\dif s.
\end{align*}
By \eqref{AQ9} and Gronwall's inequality, 
$$
\|U(t)\|^{p}_{\cC^\gamma}=0.
$$
Thus, the pathwise uniqueness holds.  The proof is now complete.
\end{proof}

\section{Existence of weak solutions to DDSDE}\label{Sec7}

In this section we study the following distribution-dependent SDE (DDSDE):
\begin{align}\label{SDE1b}
\dif X_t=\int_{\mR^d}b(t,X_t, y)\mu_{X_t}(\dif y)\dif t+\dif L^{(\alpha)}_t,
\end{align}
where $b(t,x,y): \mR_+\times\mR^d\times\mR^d\to\mR^d$ is a measurable function and $\mu_{X_t}$ stands for the law of $X_t$ on $\R^d$ with $d\geq 1$.
In the following, for a Banach space $\mB$, we shall use $\cP(\mB)$ to denote the set of all probability measures 
on $\mB$ over Borel $\sigma$-field $\cB(\mB)$.

\bd\label{Def2}
Let $\nu\in\cP(\mR^{d})$ and ${\frak F}:=(\Omega,\sF,\bP;(\sF_t)_{t\geq 0})$ be a stochastic basis, $X_t$ 
an $\mR^{d}$-valued $\sF_t$-adapted process, $L^{(\alpha)}$ a $d$-dimensional  isotropic
 $\alpha$-stable process.
We call $({\frak F}, X,L^{(\alpha)})$ a weak solution of DDSDE \eqref{SDE1b} with initial distribution $\nu$ if $\bP\circ X_0^{-1}=\nu$
and for all $t\geq 0$, 
$$
X_t=X_0+\int^t_0\int_{\mR^d}b(s, X_s,y)\mu_{X_s}(\dif y)\dif s+L^{(\alpha)}_t
\quad  \bP \hbox{-a.s.},
$$
where $\mu_{X_s}$  is the law of $X_s$.
\ed

The following is the main result of this section.

 \bt\label{Th55}
 Let $T>0$. 
Suppose that $\div_x b\equiv 0$ and for all $(t,x,y)\in[0,T]\times\mR^{d}\times\mR^{d}$,
\begin{align}\label{DR1}
|b(t,x,y)|\leq h(t,x-y) \quad \mbox{with } \quad \nor h\nor_{\widetilde\mL^{q_1}_T(\widetilde\mL^{\bbp_1}_x)}\leq\kappa_0,
\end{align}
where $(q_1,\bbp_1)\in\mI^0_\alpha$.
Then for any $\nu\in\cP(\mR^{d})$, 
 there exists at least one weak solution to DDSDE \eqref{SDE1b} with initial distribution $\nu$ in the sense of Definition \ref{Def2}.
 Moreover, for Lebesgue almost every  $t\geq 0$,
  $X_t$ admits a density $\rho(t,\cdot)$
which enjoys the following regularity: for any $\beta\in[0,\frac\alpha2)$ and $(q,\bbp)\in(1,\infty)^{1+d}$ with $\tfrac\alpha{q}+|\tfrac1{\bbp}|>d+\beta$,
\begin{align}\label{DL131}
\|\rho\|_{\mL^{q}_T(\bH^{\beta}_{\bbp})}<\infty.
\end{align}
 \et

To show the existence of a weak solution, we use the mollifying technique. 
Fix $T>0$. Let $\varGamma\in C^\infty_c(\mR^{1+2d})$ be a smooth probability density function with compact support and define 
for $n\in\mN$,  
$$
\varGamma_n(t,x,y)=n^{1+2d}\varGamma(nt,nx,ny)
$$ 
and 
\begin{align*}
b_n(t,x,y)&:=((b\b1_{[0,T]})*\varGamma_n)(t,x,y).
\end{align*}
Note that by \eqref{DR1}, 
$$
\|\nabla^j_x \nabla^k_yb_n\|_\infty<\infty,\ \ j,k\in\mN_0.
$$ 
By this, it is well known
that the following DDSDE admits a unique strong solution (cf. \cite{Szn91}):
\begin{align}\label{SDE19}
X^n_t=X_0+\int^t_0\int_{\mR^d}b_n(s, X^n_s,y)\mu_{X^n_s}(\dif y)\dif s+L^{(\alpha)}_t,
\end{align}
where
$$
\bP\circ X_0^{-1}=\nu.
$$
Now we use Theorem \ref{Th51} to derive the following crucial Krylov estimate.

\bl [Krylov's estimates]\label{Kry0}
For any $\beta\in[0,\frac\alpha2)$ and $(q_0,\bbp_0)\in\mI^\beta_\alpha$, there are $\theta\in(0,1)$ and constant $C=C(T,d,\alpha,\beta,\kappa_0,q_i,\bbp_i)>0$ 
such that for any $\delta\in(0,T)$, stopping time $\tau\leq T-\delta$ and  $f\in C^\infty_c(\mR^{1+d})$,
\begin{align}\label{Kr0}
\sup_{n\in\mN}\bE\left[ \int^{\tau+\delta}_{\tau}f(r,X^n_r)\dif r\Big|\sF_{\tau}\right]
\leq C\, \delta^\theta\nor f\nor_{\wt\mL^{q_0}_T(\wt\bH^{-\beta}_{\bbp_0})}.
\end{align}
In particular,  for Lebesgue almost every  $t\geq 0$,
 $X^n_t$ admits a density $\rho_n(t,\cdot)$ 
such  that for any $\beta\in[0,\frac\alpha2)$ and $(q,\bbp)\in(1,\infty)^{1+d}$ with $\tfrac\alpha{q}+|\tfrac1{\bbp}|>d+\beta$,
\begin{align}\label{DL135}
\sup_n\|\rho_n\|_{\mL^{q}_T(\bH^{\beta}_{\bbp})}<\infty.
\end{align}
\el
\begin{proof}
For $n\in\mN$ and $(t,x)\in[0,T]\times\mR^d$, define
\begin{align}\label{ABC5}
\wt b_n(t,x):=\int_{\mR^d} b_n(t,x,y)\mu_{X^n_t}(\dif y).
\end{align}
Note that by \eqref{DR1} and the definition \eqref{GG1},
\begin{align}\label{ABC3}
\nor \wt b_n\nor_{\widetilde\mL^{q_1}_T(\widetilde\mL^{\bbp_1}_x)}\leq \nor h\nor_{\widetilde\mL^{q_1}_T(\widetilde\mL^{\bbp_1}_x)}\leq\kappa_0,
\end{align}
and for any  $j\in\mN_0$,
$$
\|\nabla^j\wt b_n\|_\infty<\infty.
$$
Fix $0\leq t_0<t_1\leq T$. For given $f\in C^\infty_c(\mR^{1+d})$, 
let $u_n\in L^\infty_{t_1}(C^\infty_b(\mR^{d}))$ solve the following backward Kolmogorov equation:
\begin{align}\label{DC10}
\p_t u_n+\Delta^{\alpha/2} u_n+\wt b_n\cdot\nabla u_n=f,\ \ u_n(t_1)=0.
\end{align}
Since $|\tfrac1{\bbp_0}|+\tfrac\alpha {q_0}<\alpha-\beta$, one can choose $\bar q_0< q_0$ so that
$$
|\tfrac1{\bbp_0}|+\tfrac\alpha {\bar q_0}<\alpha-\beta.
$$
Thus by $\div\wt b_n=0$, \eqref{ABC3} and \eqref{EA80}, there is a constant $C=C(T,d,\alpha,\beta,\kappa_0,q_i,\bbp_i)>0$ such that
\begin{align}\label{DC0}
\sup_{n\in\mN}\|\1_{[t_0,t_1]}u_n\|_\infty\lesssim_C\nor\1_{[t_0,t_1]}f\nor_{\wt\mL^{\bar q_0}_T(\wt\bH^{-\beta_0}_{\bbp})}
\lesssim_C(t_1-t_0)^{\frac1{\bar q_0}-\frac1{ q_0}}\nor\1_{[t_0,t_1]}f\nor_{\wt\mL^{q_0}_T(\wt\bH^{-\beta}_{\bbp_0})},
\end{align}
where the second step is due to H\"older's inequality.
Now by \eqref{DC10} and It\^o's formula, we have
\begin{align*}
0=u_n(t_1, X^n_{t_1})=u_n(t_0, X^n_{t_0})+\int^{t_1}_{t_0}f(r,X^n_r)\dif r+M_{t_1}-M_{t_0},
\end{align*}
where $t\mapsto M_t$ is a local martingale.
Hence, by taking conditional expectation with respect to $\sF_{t_0}$ and \eqref{DC0},
\begin{align*}
\bE\left[ \int^{t_1}_{t_0} f(r,X^n_r)\dif r\Big|\sF_{t_0}\right] =-u_n(t_0, X^n_{t_0})
\leq C\, (t_1-t_0)^{\frac1{\bar q_0}-\frac1{ q_0}}\nor\1_{[t_0,t_1]}f\nor_{\wt\mL^{q_0}_T(\wt\bH^{-\beta}_{\bbp_0})}.
\end{align*}
By discretization stopping time approximation (see \cite[Remark 1.2]{ZZ18}), we obtain \eqref{Kr0}.
As for \eqref{DL135}, it follows by
$$
\sup_{n\in\mN}\bE\left[ \int^T_0f(r,X^n_r)\dif r\right]
\leq C\, \delta^\theta\nor f\nor_{\wt\mL^{q_0}_T(\wt\bH^{-\beta}_{\bbp_0})}\leq C\, \delta^\theta\| f\|_{\mL^{q_0}_T(\bH^{-\beta}_{\bbp_0})},
$$
and a standard duality method.
\end{proof}

Let $\mD$ be the space of all c\'adl\'ag functions from $[0,T]$ to $\mR^d$. We use $(\omega_t)_{t\in[0,T]}$ to denote a generic path in $\mD$.
By Krylov's estimate \eqref{Kr0}, we show the following tightness result.

\bl\label{Le56}
The law $\mP_n$ of $X^n$ in $\mD$ is tight. For any accumulation point $\mP$, it holds that for 
any $\beta\in[0,\frac\alpha2)$ and $(q_0,\bbp_0)\in\mI^\beta_\alpha$,  and  $f\in C^\infty_c(\mR^{1+d})$,
\begin{align}\label{ACV6}
\mE^\mP\left[ \int^T_0f(r,\omega_r)\dif r\right]
\lesssim_C\nor f\nor_{\wt\mL^{q_0}_T(\wt\bH^{-\beta}_{\bbp_0})}.
\end{align}
In particular, for Lebesgue almost every $t\in[0,T]$, $\mP\circ\omega_t^{-1}(\dif x)=\rho(t,x)\dif x$.
\el

\begin{proof}
Let $\delta<T$ and $\tau$ be any stopping time bounded by $T-\delta$. By SDE \eqref{SDE19}, we have
\begin{align*}
X^n_{\tau+\delta}-X^n_{\tau}=\int^{\tau+\delta}_\tau \wt b_n(r,X^n_r)\dif r+L^{(\alpha)}_{\tau+\delta}-L^{(\alpha)}_\tau,
\end{align*}
where $\wt b_n$ is defined by \eqref{ABC5}.
By Krylov's estimate \eqref{Kr0} and \eqref{ABC3}, there is a $\theta\in(0,1)$ such that
\begin{align}\label{ABC6}
\bE\left[ \int^{\tau+\delta}_\tau |\wt b_n(r,X^n_r)|\dif r\right] 
\lesssim\delta^\theta\nor \wt b_n\nor_{\widetilde\mL^{ q_1}_T(\widetilde\mL^{\bbp_1}_x)}\lesssim\delta^\theta\nor h\nor_{\widetilde\mL^{ q_1}_T(\widetilde\mL^{\bbp_1}_x)}.
\end{align}
Moreover, by the property of L\'evy processes, for any $\eps>0$,
\begin{align*}
\lim_{\delta\to0}\sup_{\tau}\bP\left(\sup_{t\in[0,\delta]}|L^{(\alpha)}_{\tau+t}-L^{(\alpha)}_\tau|\geq\eps\right)=0,
\end{align*}
which together with \eqref{ABC6} yields that for any $\eps>0$,
\begin{align*}
\lim_{\delta\to 0}\sup_{n}\sup_\tau\bP\left(\sup_{t\in[0,\delta]}|X^n_{\tau+t}-X^n_{\tau}|>\eps\right)=0.
\end{align*}
By Aldous' criterion (see \cite[p.356, Theorem 4.5]{JS03}), $\{\mP_n\}_{n\geq 1}$ is  tightness.
Estimate \eqref{ACV6} follow by taking weak limits for \eqref{Kr0}.
\end{proof}
\bc\label{COR1}
There is a subsequence $n_k$ so that $\mP_{n_k}$ weakly converges to some $\mP\in\cP(\mD)$ and
\begin{align}\label{ABC7}
\lim_{k\to\infty}\|(\rho_{n_k}-\rho)\1_{[0,T]\times B_m}\|_{\mL^1_t\mL^1_x}=0
\quad \hbox{for every }  m\in\mN.
\end{align}
\ec

\begin{proof}
Note that $\rho_n$ solves the following Fokker-Planck equation:
$$
\p_t\rho_n=\Delta^{\alpha/2}\rho_n+\div(\wt b_n \rho_n).
$$
Let $\chi_m=\chi^0_m$ be the cutoff function as in \eqref{CHI}. By Sobolev's embedding, it is easy to see that
\begin{align*}
\|\p_t(\rho_n\chi_m)\|_{\mL^1_T(\bH^{-3}_2)}&\leq \|(\Delta^{\alpha/2}\rho_n)\chi_m\|_{\mL^1_T(\bH^{-3}_2)}+
\|\div(\wt b_n\chi_{2m} \rho_n)\chi_m\|_{\mL^1_T(\bH^{-3}_2)}\\
&\lesssim \|\Delta^{\alpha/2}\rho_n\|_{\mL^1_T(\bH^{-3}_2)}\|\chi_m\|_{\bH^3_2}+
\|\wt b_n\chi_{2m} \rho_n\|_{\mL^1_T\mL^1_x}\|\chi_m\|_{\bH^3_2}\\
&\lesssim \|\rho_n\|_{\mL^1_T\mL^1_x}+
\|\wt b_n\chi_{2m}\|_{\mL^{q_1}_T\mL^{\bbp_1}_x}\|\rho_n\|_{\mL^{q_2}_T\mL^{\bbp_2}_x},
\end{align*}
where $\frac1{q_1}+\frac1{q_2}=1$ and $\frac1{\bbp_1}+\frac1{\bbp_2}=\1$.
By \eqref{DL135} we have
$$
\sup_n\|\p_t(\rho_n\chi_m)\|_{\mL^1_T(\bH^{-3}_2)}<\infty.
$$
Since $\bH^\beta_p$ is locally compactly embedding in $\mL^{\bbp}_x$, by Aubin-Lions' lemma and a diagalization method, 
there is a subsequence $n_k$ so that \eqref{ABC7} holds for all $m\in\mN$. Moreover, by Prohorov's theorem, one can find a common subsequence $n_k$ 
such that  $\mP_{n_k}$ weakly converges to $\mP$ as $k\to\infty$.
\end{proof}

\bl
For any $m\in\mN$ and $q_0\in[1,q_1)$, $\bbp_0\in[\1,\bbp_1)$, we have
\begin{align}\label{ACV1}
\lim_{k\to\infty}\|(\wt b_{n_k}-\wt b)\1_{[0,T]\times B_m}\|_{\mL^{q_0}_t\mL^{\bbp_0}_x}=0,
\end{align}
where $\wt b(t,x):=\int_{\mR^d}b(t,x,y)\rho(t,y)\dif y$ and $n_k$ is the subsequence in Corollary \ref{COR1}.
\el
\begin{proof}
Note that by  \eqref{DR1} and the definition \eqref{GG1},
$$
\sup_k\|(\wt b_{n_k}-\wt b)\b1_{[0,T]\times B_m}\|_{\mL^{q_1}_t\mL^{\bbp_1}_x}\lesssim \nor h\nor_{\wt\mL^{q_1}_T\wt\mL^{\bbp_1}_x}.
$$
By the uniform integrability, to show \eqref{ACV1}, it suffices to prove that
\begin{align}\label{ACV3}
\lim_{k\to\infty}\|(\wt b_{n_k}-\wt b)\b1_{[0,T]\times B_m}\|_{\mL^1_t\mL^1_x}=0.
\end{align}
By definition we can write for any $R>0$,
\begin{align*}
(\wt b_{n_k}-\wt b)(t,x)&=\int_{\mR^d}[b_{n_k}(t,x,y)\rho_{n_k}(t,y)-b(t,x,y)\rho(t,y)]\dif y\\
&=\Big( \int_{B^c_R}+\int_{B_R}\Big) \left( b_{n_k}(t,x,y)\rho_{n_k}(t,y)-b(t,x,y)\rho(t,y) \right) \dif y\\
&=:I^R_k(t,x)+J^R_k(t,x).
\end{align*}
For $I^R_k$, by \eqref{DR1} we have
$$
\|I^R_k\|_{\wt\mL^1_T\wt\mL^1_x}\leq\|h\|_{\wt\mL^1_T\wt\mL^1_x}\int_{B^c_R}[\rho_{n_k}(t,y)+\rho(t,y)]\dif y.
$$
By the tightness of $\rho_{n_k}$, we get
\begin{align}\label{AVC2}
\lim_{R\to\infty}\sup_k\|I^R_k\|_{\wt\mL^1_T\wt\mL^1_x}=0.
\end{align}
For $J^R_k$, noting that by \eqref{ABC7} and the property of convolution,
$$
\mbox{$b_{n_k}(t,x,y)\rho_{n_k}(t,y)$ converges to $b(t,x,y)\rho(t,y)$ in Lebesgue measure $\dif t\times\dif x\times\dif y$,}
$$
by the uniform integrability, we have for each $R>0$,
$$
\lim_{k\to\infty}\|J^R_k\b1_{[0,T]\times B_m}\|_{\mL^1_T\mL^1_x}
=\lim_{k\to\infty}\int^T_0\!\!\!\int_{B_m}\!\int_{B_R}|b_{n_k}(t,x,y)\rho_{n_k}(t,y)-b(t,x,y)\rho(t,y)|\dif y\dif x\dif t=0,
$$ 
which together with \eqref{AVC2} yields \eqref{ACV3}.
\end{proof}

Now we can give the
\begin{proof}[Proof of Theorem \ref{Th55}]
Let $\mP_0$ be the law of $L^{(\alpha)}_\cdot$ in $\mD$.
Set $\mQ_k:=\mP_{n_k}\times \mP_0$. By Corollary \ref{COR1},
$\mQ_k$ weakly converges to $\mQ=\mP\times\mP_0$ as $k\to\infty$. By Skorokhod's representation theorem, 
there is a probability space $(\Omega,\sF,\bP)$ and $\mD\times\mD$-valued processes $(\wt X^k,\wt L^k)$ and $(\wt X,\wt L)$ such that
\begin{align}\label{ACV4}
(\wt X^k,\wt L^k)\to (\wt X,\wt L)\  \mbox{ in $\mD\times\mD$},\ \ \bP-a.s.,
\end{align}
and
$$
\bP\circ(\wt X^k,\wt L^k)^{-1}=\mQ_k,
\quad \bP\circ(\wt X,\wt L)^{-1}=\mQ.
$$
Moreover, $\wt L^k$ and $\wt L$ are still $\alpha$-stable L\'evy processes, and
$$
\wt X_t^k=\wt X_0^k+\int^t_0\wt b_{n_k}(s,\wt X^k_s)\dif s+\wt L^k_t.
$$
By \eqref{ACV4}, \eqref{ACV1} and taking limits, one sees that 
$$
\wt X_t=\wt X_0+\int^t_0\wt b(s,\wt X_s)\dif s+\wt L_t.
$$
Indeed, for each $m\in\mN$, by \eqref{ACV4} and the dominated convergence theorem, we have
\begin{align}\label{CV4}
\lim_{k\to\infty}\int^t_0|\wt b_{n_m}(s,\wt X^k_s)-\wt b_{n_m}(s,\wt X_s)|\dif s=0,\ \ \bP-a.s.
\end{align}
Since $\mI^0_\alpha$ is open, there is a $\theta>1$ so that 
$$
(q_0,\bbp_0):=(\tfrac{q_1}{\theta},\tfrac{\bbp_1}{\theta})\in\mI^0_\alpha.
$$
By H\"older's inequality and Krylov's estimate \eqref{Kr0} we have
\begin{align*}
&\bE\left[ \int^t_0|\wt b_{n_m}(s,\wt X^k_s)-\wt b(s,\wt X^k_s)|\b1_{|\wt X^k_s|\geq R}\dif s\right] \\
&\quad\leq\left(\bE\int^t_0|\wt b_{n_m}(s,\wt X^k_s)-\wt b(s,\wt X^k_s)|^\theta\dif s\right)^{1/\theta}\left(\bE\int^t_0\b1_{|\wt X^k_s|\geq R}\dif s\right)^{1-1/\theta}\\
&\quad\lesssim\left(\nor\wt b_{n_m}\nor^\theta_{\wt\mL^{q_1}_T\wt\mL^{\bbp_1}_x}+\nor\wt b\nor^\theta_{\wt\mL^{q_1}_T\wt\mL^{\bbp_1}_x}\right)^{1/\theta}
\left(\int^t_0\bP\{|\wt X^k_s|\geq R\}\dif s\right)^{1-1/\theta}\\
&\quad\lesssim\nor h\nor_{\wt\mL^{q_1}_T\wt\mL^{\bbp_1}_x}\left(\int^t_0\int_{B^c_R}\rho_{n_k}(s,y)\dif y\dif s\right)^{1-1/\theta}.
\end{align*}
By the tightness of $\rho_{n_k}$, we get
\begin{align}\label{CV5}
\lim_{R\to\infty}\sup_{m,k}\bE\left[ \int^t_0|\wt b_{n_m}(s,\wt X^k_s)-\wt b(s,\wt X^k_s)|\b1_{|\wt X^k_s|\geq R}\dif s\right] =0.
\end{align}
Moreover, for each $R>0$, by \eqref{Kr0} and \eqref{ACV1} we have
\begin{align*}
\sup_k\bE\left[ \int^t_0|\wt b_{n_m}(s,\wt X^k_s)-\wt b(s,\wt X^k_s)|\b1_{|\wt X^k_s|\leq R}\dif s\right] 
\lesssim \|(\wt b_{n_m}-\wt b)\1_{[0,T]\times B_R}\|_{\mL^{q_0}_t\mL^{\bbp_0}_x}\to 0
\end{align*}
as $m\to\infty$, which combining \eqref{CV4} and \eqref{CV5} yields 
$$
\int^t_0\wt b_{n_k}(s,\wt X^k_s)\dif s\to \int^t_0\wt b(s,\wt X_s)\dif s
$$
in probability as $k\to\infty$. The proof is complete.
\end{proof}
As a corollary, we have 
 the following existence of solutions to the
 fractional
  Navier-Stokes equations in $\mR^2$ with measure as initial vorticity.

\bc
Let $\alpha\in(1,2)$.
For any $\nu\in\cP(\mR^2)$, there exists a weak solution $\rho$ to equation \eqref{2VNS} in the sense that for each $\varphi\in C^\infty_c(\mR^d)$,
\begin{align*}
\int_{\mR^2}\rho(t,y)\varphi(y)\dif y
&=\int_{\mR^2}\varphi(y)\nu(\dif y)+\int^t_0\int_{\mR^2}\Delta^{\alpha/2}\varphi(y)\rho(s,y)\dif y\dif s\\
&\quad-\int^t_0\int_{\mR^2}\rho(s,y) (K_2*\rho(s))(y)\cdot\nabla\varphi(y)\dif y\dif s.
\end{align*}
Moreover, $\rho$ enjoys the regularity \eqref{DL131}.
\ec
\begin{proof}
Note that for any $p\in[1,2)$,
$$
\int_{B_1}|K_2(y)|^p\dif y<\infty.
$$
For $\alpha\in(1,2)$, one can choose $p\in(1,2)$ and $q$ close to infinity so that
$$
\tfrac{\alpha}q+\tfrac{2}{p}<\alpha.
$$
Thus we can use Theorem \ref{Th55} with $b(t,x,y)=K_2(x-y)$ to conclude the proof.
\end{proof}

\section{Well-posedness of generalized martingale problems in critical case}
In this section we show the well-posedness for a class of generalized martingale problems of SDE \eqref{SDE1}
in critical case. Fix $T>0$. Suppose that
\begin{align}\label{CN12}
\|\div b\|_{\mL^\infty_T\mL^\infty_x}+\nor b\nor_{\wt\mL^\infty_T\wt\mL^1_x}+\|b\|_{\rm crit}<\infty.
\end{align}
For given $f\in  C^\infty_c(\mR^{d})$, consider the following backward Kolmogorov equation:
\begin{align}\label{BPDE}
\p_t u+\Delta^{\alpha/2} u+b\cdot\nabla u=f \quad \hbox{with} \quad u(T)=0.
\end{align}
By reversing the time variable and Theorem \ref{Th42},
there is a unique weak solution $u_f$ to \eqref{BPDE} with
$$
u_f\in C_b([0,T]\times\mR^{d}).
$$
Here the continuity of $u_f$ is crucial
for the following notion of generalized martingale solutions.
\bd[Generalized martingale problem]\label{MP81}
Let $s\geq 0$ and $\nu\in\cP(\mR^{d})$. A probability measure $\mP\in\cP(\mD)$ is called a generalized martingale solution of SDE \eqref{SDE1} starting from $\nu$ 
at time $s$ if $\mP\circ \omega^{-1}_s=\nu$ and for any $T>s$ and $f\in C^\infty_c(\mR^{d})$, the process
\begin{align}\label{MM17}
M_t:=u_f(t,\omega_t)-u_f(s,\omega_s)-\int^t_sf(\omega_r)\dif r,\ \ t\in[s,T],
\end{align}
is a $\cB_t$-martingale with respect to $\mP$, where $u_f$ is the unique solution of \eqref{BPDE}.
The set of all the generalized martingale solutions $\mP\in\cP(\mD)$  with initial distribution $\nu$ at time $s$ is denoted by $\cM^{b}_{s,\nu}$.
\ed
\br\rm
The above notion of martingale solutions was originally introduced in \cite[Chapter 4]{EK86}.
It should be noticed that if \eqref{BPDE} has a $C^2$-solution $u$, then by It\^o's formula, any weak solution 
 of \eqref{SDE1} must be a generalized martingale solution. In general, 
these two notions are not equivalent.
\er

The following is the main result of this section. 

\bt
Under \eqref{CN12}, for each $s\geq 0$ and $\nu\in\cP(\mR^d)$, there is a unique generalized martingale solution 
$\mP\in\cM^{b}_{s,\nu}$ to SDE \eqref{SDE1} in the sense of Definition \ref{MP81}, where $\mP$ is also the distribution of a weak solution of SDE \eqref{SDE1} in the sense of Definition \ref{Def2}. 
Moreover, 
 for Lebesgue almost every $t\geq s$, $\mP\circ \omega_t^{-1}(\dif x)=\rho(t,x)\dif x$ and $\rho(t,x)$ solves the following Fokker-Planck equation
 in the distributional sense
 \begin{align}\label{FPK1}
 \p_t\rho=\Delta^{\alpha/2}\rho-\div(b\rho)
 \quad \hbox{with} \quad 
 \lim_{t\downarrow s}\rho(t,\cdot)=\nu\ \mbox{ \rm weakly},
\end{align}
 and $\|\rho\|_{\mL^{q}_T(\bH^{\beta}_{\bbp})}<\infty$ for any $\beta\in[0,\frac\alpha2)$ and $(q,\bbp)\in(1,\infty)^{1+d}$ with $\tfrac\alpha{q}+|\tfrac1{\bbp}|>d+\beta$.
\et
\begin{proof}
{\bf (Existence)} Without loss of generality we assume $s=0$.
Consider the approximating SDE
$$
X^n_t=X_0+\int^t_0b_n(s, X^n_s)\dif s+L^{(\alpha)}_t,
$$ 
where $\nu=\bP\circ X^{-1}_0$ and
$$
b_n(t,x):=[b(t,\cdot)*\varGamma_n](x).
$$
Let $\mP\in\cP(\mD)$ be an accumulation point of $(\mP_n)_{n\in\mN}$.
By 
Lemma \ref{Le56},
$$
(\mP_n\circ \omega_t^{-1})(\dif x)=\rho_n(t,x)\dif x,\ \ \mP\circ \omega_t^{-1}(\dif x)=\rho(t,x)\dif x.
$$
Since for each $\varphi\in C^\infty_c(\mR^{d})$, by It\^o's formula,
$$
\p_t \int \varphi\rho_n=\int \rho_n\Delta^{\alpha/2}\varphi+\int (b_n \cdot\nabla \varphi) \rho_n,
$$
taking limits along the subsequence $n_k$ in Corollary \ref{COR1} and by \eqref{ABC7} and \eqref{ACV1}, it is easy to see that
$$
\p_t \int \varphi\rho=\int \rho\Delta^{\alpha/2}\varphi+\int (b\cdot\nabla \varphi) \rho.
$$
In other words, $\rho$ solves \eqref{FPK1} in the distributional sense.
Fix $T>0$ and $f\in C^\infty_c(\mR^{d})$. Let $u_n\in L^\infty_T(C^\infty_b(\mR^{d}))$ be the unique smooth solution of the following backward PDE:
\begin{align}\label{CC10}
\p_t u_n+ \Delta^{\alpha/2}u_n+b_n\cdot\nabla u_n=f,\ \ u_n(T)=0.
\end{align}
By Theorem \ref{Th42}, we have for any $m\in\mN$,
\begin{align}\label{CC1}
\lim_{n\to\infty}\sup_{(t,x)\in [0,T]\times B_m}|u_{n}(t,x)-u_f(t,x)|=0,
\end{align}
where $u_f$ is the unique weak solution of PDE \eqref{BPDE}.
We now show that for each $T\geq t_1>t_0\geq 0$ and bounded continuous $\cB_{t_0}$-measurable functional $G_{t_0}$,
\begin{align}\label{MM35}
\mE^{\mP}(M_{t_1} G_{t_0})=\mE^{\mP}(M_{t_0} G_{t_0}),
\end{align}
where $M_t$ is defined by \eqref{MM17}.
By SDE \eqref{SDE19}, It\^o's formula and \eqref{CC10}, it is easy to see that
\begin{align}\label{MM45}
\mE^{\mP_n}(M^n_{t_1} G_{t_0})=\mE^{\mP_n}(M^n_{t_0} G_{t_0}),
\end{align}
where	
$$
M^n_t:=u_n(t,\omega_t)-u_n(0,\omega_0)-\int^t_0f(\omega_r)\dif r,
$$
Sending $n\to \infty$ in  \eqref{MM45}, 
 and using the pointwise convergence \eqref{CC1}, we obtain \eqref{MM35}.
Moreover, as in the proof of Theorem \ref{Th55}, $\mP$ is also the distribution of a weak solution of SDE \eqref{SDE1} in the sense of Definition \ref{Def2}. 

\medskip

{\bf (Uniqueness)} The uniqueness is a direct consequence of the definition. In fact,
let $\mP_1,\mP_2\in\sM^{b}_{s,\nu}$ be two solutions of the generalized martingale problem.  Fix $T>s$ and $f\in C^\infty_c(\mR^{d})$. 
Let $u_f$ be the unique weak solution of (see Theorem  \ref{Th42}),
$$
\p_t u_f+\Delta^{\alpha/2} u_f+b\cdot\nabla u_f=f,\ \ u_f(T)=0.
$$
By Definition \ref{MP81} and $u_f(T)=0$, we have
$$
\int_{\mR^{d}}u_f(s,x)\nu(\dif x)=-\mE^{\mP_i}\int_s^Tf(\omega_r)\dif r,\quad i=1,2,
$$
which means that for each $T>s$,
\begin{align*}
\int_s^T\mE^{\mP_1} f(\omega_r)\dif r=\int_s^T\mE^{\mP_2} f(\omega_r)\dif r.
\end{align*}
Hence, for any $f\in C^\infty_c(\mR^{d})$
\begin{align}\label{SD91}
\mE^{\mP_1} f(\omega_T)=\mE^{\mP_2} f(\omega_T),\ \ \forall T>s.
\end{align}
From this, by a standard way (see Theorem 4.4.2 in \cite{EK86}), we derive that
$$
\mP_1=\mP_2.
$$
Indeed, it suffices to prove the following claim by induction: 

{\bf (C$_n$)} for given $n\in\mN$, and for any $s\leq t_1<t_2<t_n<T$ and strictly positive and bounded measurable functions $g_1,\cdots, g_n$ on $\mR^{d}$,
\begin{align}\label{GK1}
\mE^{\mP_1}(g_1(\omega_{t_1})\cdots g_n(\omega_{t_n}))=\mE^{\mP_2}(g_1(\omega_{t_1})\cdots g_n(\omega_{t_n})).
\end{align}
By \eqref{SD91}, {\bf (C$_1$)} holds.
Suppose now that {\bf (C$_n$)} holds for some $n\geq 2$. For simplicity we write
$$
\eta:=g_1(\omega_{t_1})\cdots g_n(\omega_{t_n})>0,
$$
and for $i=1,2$, we define new probability measures
$$
\dif\widetilde\mP_i:=\eta\dif\mP_i/(\mE^{\mP_i}\eta)\in\cP(\mD),\ \ \widetilde\nu_i:=\widetilde\mP_i\circ \omega^{-1}_{t_n}\in\cP(\mR^d).
$$
Now we show
$$
\widetilde\mP_i\in\cM^b_{t_n;\widetilde\nu_i},\ \ i=1,2.
$$
For any $f\in C^\infty_c(\mR^{d})$, let 
$$
M_t:=u_f(t,\omega_t)-u_f(t_n,\omega_{t_n})-\int^t_{t_n}f(\omega_r)\dif r,\ \ t\in[t_n,T].
$$
We only need to prove that for any $T\geq t'>t\geq t_n$ and bounded $\cB_t$-measurable $\xi$,
$$
\mE^{\widetilde\mP_i}\left(M_{t'}\xi\right)=\mE^{\widetilde\mP_i}\left(M_t\xi\right)\Leftrightarrow \mE^{\mP_i}(M_{t'}\xi\eta)=\mE^{\mP_i}(M_t\xi\eta),
$$
which follows from $\mP_i\in\cM^b_{s,\nu}$, $i=1,2$. Thus, by induction hypothesis and \eqref{SD91},
$$
\widetilde\nu_1=\widetilde\nu_2\Rightarrow\widetilde\mP_1\circ \omega^{-1}_{t_{n+1}}=\widetilde\mP_2\circ \omega^{-1}_{t_{n+1}},\ \ \forall T\geq t_{n+1}>t_n.
$$
which in turn implies that {\bf (C$_{n+1}$)} holds.
\end{proof}

\end{document}